\documentclass[10pt,letterpaper]{article}
\usepackage{amsmath,amsthm,amsfonts,color,hyperref,anysize,amssymb,enumitem,color,graphicx,epstopdf,mathtools}
\usepackage[section]{placeins}
\usepackage[usenames,dvipsnames]{xcolor}
\usepackage[normalem]{ulem}

\usepackage[margin=0.8in]{geometry}

\DeclareMathOperator*{\argmax}{arg\,max}
\DeclareMathOperator*{\res}{res}

\newcommand{\floor}[1]{\left\lfloor #1 \right\rfloor }
\newcommand{\ceil}[1]{\left\lceil #1 \right\rceil }

\newcommand{\ep}{\epsilon}
\newcommand{\rr}{\mathbb{R}}
\newcommand{\rrplus}{\mathbb{R}^+}

\newcommand{\ra}{\rightarrow}

\newcommand{\pc}{p^\circ}
\newcommand{\xc}{x^\circ}
\newcommand{\ps}{p^*}
\newcommand{\tg}{\tilde{g}}

\newcommand{\tp}{\tilde{p}}
\newcommand{\tx}{\tilde{x}}
\newcommand{\tP}{\tilde{P}}
\newcommand{\comm}[1]{\qquad\mbox{(#1)}}
\newcommand{\xs}{x^*}

\newcommand{\G}{\Gamma}
\newcommand{\pt}{x^t}
\newcommand{\ptone}{x^{t-1}}

\newcommand{\pj}{p_j}
\newcommand{\Dp}{\Delta x}
\newcommand{\Dx}{\Delta x}

\newcommand{\xminus}{x^{-}}
\newcommand{\xplus}{x^{+}}
\newcommand{\pk}{p_k}
\newcommand{\Lj}{L_{j}}
\newcommand{\Ljj}{L_{jj}}

\newcommand{\Ljk}{L_{jk}}
\newcommand{\Dt}{\Delta t}
\newcommand{\td}{\tilde{d}}
\newcommand{\PRG}{\text{\textsf{PRG}}}
\newcommand{\PRGe}{\text{\emph{\textsf{PRG}}}}
\newcommand{\hW}{\widehat{W}}
\newcommand{\hWj}{\widehat{W}_j}
\newcommand{\hWk}{\widehat{W}_k}
\newcommand{\hd}{\widehat{d}}
\newcommand{\hdj}{\widehat{d}_j}

\newcommand{\Wj}{W_j}

\newcommand{\ve}{\vec{e}}
\newcommand{\kt}{k_t}
\newcommand{\ktau}{k_\tau}
\newcommand{\bG}{\overline{\G}}
\newcommand{\calR}{\mathcal{R}}

\newcommand{\Gjt}{\Gamma_j^t}

\newcommand{\Lmax}{L_{\max}}

\newcommand{\Lres}{L_{\res}}

\newcommand{\hWkt}{\widehat{W}_{k_t}}
\newcommand{\hWktau}{\widehat{W}_{k_{\tau}}}
\newcommand{\Psij}{\Psi_j}
\newcommand{\Psik}{\Psi_k}
\newcommand{\Psikt}{\Psi_{k_t}}

\newcommand{\gej}{g_j}
\newcommand{\gejp}{g'_j}
\newcommand{\gkt}{g_{k_t}^t}

\newcommand{\gkst}{g_{k}^t}

\newcommand{\gktt}{g_{k_t,2}^t}

\newcommand{\pks}{x_{k_s}}
\newcommand{\pkt}{x_{k_t}}

\newcommand{\xktautau}{x_{k_{\tau}}^{\tau}}
\newcommand{\pktp}{x_{k'_t}}

\newcommand{\pktt}{x_{k_t}^t}
\newcommand{\pktm}{x_{k}^{t-1}}
\newcommand{\pkttm}{x_{k_t}^{t-1}}
\newcommand{\pkttau}{x_{k_t}^{\tau}}
\newcommand{\tgkt}{\tilde{g}_{k_t}^t}

\newcommand{\gjo}{g_j^1}
\newcommand{\gjtw}{g_j^2}

\newcommand{\kmax}{\kappa_{\max}}

\newcommand{\muF}{\mu_F}
\newcommand{\muf}{\mu_f}

\newcommand{\Xs}{X^*}

\newcommand{\expect}[1]{\mathbb{E}\left[ #1 \right]}
\newcommand{\expectpi}[1]{\mathbb{E}_{\pi}\left[ #1 \right]}

\renewcommand{\O}{\mathcal{O}}

\def\pjq{p_{j,q}}
\def\tzjq{\tilde z_{j,q}}
\def\D{\Delta}

\newtheorem{theorem}{Theorem}
\newtheorem{lemma}{Lemma}

\newtheorem{cor}[lemma]{Corollary}

\newtheorem{assume}{Assumption}
\newtheorem{defn}{Definition}

\newtheorem{notation}{Notation}

\newenvironment{rlemma}[2][Lemma]{%
\medskip
\noindent\textbf{#1 \ref{#2}.} \begin{itshape}}{\end{itshape} 
}

\newenvironment{pf}{\begin{proof}[\emph{\textbf{Proof: }}]}{\end{proof}}
\newenvironment{pfof}[1]{\begin{proof}[\emph{\textbf{Proof of #1: }}]}{\end{proof}}

\newcommand{\vuppp}{\vspace{-.06in}}

\title{
A Unified Approach to Analyzing Asynchronous Coordinate Descent and Tatonnement}
\author{Yun Kuen Cheung\thanks{Part of the work done while this author was at
		the Courant Institute, NYU and at the Faculty of Computer Science, University of Vienna.
		He was supported in part by NSF Grant CCF-1217989, the Vienna Science and Technology Fund (WWTF) project ICT10-002,
		and the European Research Council under the European Union's Seventh Framework Programme
		(FP7/2007-2013) / ERC Grant Agreement no.~340506.}\\
		Max Planck Institute for Informatics\\Saarland Informatics Campus 
\and
		Richard Cole\thanks{This work was supported in part by NSF Grants CCF-1217989 and CCF-1527568.}\\
		Courant Institute, NYU
}
\date{}

\begin{document}\setlength{\parindent}{0.2in}
\maketitle

\begin{abstract}\thispagestyle{empty}
This paper concerns asynchrony in iterative processes, focusing on gradient descent and tatonnement, a fundamental price dynamic.

Gradient descent is an important class of iterative algorithms for minimizing convex functions.
Classically, gradient descent has been a sequential and synchronous process,
although distributed and asynchronous variants have been studied since the 1980s.
Coordinate descent is a commonly studied version of gradient descent.
In this paper, we focus on asynchronous coordinate descent on convex functions $F:\rr^n\ra\rr$ of the form
$$F(x) = f(x) + \sum_{k=1}^n \Psi_k(x_k),$$
where $f:\rr^n\ra\rr$ is a smooth convex function, and each $\Psi_k:\rr\ra\rr$ is a univariate and possibly non-smooth convex function.
Such functions occur in many data analysis and machine learning problems.

We give new analyses of cyclic coordinate descent, a parallel asynchronous stochastic coordinate descent,
and a rather general worst-case parallel asynchronous coordinate descent. 
For all of these, we either obtain sharply improved bounds, or provide the first analyses.
Our analyses all use a common amortized framework.
The application of this framework to the asynchronous stochastic version requires some new ideas,
for it is not obvious how to ensure a uniform distribution
where it is needed in the face of asynchronous actions that may undo uniformity.
We believe that our approach may well be applicable to the analysis of other iterative asynchronous stochastic processes.

We extend the framework to show that an asynchronous version of tatonnement,
a fundamental price dynamic widely studied in general equilibrium theory,
converges toward a market equilibrium for Fisher markets with CES utilities or Leontief utilities,
for which tatonnement is equivalent to coordinate descent.
\end{abstract}

\newpage
\thispagestyle{plain}\setcounter{page}{1}

\section{Introduction}
Gradient descent, an important class of iterative algorithms for minimizing convex functions, is a key subroutine in many computational problems.
Broadly speaking, gradient descent proceeds by moving iteratively in the direction of the negative gradient of a convex function.
Classically, gradient descent is a sequential and synchronous process.
Distributed and asynchronous variants have also been studied, starting with the work of Tsitsiklis et al.~\cite{TBA1986} in the 1980s;
these variants have been experiencing a resurgence of attention due to recent advances in multi-core parallel processing technology
and a strong demand for speeding-up large-scale gradient descent problems via parallelism.

In this paper, we consider the problem of (approximately) finding a minimum point of a convex function $F:\rr^n\ra\rr$ of the following form:
$$F(x) = f(x) + \sum_{k=1}^n \Psi_k(x_k),$$
where $f:\rr^n\ra\rr$ is a smooth convex function, and each $\Psi_k:\rr\ra\rr$ is a univariate convex function, but may be non-smooth.
Such functions occur in many data analysis and machine learning problems, such as linear regression
(e.g., the Lasso approach to regularized least squares~\cite{Tib94}) where $\Psi_k(x_k) = |x_k|$,
logistic regression~\cite{MVdGB2008}, ridge regression~\cite{SGV1998} where $\Psi_k(x_k)$ is a quadratic function,
and Support Vector Machines~\cite{CV1995} where $\Psi_k(x_k)$ is often a quadratic function or a hinge loss (essentially, $\max\{0,x_k\}$).
Neither $|x_k|$ nor the hinge loss are smooth at $x_k=0$.

Coordinate descent is a commonly studied version of gradient descent.
It proceeds by repeatedly selecting and updating a single coordinate of the argument to the convex function.
Coordinate descent is appealing because in a single iteration one need compute the gradient only in the direction of the update,
which can be much less expensive than the usual gradient descent, which computes the full gradient.
The offsetting cost, of course, is that in any iteration, little progress may be made.

Until fairly recently, convergence could be shown only for versions of coordinate descent in which
each update was along the ``best'' coordinate, whose determination required the computation of the full gradient,
rather nullifying the advantage of this approach~\cite{Nesterov2012}.
Futhermore, the presence of the non-smooth terms in $F(x)$ has meant that
until recently there were few analyses of the convergence rate of coordinate descent on such functions.
Of late, the advent of very large scale problems in which the full gradient of the smooth part might be unavailable or too expensive to compute,
has resulted in the use of coordinate descent, and the observation that it works well in practice.
Analyses have emerged to justify this observation.

The first analyses were for sequential stochastic updates~\cite{Nesterov2012,RiT2014, LX2013, LWRBS2015}:
at each iteration the next coordinate to update is chosen uniformly at random
(there are also versions in which different coordinates can be selected with different probabilities).
The expected bounds in this case provide a benchmark against which the bounds for other update sequences can be compared.

The next set of analyses were for cyclic coordinate descent~\cite{Tseng2001, ST2013, BT2013, HuaY2015}.
Here the coordinates are repeatedly updated one at a time in a fixed order.
We improve the best currently known worst case bound on the rate of convergence due to Hua and Yamashita~\cite{HuaY2015}
by a factor of $\Theta(\sqrt n/ \log n)$ where $n$ is the number of coordinates.

Of late, there has also been considerably interest in parallel versions of coordinate descent.
One important issue in parallel implementations is whether to ensure the different processors
are all using up-to-date information for their computations.
To ensure this requires considerable synchronization, locking, and consequent waiting.
The two analyses to date that do not make this assumption are those
by Avron, Druinsky and Gupta~\cite{ADG2014}, and Liu and Wright~\cite{LiuW2015};
they called this the ``inconsistent read'' model. We follow this approach;
its advantage is that it reduces and potentially eliminates the need for waiting.
At the same time, as some of the data being used in calculating updates will be out of date,
one has to ensure that the out-of-datededness is bounded in some fashion.
Again, following Liu and Wright, we assume there is a bounded amount of overlap between the various updates,
and we show how to implement this property in a lightweight manner.
By means of a new analysis we considerably improve the bounds given by Liu and Wright, 
achieving linear speed-up for a considerably wider range of parameters.
Further, as we will explain, Liu and Wright's analysis is in fact incomplete and it is not clear how to complete it.

We also consider worst case bounds in the asynchronous parallel setting.
This analysis covers the following natural case, among many others:
suppose the coordinates are evenly partitioned among the processors, with each processor updating its coordinates in cyclic order.
Again, for a considerable parameter range, we achieve linear speedup compared to the sequential cyclic coordinate case.
As far as we know, there was no previous analysis for this type of setting.

The challenges produced by an asynchronous setting are the lack of control over the timing of updates
(as we shall see, even in a stochastic setting one cannot assume the update commits are uniformly distributed across the coordinates),
and the possibility of having out-of-date information, as already mentioned.
This would mean that bad individual updates (updates that increase the value of the convex function) are unavoidable in general.
As we shall see, in our Asyncrhonous Coordinate Descent (ACD) algorithms,
every update leads to errors in subsequent gradient measurements at other cores (or processors).
A natural question to ask is whether such errors can propagate and be persistent and whether they might,
in the worst case, prohibit convergence toward a minimal point.
Our amortized analysis shows that this will not happen when the step sizes used in the ACD algorithms and the degree of parallelism are suitably bounded.
The following observation forms a key part of the analysis: if there is a bad update to one coordinate,
it can only be due to some recent good updates to other coordinates, or to \emph{chaining} of this effect.

Our analyses all follow a common framework.
We use an amortized analysis which relates the actual progress to the desired progress,
where the desired progress is a constant fraction of the progress achieved with sequential stochastic updating.
The amortization is used to hide the difference between these two measures of progress by amortizing it over multiple updates.
As we shall see, this difference is bounded by the squares of appropriate gradient differences,
and using Lipschitz parameters, these can in turn be bounded by changes to the squares of recent changes to the coordinates.
The final ingredient is to show that the progress is an upper bound on a  (small) multiple of the squares of the change to the updated coordinate.
Combining these ingredients yields a lower bound on the rate of progress.
Of course, the details of the combination vary from one version of coordinate descent to another,
and the asynchronous stochastic analysis in particular is quite non-trivial.

We finish by extending this framework to analyze an asynchronous tatonnement price update rule for some classes of market economies.
The reason for considering an asynchronous rule is that in the context of the Ongoing Market introduced by Cole and Fleischer~\cite{CF2008}
it allows for a more natural behavior than the classic auctioneer model of tatonnement, as we discuss in more detail in Section \ref{sect:CES}.

\paragraph{Related Work}

Convex optimization is one of the most widely used methodologies in applications across multiple disciplines.
Unsurprisingly, there is a vast literature studying convex optimization, with various assumptions and in various contexts.
We refer readers to Nesterov's text~\cite{Nesterov2004} for an excellent overview of the development of optimization theory.
Distributed and asynchronous computation has a long history in optimization,
initiated by the seminal works of Tsitsiklis, Bertsekas and Athans \cite{TBA1986,BT1989};
more recent results include~\cite{Borkar1998,BT2000}. See Frommer and Szyld \cite{FS2000} for a fairly recent review,
Liu et al.~\cite{LWRBS2015} for an account of recent developments, and Wright~\cite{Wright2015} for a recent survey on coordinate descent.

While many convergence results were proved for synchronous gradient descent,
only a few convergence results are known for cyclic coordinate descent,
including Tseng~\cite{Tseng2001}, Saha and Tewari~\cite{ST2013}, Beck and Tetruashvili~\cite{BT2013} and Hua and Yamashita~\cite{HuaY2015}.
We note that before the year 2012, as pointed out by Nesterov~\cite{Nesterov2012},
no global convergence rate for cyclic coordinate descent had been established,
and he commented that establishing such a bound is ``almost impossible''.
Although the recent work~\cite{ST2013,BT2013,HuaY2015} proved his comment wrong,
the late appearance of analyses for such a natural coordinate descent methodology suggests it is difficult to analyze,
and further serves as an indicator of the strength of our amortized approach which handles a wide range of scenarios.

Stochastic coordinate descent, in which coordinates are updated in random order,
has recently attracted attention. Relevant works include Nesterov \cite{Nesterov2012},
Richt{\'{a}}rik and Tak{\'{a}}c \cite{RiT2014} and Lu and Xiao \cite{LX2013}.

Two versions of asynchronous stochastic coordinate descent were analyzed by Liu et al.~\cite{LWRBS2015} and by Liu and Wright~\cite{LiuW2015}.
Both obtained bounds for both convex and ``optimally'' strongly convex functions\footnote{This is a weakening of the standard strong convexity.},
attaining speed-up more or less linear in the number of cores so long as they are not too numerous.
Liu et al.~\cite{LWRBS2015} obtained bounds similar to ours (see their Corollary 2 and our Theorem \ref{thm:main-SACD}),
but the version they analyzed is more restricted than ours in two aspects:
first, they imposed the unrealistic assumption of consistent reads, and second,
they considered only smooth functions (i.e., without the non-smooth univariate components $\Psi_k$).
The version analyzed by Liu and Wright~\cite{LiuW2015} is the same as ours.
While not explicitly stated, the latter bound degrades when the parallelism exceeds $\Theta(n^{1/4})$.\footnote{This is
expressed in terms of a parameter $\tau$, renamed $q$ in this paper, which is essentially the possible parallelism;
this connection depends on the relative times to calculate different updates. Their bound also depends on the ratio of two Lipschitz parameters.}
Our bound has a similar flavor but with a limit of $\Theta(n^{1/2})$.
Further, both analyses overlooked the fact that asynchronous actions may undo the distributional uniformity of the updates
(we will explain why this happens in Section~\ref{sect:stochastic-ACD}), and thus their analyses are incomplete.

In statistical machine learning, the objective functions to be minimized typically have the form
$\sum_i \ell(x,Z_i)$, where the $Z_i$ are samples; updates are typically sample-wise but not coordinate-wise,
so our model will not cover these update algorithms.
It is a very interesting problem to investigate if there is an adaption of our model and amortized analysis for these algorithms.
Many of the assumptions made for asynchronous sample-wise updates share similarities with ours.
For instance, Tsianos and Rabbat \cite{TR2012} extended the analysis of Duchi, Agarwal and Wainwright~\cite{DAW2012}
to analyse distributed dual averaging (DDA) with communication delay;
the same authors \cite{RT2014} studied DDA with heterogeneous systems, i.e., distributed computing units with different query and computing speeds.
Langford, Smola and Zinkevich~\cite{LSZ2009} and Liu et al.~\cite{LWRBS2015} also studied problems with bounded communication delay.

Avron, Druinsky and Gupta~\cite{ADG2014} recently studied an asynchronous and randomized version of the Gauss-Seidel algorithm
for solving symmetric and positive definite matrix systems.
They show that being less aggressive, i.e., reducing step sizes, can actually improve the guaranteed convergence rate.

In a similar spirit to our analysis, Cheung, Cole and Rastogi~\cite{CCR2012} analyzed asynchronous tatonnement in certain Fisher markets.
This earlier work employed a potential function which drops continuously when there is no update and does not increase when an update is made.
This approach could be followed for the current market setting, and in fact is the approach we took in an earlier version of this work
(this did require the design of a completely new potential function compared to the one used in~\cite{CCR2012}).
But in the current work, we instead use a discrete analysis which has more in common with our coordinate descent analyses.

For the closely related topic of learning dynamics in games, where updates are based on the payoffs received by agents,
again, the classical approach assumes synchronous or round-robin updates with up-to-date payoffs;
models with stochastic update schedules were also studied previously (e.g., in~\cite{Blume1993,AN2010,MS2012}),
while learning dynamics with delayed payoffs~\cite{PalLa2015} were studied recently.

\paragraph{Organization of This Paper}

In Section \ref{sect:model}, we formally describe our models of  coordinate descent,
and state our main results on cyclic coordinate descent, asynchronous coordinate descent and its stochastic variant.
In Section \ref{sect:keyidea}, we highlight the key observations and lemmas
which underlie our three main results.
In Sections \ref{sect:cyclic-CD}--\ref{sect:parallel-ACD}, we delve further into each of these results.
Finally, in Section \ref{sect:CES}, we describe our model of asynchronous tatonnement and state the associated convergence result.
Omitted proofs are given in the appendix.

\section{Model and Main Coordinate Descent Results}\label{sect:model}

\newcommand{\Ap}{A^+}
\newcommand{\Am}{A^-}
\newcommand{\pktaum}{x_k^{\tau-1}}
\newcommand{\ptauo}{x^{\tau-1}}
\newcommand{\ptauprm}{x^{\tau+n-1}}
\newcommand{\ptp}{x^{t+1}}

We consider the problem of (approximately) finding a minimum point of a convex function $F:\rr^n\ra\rr$ of the form
$F(x) = f(x) + \sum_{k=1}^n \Psi_k(x_k)$, where $f:\rr^n\ra\rr$ is a smooth convex function,
and each $\Psi_k:\rr\ra\rr$ is a univariate and possibly non-smooth convex function.
Let $X^*$ denote the set of minimum points of $F$; we use $x^*$ to denote a minimum point of $F$.
Without loss of generality, we assume that $F^*$, the minimum value of $F$, is zero.

We recap a few standard terminologies. Let $\ve_j$ denote the unit vector along coordinate $j$.

\begin{defn}\label{def:Lipschitz-parameters}
The function $f$ is $L$-Lipschitz-smooth if for any $x,\Dx\in\rr^n$, $\|\nabla f(x+\Dx) - f(x)\| ~\le~ L\cdot\|\Dx\|$.

For any coordinates $j,k$, the function $f$ is $\Ljk$-Lipschitz-smooth if for any $x\in\rr^n$ and $r\in\rr$,
$|\nabla_k f(x+r\cdot \ve_j) - \nabla_k f(x)| ~\le~ \Ljk\cdot |r|$. Also, as is standard, $\Lj$ denotes $\Ljj$.
	
Let $\Lmax ~:=~ \max_{j,k} \Ljk$, and let $\Lres ~:=~ \max_j \left(\sum_k (\Ljk)^2\right)^{1/2}$.
\end{defn}

Next, we introduce the following notation.
\begin{notation}\label{not:nat-terms}
For $d,g,x\in\rr$, $\G\in\rrplus$ and $\Psi:\rr\ra\rr$ a univariate, proper, convex and lower semi-continuous function, let
$$
W(d,g,x,\G,\Psi) ~:=~ -gd ~-~ \frac{\G}{2}\cdot d^2 ~+~ \Psi(x) ~-~ \Psi(x+d);
$$
the function $W$ is often called a \emph{proximal function}. Also, let
$$
\hW(g,x,\G,\Psi) ~:=~ \max_{d\in\rr}~W(d,g,x,\G,\Psi)~~~~~~~~\text{and}~~~~~~~~\hd(g,x,\G,\Psi) ~:=~ \argmax_{d\in\rr}~W(d,g,x,\G,\Psi).
$$
We note that $\hW(g,x,\G,\Psi) \geq W(0,g,x,\G,\Psi) = 0$.

We will use $k_t$ to denote the coordinate being updated at time $t$.
We will also use $\Delta x$ to denote $\hd(g,x,\G,\Psi)$, and more specifically $\Delta \pkt^t$ when we want to focus on the
coordinate $\kt$ updated at time $t$.
\end{notation}

\paragraph{Update Rule} 

If the implementation is sequential, or parallel and synchronous, the standard coordinate descent update rule, applied to coordinate $j$,
first computes the \emph{accurate} gradient $g_j := \nabla_j f(\ptone)$,
and then performs the update given below with a suitable parameter $\G_j$:
$$
\pt_j \leftarrow \ptone_j + \hd(g_j,\ptone_j,\G_j,\Psi_j)
~~~~~~~~\text{and}~~~~~~~~\forall k\neq j,~\pt_k \leftarrow \ptone_k.
$$

However, in an asynchronous environment, an updating processor (or core) might possess outdated information $\tx$ instead of $\ptone$,
so the gradient the core computes will be $\tg_j := \nabla_j f(\tx)$, instead of the accurate value $\nabla_j f(\ptone)$.
Our update rule, which is naturally motivated by its synchronous counterpart, is
\begin{equation}\label{eq:update-rule}
\pt_j \leftarrow \ptone_j + \hd(\tg_j,\ptone_j,\G_j,\Psi_j)~=~\ptone_j+\Delta \pt_j~~~~~~~~\text{and}~~~~~~~~\forall k\neq j,~\pt_k \leftarrow \ptone_k.
\end{equation}

\subsection{Results}

The basis form of all our results is given by the following meta-theorem (or a slight variant of it for the stochastic case),
where we use a common step size $\G$ for all the coordinates. 

\begin{theorem}\label{thm:meta-new}
Let $\G$ be a sufficiently large step size for the update rule, and let $r,q$ be two fixed integer parameters.
Let $A(t)$ be a non-negative function with $A(0) = 0$, and let $H(t):=F(\pt) + A(t)$.
Suppose that
\begin{itemize}
\item for all $t\ge 1$, $H(t) \le H(t-1)$, i.e., $H(t)$ is a decreasing function of $t$;
\item there exists constants $\alpha,\beta > 0$ such that for any $t\ge 2r$,
$$
\sum_{i=t-2r+1}^t \left[ H(i-1) - H(i) \right]
~\ge~
\sum_{i=t-2r+1}^{t-r} \left[~\frac \alpha n ~ \sum_{k=1}^n ~ \hWk(\nabla_k f(x^{i-1}),x_k^{i-1},\G,\Psi_k) ~+~ \frac \beta q \cdot A(i-1)~\right].
$$
\end{itemize}

\noindent
\emph{(i)} Then, if $F$ is strongly convex with parameter $\muF$,\footnote{i.e.,
for all $x,y \in \rr^n$ and $F'(x)$ which is any subgradient of $F$ at $x$,
$F(y) \ge F(x) + \langle F'(x), y-x \rangle + \frac 12 \muF ||y-x||^2$.}
and $f$ has strongly convex parameter $\muf$,
$$
F(\pt) ~\le~ \left[ 1 - \min\left\{\frac {\alpha}{2n} \cdot \frac{\muF} {\muF + \G - \muf}~,~ \frac{\beta}{2q} \right\} \right]^{t-2r+1} \cdot F(\xc).
$$

\noindent
\emph{(ii)} Now suppose that $F$ is convex.
Let $R$ be the radius of the level set for $\xc$. Formally,
let $X = \{x \,|\, F(x) \le F(\xc)\}$;
then
$R = sup_{x\in X} \min_{\xs \in \Xs} ||x - \xs||$.
Then, for $t \ge r$,
$$
F(\pt) ~\le~ \frac{F(\pc)}{1 + \min \left\{ \frac{\beta}{4q~F(\pc)} ~,~ \frac{\alpha}{8n~F(\pc)} ~,~ \frac{1}{8\G \calR ^2} \right\}
\cdot H(2r-1) \cdot (t-2r+1)}.
$$
\end{theorem}

The existing form of this theorem does not involve the function $A$ and the parameter $r$ had been limited to $r=1$.
The definition of $A$ varies from analysis to analysis though they are somewhat similar. 
The theorem is proved in the appendix for the main part of the analysis concerns the proof of the premise of this theorem.

For the sequential stochastic version of coordinate descent, $\alpha=1$ and $\G \ge \Lmax$ suffices.
Clearly, the minimum possible value of $\G$ and the value of $\alpha$ control the convergence rate.
Accordingly, in Table~\ref{tbl:results}, we report our new results in terms of $\alpha$ and the minimum $\G$.

\begin{table}[tbh]
\centering 
\begin{tabular}{|l|c|c|} \hline
Algorithm & Known Result & New Result \\ \hline
cyclic  & $\G \ge \max\{M,L\}$; $M \le L \sqrt n$ & $\G \ge \frac{4}{\sqrt 3}L \ceil{\log_2 n}$ \\
          & $\alpha = \O(1)$ \footnotesize{ \cite[2015]{HuaY2015} }\footnotemark
          & $r=n$ and $\alpha = \frac 13$ ($\beta$ not needed) \\ \hline
Asynchronous stochastic (SACD) & $\G \ge 2  \Lmax$ & $\G \ge \Lmax$ \\
~~~~linear speedup if: & $q \le \left( \frac {\G \sqrt n } {4e \Lres} \right)^{1/2}$
                                                & $q \le \min\left\{ \frac {\G \sqrt{n-q} } {8 \sqrt{10} \Lres}, \frac{9}{100}n \right\}$ \\
& \footnotesize{ \cite[2015]{LiuW2015} }\footnotemark & $r=1$, $\alpha = \frac 12$, and
$\beta = \frac{\alpha \cdot q}{n}$\\ \hline
Worst case asynchronous\ (PACD) & & $\G \ge \frac {16}{\sqrt 3} L \sqrt{\kmax}\ceil{\log_2 r}$ \\
~~~~linear speedup if: & no prior work
 & $q \le \frac {\G\sqrt3} {8 \cdot \Lmax }$ \\ 
&& $r \ge n$, $\alpha = \frac {n} {3 r}$, and $\beta = \frac{1}{2}$ \\ \hline
\end{tabular}
\caption{\label{tbl:results}The Coordinate Descent Results as a Function of $\G$, $\alpha$ and $\beta$.}
\end{table}
\footnotetext[4]{This analysis was not expressed precisely.} 
\footnotetext[5]{This analysis was incomplete. Also, it used a quite different approach so it is not evident what corresponded to $\alpha$.}

\paragraph{Cyclic Coordinate Descent (CCD)}
Here the $n$ coordinates are repeatedly updated in a fixed order. We have the following result.

\begin{theorem}\label{thm:main-cyclic}
If $\G \ge \frac {4}{\sqrt 3} L \ceil{\log_2 n}$, 
then the bound of Theorem~\ref{thm:meta-new} holds with $ r=n$ and $\alpha = \frac 13$.
($\beta$ is not needed.)
\end{theorem}

\paragraph{Asynchronous Coordinate Descent --- General Setting.}
The coordinate descent process starts at an initial point $\xc = (\xc_1,\xc_2,\cdots,\xc_n)$.
Multiple cores then iteratively update the coordinate values.
We assume that at each time, there is exactly one coordinate value being updated.
In practice, since there will be little coordination between cores,
it is possible that multiple coordinate values are updated at the same \emph{moment};
but by using an arbitrary tie-breaking rule in our analysis, we can immediately extend the analysis to these scenarios.

In our algorithms, the cores use a shared memory which stores the coordinate values.
Each core iteratively performs the following tasks without global coordination:
(1) chooses a coordinate $k$ according to \emph{some} rule;
(2) retrieves coordinate values from the shared memory; let $\tx = (\tx_1,\tx_2,\cdots,\tx_n)$ denote the retrieved values;
(3) computes $\nabla_k f(\tx)$;
(4) requests a lock on the memory that stores the value of the $k$-th coordinate;
(5) updates the $k$-th coordinate using rule \eqref{eq:update-rule};~\footnote{
	An implementation detail: even if the core had retrieved the value of the $k$-th coordinate from the shared memory in Step (2),
	the core needs to retrieve it \emph{again} in Step (5),
	because it needs the most updated value when applying the update rule \eqref{eq:update-rule}.}
(6) releases the lock.
We note that if the rule in Step (1) enforces that each coordinate is always updated by exactly one core,
then Steps (4) and (6) can be omitted.

The retrieval times for Step (2) plus the gradient-computation time for Step (3) can be non-trivial,
and also in Step (4) a core might need to wait if the coordinate it wants to update is locked by another core.
Thus, during this period of time other coordinates are likely to be updated.
For each update, we call the period of time spent for performing the above tasks the \emph{interim interval} of the update.
We say that update $A$ \emph{interferes} with update $B$ if the finishing time of update $A$ lies in the interim interval of update $B$.

\paragraph{Parallel Asynchronous Coordinate Descent (PACD)}
The time of an update refers to the moment it finishes.
We assume that the rule in Step (1) satisfies the following two assumptions:

\begin{assume}\label{assume:PACD-q}
There exists a non-negative integer $q$ such that for any update at time $t$, the only updates that interfere with it are those at times $t-1,t-2,\cdots,t-q$.
\end{assume}

\begin{assume}\label{assume:PACD-frequent}
There exists a parameter $r$ such that each coordinate is updated at least once in any time interval of length $r$.
Also, in any time interval of length $r$, each coordinate is updated at most $\kmax$ times.
\end{assume}

We note that even without any global coordination, in many circumstances,
$n \le r =\O(n)$ and $\kmax$ would be a small constant.
In Appendix \ref{sect:appendix-parallel-ACD}, we provide a simple implementation method
with modest overhead for strictly enforcing the two assumptions.

\begin{theorem}\label{thm:main-PACD}
Suppose that Assumptions \ref{assume:PACD-q} and \ref{assume:PACD-frequent} hold.
Then, if $\G \ge \frac {16}{\sqrt 3} L \sqrt{\kmax}\ceil{\log_2 r}$, Theorem~\ref{thm:meta-new}
holds with $\alpha= \frac{n} {3r}$ and $\beta = \frac{1}{2}$, when $q \le \frac {\G\sqrt 3} {8 \cdot \Lmax}$.
\end{theorem}

\paragraph{Stochastic Asynchronous Coordinate Descent (SACD)}
Here each core repeatedly selects a coordinate to update uniformly at random.
A difficulty we face is that the resulting orderings of the updates by start time and by commit time need not be the same.
For our result we assume that the algorithm selects exactly $\bar t$ 
coordinates to update for some prespecified $\bar t$, and that these coordinates are all updated,
with the commit times constrained by the following assumption.

\begin{assume}\label{assume:SACD-q}
There exists a non-negative integer $q$ such that
for any update at time $t$, the only updates that interfere with it are those at times $t-1,t-2,\cdots,t-q$ and $t+1,t+2,\cdots,t+q$.
\end{assume}
This can be enforced in almost the same way as Assumption~\ref{assume:PACD-q}.

\begin{theorem}\label{thm:main-SACD}
Suppose that Assumption \ref{assume:SACD-q} holds.
If $\G \ge \max_j \Lj$, then in expectation, the bound of Theorem~\ref{thm:meta-new}
holds for $t = \bar t$, with $r=1$ and $\alpha = \frac 12$,
when $q \le \min\left\{ \frac {\G \sqrt{n -q} } {8 \sqrt {10} \Lres}, \frac {9}{100} n\right\}$, and $\beta = \frac{\alpha \cdot q}{n}$.
\end{theorem}

\section{Key Ideas}\label{sect:keyidea}

\newcommand{\xktau}{x_{k_\tau}}
\newcommand{\pktautm}{x_{k_{\tau}}^{t-1}}
\newcommand{\gktaut}{g_{k_{\tau}}^t}

\subsection{Technical Overview}

Let $\kt$ denote the index of the coordinate that is updated at time $t$, let $\pkt$ denote its updated value,
let $\gkt ~:=~ \nabla_{k_t} f(\ptone)$ be the value of the gradient along coordinate $\pkt$ computed at time $t$
using up-to-date values of the coordinates, and let $\tgkt$ be the actual value computed, which may use some out-of-date values.

As we shall see, an update guarantees the following progress
\begin{equation}\label{eqn:basic-prog}
F(\ptone) - F(\pt) ~\ge~ \hWkt(\tgkt,\pkttm,\G,\Psikt).
\end{equation} 

In the case of a sequential stochastic update this implies
\begin{equation}
\label{eqn:seq-stoc-prog}
\expect{ F(\ptone) - F(\pt) } ~\ge~ \frac 1n \sum_{k=1}^n \hWk(\gkst,\pktm,\G,\Psik),
\end{equation}  
from which one can deduce the convergence results in Theorem~\ref{thm:meta-new} in expectation.

\subsection{Technical Approach}

All our results build on the following lemmas.

\begin{lemma}\label{lem:F-prog-one}
$F(\ptone) - F(\pt) \ge \hWkt(\gkt,\pkttm,\G,\Psikt) - \frac {1}{\G}(\gkt - \tgkt)^2.$
\end{lemma}

\begin{lemma}\label{lem:F-prog-two}
$F(\ptone) - F(\pt) \ge \frac 14 \G \left( \Delta \pkt \right)^2 - \frac {1}{\G} (\gkt - \tgkt)^2$.
If $\gkt = \tgkt$, then 
$F(\ptone) - F(\pt) \ge \frac 12 \G \left( \Delta \pkt \right)^2 $.
\end{lemma}

Combining Lemmas~\ref{lem:F-prog-one} and~\ref{lem:F-prog-two} yields
\begin{align}
\label{eqn:F-prog}
F(\ptone) - F(\pt) & \ge \frac 12 \hWkt(\gkt,\pkttm,\G,\Psikt) + \frac 18 \G \left( \Delta \pkt \right)^2 - \frac {1}{\G} (\gkt - \tgkt)^2 \\
\label{eqn:F-prog-better}
\text{or}~~~~~~~~F(\ptone) - F(\pt) & \ge \frac 12 \hWkt(\gkt,\pkttm,\G,\Psikt) + \frac 14 \G \left( \Delta \pkt \right)^2
~~~~~~~~~~\mbox{if $\gkt = \tgkt$}.
\end{align}
However, the progress we would like is the term $\frac 1n \sum_{k=1}^n  \hWk(\gkst,\pktm,\G,\Psik)$ from \eqref{eqn:seq-stoc-prog}.
The following lemma will allow us to relate the progress we have in \eqref{eqn:F-prog} to the progress we desire.

\begin{lemma}
\label{lem:W-shift}
For any $\gej$, $\gejp$,
$
\hWj(\gej,x_j,\G,\Psij) ~\ge~ \frac 23 \cdot \hWj(\gejp,x_j,\G,\Psij) ~-~ \frac{4}{3\G} \cdot (\gej - \gejp)^2.
$
\end{lemma}

The final issue will be to bound the differences $(\gjo - \gjtw)^2$ in terms of the $(\Delta x_{k_t}^t)^2$,
which we approach by using the Lipschitz parameters for the gradients as defined in Definition \ref{def:Lipschitz-parameters}.
Specifically, for $x^1,x^2\in \rr^n$,
for any $k$, let $\Dx_k ~:=~ x^1_k - x^2_k$, and
for $i=1,2$, let $g_j^i ~:=~ \nabla_j f(x^i)$. Then
\begin{align}
\label{eqn:g-diff-multi-L-bound}
\left( g_j^1 - g_j^2 \right)^2 &\le~ \left[\sum_{k=1}^n L_{kj} \left|\Dx_k\right| \right]^2\\
\label{eqn:g-diff-plain-L-bound}
\text{and}~~~~~~~\sum_{j=1}^n\left(g_j^1 - g_j^2\right)^2 & \le~ L^2 \sum_{k=1}^n \left(\Dx_k\right)^2.
\end{align}

The challenge is to combine these equations while minimizing $\G$. 

\section{Cyclic Coordinate Descent (CCD)}\label{sect:cyclic-CD}

\newcommand{\gktaunu}{g_{k_{\tau}}^{\nu}}
\newcommand{\gkttau}{g_{k_t}^{\tau}}
\newcommand{\gktautau}{g_{k_{\tau}}^{\tau}}
\newcommand{\gktausig}{g_{k_{\tau}}^{\sigma}}
\newcommand{\gktsig}{g_{k_{t}}^{\sigma}}
\newcommand{\pktautau}{x_{k_{\tau}}^{\tau}}
\newcommand{\pktautaum}{x_{k_{\tau}}^{\tau-1}}
\newcommand{\pkttaum}{x_{k_{\tau}}^{\tau-1}}
\newcommand{\pktausigm}{x_{k_{\tau}}^{\sigma-1}}
\newcommand{\pktsigm}{x_{k_{t}}^{\sigma-1}}
\newcommand{\pktaunum}{x_{k_{\tau}}^{\nu-1}}

In the CCD case, accurate gradients are used for all updates, so we will use equation \eqref{eqn:F-prog-better}.
To use Theorem \ref{thm:meta-new}, for the CCD case,
we set $r=n$ and $A(t) \equiv 0$, and hence $H(t) = F(\pt)$.
Then by \eqref{eqn:F-prog-better}, the first condition required in Theorem \ref{thm:meta-new} is satisfied.

Next, we establish the second condition required for using Theorem \ref{thm:meta-new}.
For any $t\ge 2n$,
$$
\sum_{i=t-2n+1}^t \left[ H(i-1) - H(i) \right]
~\ge~ \sum_{i=t-2n+1}^t \left[ ~\frac 12 \hW_{k_i} ( g_{k_i}^i , x_{k_i}^{i-1},\G,\Psi_{k_i}) ~+~ \frac \G 4 (\Delta x_{k_i})^2~\right].
$$

Applying Lemma \ref{lem:W-shift} yields
\begin{align}
&\sum_{i=t-2n+1}^t \left[ H(i-1) - H(i) \right]\nonumber\\
&~\ge~ \sum_{i=t-2n+1}^t \left[ ~
\sum_{j=\max\{t-2n+1,i-n+1\}}^i \left[~\frac{1}{3n} \hW_{k_i} (g_{k_i}^j,x_{k_i}^{j-1},\G,\Psi_{k_i}) - \frac{2}{3\G n}  \left( g_{k_i}^j - g_{k_i}^i \right)^2 ~\right]
~+~ \frac \G 4 (\Delta x_{k_i})^2~\right]\nonumber\\
&~\ge~ \sum_{j=t-2n+1}^{t-n} \frac{1}{3n} ~ \sum_{k=1}^n ~ \hWk (g_k^j,x_k^{j-1},\G,\Psi_k) \nonumber\\
&~~~~~~~~~~~~~~~~~
~+~ \underbrace{\left(\frac \G 4\sum_{i=t-2n+1}^t (\Delta x_{k_i})^2 ~-~ \frac{2}{3\G n} \sum_{i=t-2n+1}^t ~\sum_{j=\max\{t-2n+1,i-n+1\}}^i \left( g_{k_i}^j - g_{k_i}^i \right)^2 \right)}_{Q}\label{eqn:target-cyclic-final}
\end{align}
Thus, once we can prove that $Q\ge 0$ for some sufficiently large $\G$,
we can apply Theorem \ref{thm:meta-new} with $\alpha = 1/3$ and $\beta = 1$.
(Since $A(t)\equiv 0$, we can use any positive constant $\beta$.)
The rest of this section is devoted to showing that 
$\G ~\ge~ \frac{4}{\sqrt{3}} L \ceil{\log_2 n}$ suffices.

For each $i=0,1,2,\cdots,\floor{\lg (n-1)} = \ceil{\lg n}-1$, we divide the time interval
$[t-2n+1,t]$ into consecutive intervals of length $2^i$.
We call such time intervals the \emph{doubling intervals},
and let $\mathcal{D}$ denote the collection of all such doubling intervals.
By Definition \ref{def:Lipschitz-parameters}, for each doubling interval $[\tau_b,\tau_c]$,
\begin{equation}\label{eqn:cyclic-inter1}
\sum_{k=1}^n \left( g_k^{\tau_b} - g_k^{\tau_c} \right)^2 ~\le~ L^2 \sum_{\tau=\tau_b}^{\tau_c-1} \left( \Dx_{k_\tau}^{\tau} \right)^2.
\end{equation}

On the other hand, any $\left(g_{k_i}^j - g_{k_i}^i\right)$, with $0 \le i-j < n$, can be  written as a telescoping summation
$
\left(g_{k_t}^j - g_{k_i}^i\right)  ~=~ \sum_{a=1}^{\ell} ~ \left(g_{k_i}^{\tau_{a-1}} - g_{k_i}^{\tau_a}\right)
$
for some $\ell = \ell(j,i) \le 2\ceil{\lg n}$, where each $[\tau_{a-1},\tau_a]$ is a doubling interval,
$\tau_0 = j$, and $\tau_\ell = i$.
Then, by the Cauchy-Schwarz inequality,
\begin{equation}\label{eqn:cyclic-inter2}
\left(g_{k_i}^j - g_{k_i}^i\right)^2 ~\le~ 2\ceil{\lg n} ~ \sum_{a=1}^{\ell} ~ \left(g_{k_i}^{\tau_{a-1}} - g_{k_i}^{\tau_a}\right)^2.\vuppp\vuppp
\end{equation}
\begin{align}
\text{Then,}~~~~~~~~~~~~~~~~~~~~~&\frac{2}{3\G n} \sum_{i=t-2n+1}^t ~\sum_{j=\max\{t-2n+1,i-n+1\}}^i \left( g_{k_i}^j - g_{k_i}^i \right)^2
~~~~~~~~~~~~~~~~~~~~~~~~~~~~~~~~~~~~~~~~~~~~~~~~~~~~~~~~~~~\nonumber\\
&~\le~ \frac{4\ceil{\lg n}}{3\G n} \sum_{i=t-2n+1}^t ~\sum_{j=\max\{t-2n+1,i-n+1\}}^i \sum_{a=1}^{\ell(j,i)}
 ~ \left(g_{k_i}^{\tau_{a-1}} - g_{k_i}^{\tau_a}\right)^2\comm{by Eqn.~\eqref{eqn:cyclic-inter2}}\nonumber\\
&~\stackrel{(*)}{\le}~ \frac{4\ceil{\lg n}}{3\G n} \sum_{[\tau_b,\tau_c] \in \mathcal{D}} ~ \sum_{k=1}^n ~n \cdot \left(g_k^{\tau_b} - g_k^{\tau_c}\right)^2\nonumber\\
&~\le~ \frac{4\ceil{\lg n}}{3\G} \sum_{[\tau_b,\tau_c] \in \mathcal{D}} ~ \sum_{a=\tau_b}^{\tau_c-1} L^2 (\Delta x_{k_a}^a)^2~~~~~~~~\comm{by Eqn.~\eqref{eqn:cyclic-inter1}}\nonumber\\
&~\stackrel{(**)}{\le}~ \frac{4\ceil{\lg n}L^2}{3\G} \sum_{i=t-2n+1}^t \ceil{\lg n} \left(\Delta x_{k_i}^i\right)^2
~=~ \frac{4\ceil{\lg n}^2 L^2}{3\G} \sum_{i=t-2n+1}^t \left(\Delta x_{k_i}^i\right)^2.\label{eqn:g-to-x-log-bound}
\end{align}
We explain why inequalities $(*)$ and $(**)$ hold here.
For $(*)$, observe that in the triple summation above, for any doubling interval $[\tau_b,\tau_c]$ and any $k$,
each $\left(g_k^{\tau_b} - g_k^{\tau_c}\right)^2$ term appears at most $n$ times.
For $(**)$, observe that in the double summation above, for any $i$,
each $\left(\Delta x_{k_i}^i\right)^2$ term appears at most $\ceil{\lg n}$ times.

Thus, for $Q\ge 0$, it suffices that $\frac{4\ceil{\lg n}^2 L^2}{3\G} ~\le~ \frac \G 4$,
or equivalently, $\G ~\ge~ \frac{4}{\sqrt{3}} L \ceil{\lg n}$.

In Appendix \ref{app:eigenvalues}, we describe a family of convex functions for which the above relationship between
the squares of the gradient differences and the squares of the coordinate differences 
has a $\Theta(\log^2 n)$ gap, indicating that the bound in this part of the analysis is tight.

\section{Stochastic Asynchronous Coordinate Descent (SACD)}\label{sect:stochastic-ACD}

\newcommand{\node}{\mathcal{N}}
\newcommand{\ideal}{\Delta^{\text{I}}~}
\newcommand{\idealt}[1]{\Delta^{\text{I}}_{#1}~}
\newcommand{\differ}{\Delta^{\text{D}}~}
\newcommand{\FE}{\Delta^{\text{FE}}~}
\newcommand{\FEI}{\Delta^{\text{FEI}}~}
\newcommand{\FEIt}[1]{\Delta^{\text{FEI}}_{#1}~}
\newcommand{\FEs}{\Delta^{\text{FE}}}
\newcommand{\Lksk}{L_{k_s,k}}
\newcommand{\Lksksq}{(L_{k_s,k})^2}
\newcommand{\Dxks}{\Delta~x_{k_s}^{\pi,s}}
\newcommand{\Dxku}{\Delta~x_{k_u}^{\pi,u}}
\newcommand{\Dxkt}{\Delta~x_k^{\pi,t}}
\newcommand{\Dxktnopi}{\Delta~x_k^{t}}
\newcommand{\xkpit}{x_k^{\pi,t}}
\newcommand{\xkspis}{x_{k_s}^{\pi,s}}
\newcommand{\xktaupitau}{x_{k_\tau}^{\pi,\tau}}
\newcommand{\DRIS}[2]{\Delta^{\text{RI},#1}_{#2}~}
\newcommand{\DRES}[2]{\Delta^{\text{RE},#1}_{#2}~}
\newcommand{\emp}{\emptyset}
\renewcommand{\xs}{x_s}
\newcommand{\xu}{x_u}

\newcommand{\DI}[1]{\Delta^I_{#1}}
\newcommand{\DFEIt}{\Delta_t^{\text{FEI}}~}
\newcommand{\DRI}[1]{\Delta^{\text{RI},#1}~}
\newcommand{\DRIs}[1]{\Delta_s^{\text{RI},#1}~}
\newcommand{\DRIt}[1]{\Delta_t^{\text{RI},#1}~}
\newcommand{\DRE}[1]{\Delta^{\text{RE},#1}~}
\newcommand{\DREs}[1]{\Delta_s^{\text{RE},#1}~}
\newcommand{\DREt}[1]{\Delta_t^{\text{RE},#1}~}
\newcommand{\DRs}[1]{\Delta_s^{#1}~}
\newcommand{\DRt}[1]{\Delta_t^{#1}~}
\newcommand{\singlet}{\{t\}}
\newcommand{\singles}{\{s\}}
\newcommand{\Lkskt}{L_{k_s,k_t}}
\newcommand{\Lksku}{L_{k_s,k_u}}
\newcommand{\Lktks}{L_{k_t,k_s}}
\newcommand{\Lkukt}{L_{k_u,k_t}}
\newcommand{\Lkuks}{L_{k_u,k_s}}
\newcommand{\Lkukv}{L_{k_u,k_v}}
\newcommand{\Lkvks}{L_{k_v,k_s}}
\newcommand{\Lkvkt}{L_{k_v,k_t}}
\newcommand{\Lkvku}{L_{k_v,k_u}}
\newcommand{\Lkvkw}{L_{k_v,k_w}}
\newcommand{\xkspia}{x_{k_s}^{\pi,s}}
\newcommand{\xktpit}{x_{k_t}^{\pi,t}}
\newcommand{\xkupiu}{x_{k_u}^{\pi,u}}
\newcommand{\xkvpiv}{x_{k_v}^{\pi,v}}
\newcommand{\xkwpiw}{x_{k_w}^{\pi,w}}
\newcommand{\Rbar}{\overline{R}}
\newcommand{\Rsbar}{\overline{R\cup\{s\}}}
\newcommand{\Rubar}{\overline{R\cup\{u\}}}
\newcommand{\Lkuktsq}{\left(L_{k_u,k_t}\right)^2}

\newcommand{\Dmax}{\Delta_{\max}}
\newcommand{\Dmin}{\Delta_{\min}}
\newcommand{\Es}{E_s^x}
\newcommand{\Et}{E_t^x}
\renewcommand{\FEs}{\Delta^{\text{FE}}_s}
\newcommand{\FEt}{\Delta^{\text{FE}}_t}
\newcommand{\gmaxktt}{g_{\max,k_t}^t}
\newcommand{\gminktt}{g_{\min,k_t}^t}
\newcommand{\exkss}{x_{k_s}^s}
\newcommand{\xktt}{x_{k_t}^t}

The analysis in the stochastic case is quite subtle.
The difficulty is that randomization occurs at the time the coordinates are selected,
and in this selection order the coordinates are uniformly distributed.
However, the updates occur at commit time need \emph{not} be in the same order,
owing to possible asynchronous effects such as varying computation times for different coordinates\footnote{
The amount of computation required for one coordinate can be quite different from that for another coordinate.},
varying communication delay, interference from other computations 
(e.g., due to mutual exclusion\footnote{This is needed when commiting updates to the same coordinate.}),
and interference from the operating system.
The fact that these orderings might be different was overlooked in the prior analysis by Liu and Wright~\cite{LiuW2015},
and consequently their analysis is incomplete. The present analysis is significantly different, and furthermore,
we show that the requirement of Liu and Wright for guaranteeing linear speed-up can be relaxed almost quadratically.

We will be using the starting time of the updates as reference points,
and thus \emph{future} updates in this ordering might interfere with the current update.
However, in any standard stochastic analysis, the progress is analyzed conditioning on the \emph{current} information available.
Our high-level approach to achieve this is: \emph{with the current information in hand},
give a worst-case estimate on how future updates can interfere with the current update.

While the above high-level approach seems natural, its implementation is quite non-trivial.
We think it is plausible that our approach may be also effective in analyzing other asynchronous stochastic iterative systems.

Suppose there are a total of $T$ updates.
We view the whole stochastic process as a branching tree of height $T$.
Each node in the tree corresponds to the \emph{moment} when some core randomly picks a coordinate to update,
and each edge corresponds to a possible update.
We will use $\pi$ to denote a path from the root down to some edge of the tree.

At this point, it is helpful to introduce the concept of a \emph{history}.
Suppose $\pi$ is a path of length $t$, and let $\node$ be $\pi$'s final node.
What had \emph{really} happened before $\node$, or in other words, what is the history before $\node$?
By our timing scheme and by Assumption \ref{assume:SACD-q},
we are sure that all updates strictly before time $t-q$ had committed before $\node$,
and thus \emph{all} information about such updates belong to the history.
Also, the coordinates $k_s$ for $s\in[t-q,t-1]$ were already chosen, so their identities also belong to the history;
however, \emph{some or all of their updated values might not yet belong to the history}.

We first give an analysis that deals with the simpler scenario when all updates are consistent,
i.e., no future update interferes with a current update.
This analysis will then be generalized to handle the scenario with inconsistent updates.

\subsection{SACD with Consistent Updates}\label{sect:SACD-consistent}

\eqref {eqn:F-prog} gives in expectation 
\begin{equation*}
F(\ptone) - F(\pt) ~ \ge ~\frac {1}{2n} \sum_{k=1}^n\hWk(\gkt,\pktm,\G,\Psik) + \frac 18 \G \left( \Delta \pkt \right)^2 - \frac {1}{\G} (\gkt - \tgkt)^2.
\end{equation*}
We will show that in expectation,  that for any $t'$,
$$\sum_{1\le t\le t'}  \frac {1}{\G} (\gkt - \tgkt)^2 \le \sum_{1\le t\le t'} \frac 18 \G \left( \Delta \pkt \right)^2.$$
This is not quite good enough for proving a convergence bound, for this only gives progress in an amortized sense.
We leave to the appendix and to the general non-consistent update case the additional elements.
(Basically, we consider the function $F(\pt)$ plus additional terms  that are fractions of $\G \left( \Delta \pks \right)^2$
for $s\le t$, and show that for this function the progress is of the form
$\frac {1}{2n} \sum_{k=1}^n\hWk(\gkt,\pktm,\G,\Psik) $ plus a $\frac {1}{2n}$ portion of the the additional term.
To achieve this we will need to increase $\G$ a little from the value the derivation below obtains.)

Fix a path $\pi$ of length $t$.
Let ${\Dmax} \xktt $ denote the maximum value that $\Delta \xktt$ can assume
when the first $t-q$ updates on path $\pi$ have been fixed.
Let ${\Dmin} \xktt $ denote the analogous minimum value.
Let $\gmaxktt$ denote the value of $\gkt$ used to evaluate ${\Dmax} \xktt $
and $\gminktt$ denote the value of $\gkt$ used to evaluate ${\Dmin} \xktt $.

Then, by Lemma~\ref{lem:change-of-Dp-vs-change-of-g} in Appendix \ref{sect:appendix-general},
\begin{align*}
& \left( \Dmax \xktt - \Dmin \xktt \right)^2
 \le ~ \frac {1} {\G^2} \left( \gmaxktt - \gminktt \right)^2
 \\
& \le ~ \frac {1} {\G^2} 
\left[ \sum_{ s\in [t-q,t-1]}  
\Lkskt \max \left\{ 
        \left|  \Dmax \exkss - \Dmin \exkss \right|, 
        \left| \Dmax \exkss \right|,  
        \left| \Dmin \exkss \right|
       \right\}
\right]^2
\\
& \le ~ \frac {q} {\G^2} 
\cdot \sum_ {s\in [t-q,t-1]}
\Lkskt^2 \max \left\{   
         \left(  \Dmax \exkss - \Dmin \exkss \right)^2, 
         \left( \Dmax \exkss \right)^2, 
          \left( \Dmin \exkss \right)^2 
   \right\}.
\end{align*}

Define $\left(\FE \xkpit\right)^2 = \expect{ \left( \Dmax \xktt - \Dmin \xktt \right)^2}$. Then
\begin{align}
\nonumber
\left(\FE \xkpit\right)^2
& \le ~ \frac {q} {\G^2} 
\cdot \sum_ {s\in [t-q,t-1]} \frac 1n\sum_{k_t}
    \Lkskt^2 \max \left\{   
             \left(  \Dmax \exkss - \Dmin \exkss \right)^2, 
             \left( \Dmax \exkss \right)^2, 
              \left( \Dmin \exkss \right)^2 
       \right\} \\
\label{eqn:expected-FE-bound}
& \le ~  \frac {q} {n\cdot \G^2} 
\cdot \sum_ {s\in [t-q,t-1]} \ \sum_{k_t}
    \Lres^2 \max \left\{   
             \left(  \Dmax \exkss - \Dmin \exkss \right)^2, 
             \left( \Dmax \exkss \right)^2, 
              \left( \Dmin \exkss \right)^2 
       \right\}.
\end{align}

Since $\left| \Dmin \exkss \right|, \left| \Dmax \exkss \right| \le \left| \Delta  \exkss \right| +\left(  \Dmax \exkss - \Dmin \exkss \right)$,
we have
$$\left(\Dmin \exkss \right)^2, \left(\Dmax \exkss \right)^2 \le 2\left( \Delta  \exkss \right)^2 + 2\left(  \Dmax \exkss - \Dmin \exkss \right)^2.$$

Thus
\begin{align*}
\left(\FE \xkpit\right)^2
& \le ~ \frac {q\Lres^2} {n\cdot \G^2} 
\sum_{s\in [t-q,t-1]} 
     \left[
             2\left(  \Dmax \exkss - \Dmin \exkss \right)^2
+ 2\left( \Delta  \exkss \right)^2
     \right].
\end{align*}

Now average over all nodes in a level, i.e.\ over all paths $\pi$.
Let $\left(\FEt \pkt\right)^2 = \expect{\left(\FE \xktt\right)^2}$.
Also, let $\left( \Es \right)^2 = \expect{\left( \Delta  \exkss \right)^2}$.
Then
\begin{align}
\label{eqn:simple-FE-bound}
\left(\FEt \xktt\right)^2 & \le 
~\frac {q\Lres^2} {n \cdot\G^2} 
\sum_{s\in [t-q,t-1]} 
     \left[
             \left(  \FEs \exkss \right)^2 + \left(\Es \right)^2
     \right].
\end{align}

Let $\nu = \frac  {q^2\Lres^2} {n \cdot \G^2} $ and choose $\G$ so that $\nu <1 $.
Then we can repeatedly redistribute the terms $(\FEs)^2$ to level $s$ and apply the upper bound in~\eqref{eqn:simple-FE-bound} recursively,
ad infinitum.
We see that level $s$ receives a total charge of at most 
\begin{equation}
\label{eqn:rec-FE-defn}
\nu \cdot \left(\Es \right)^2 \frac {1}{1 - \nu}.
\end{equation}

\subsection{General SACD}

When the updates might be inconsistent, the argument becomes more intricate.
We highlight the main ideas here, and defer the full proof to Appendix \ref{sect:appendix-stochastic-ACD}.
For now we assume that for each individual coordinate its updates remain consistent w.r.t.~each other.
In the appendix, we show the very minor change in our analysis needed to remove this assumption (increasing $q$ to $2q$ suffices).

The difficulty is that we wish to use the uniformly random distribution based on the coordinate selection times so as to be able to take expectations.
But from the perspective of this ordering, the update to $\pkt$ at ``time'' $t$,
may be affected by updates to some or all of the next $q$ coordinates in this ordering.
Further, depending on the random choices, which of the next $q$ coordinates have an effect may vary.

The key issue is how to do the averaging in \eqref{eqn:expected-FE-bound} for now
we have to consider future as well as past coordinate updates (i.e.\ $s>t$ on the RHS).
For the RHS terms need not be independent of the choice of $k$ any more.
To avoid this difficulty we will average over a smaller collections of \emph{substitutable} paths.
For the averaging at time $t$, we will be considering paths that extend to time $t+q$,
and the coordinates we will be allowing to vary will be those being updated in the time range $[t,t+q]$.

A particular collection will be identified by having exactly $\ell$ updates to $\pkt$ in the time interval $[t,t+q]$, for $1\le \ell \le q$,
at some specific times $t=t_{i_1}, t_{i_2}, \ldots, t_{i_{\ell}} \le t+q$.
A \emph{substitutable} sequence --- replaces $\pkt$ by $\pktp$ at these $\ell$ time steps
and leaves the other coordinates unchanged (in the sense of which coordinates are updated, not what are the values of the updates).

The key point is that if the update of $\pkt$ at time $t$ depends on one of the other coordinate updates
in the range $[t+1,+t+q]$ then the latter coordinate's update does not depend on the update to $\pkt$ nor on any of the later updates to $\pkt$.
This property remains true for all the substitutable coordinates $\pktp$.

As the interval $[t+1,t+q]$ includes at most $q$ distinct coordinates,
there are at least $n-q$ substitutable coordinates in any collection of substitutable coordinates.
Further, these collections partition the coordinate updates at level $t$.

Following the averaging over substitutable coordinates, we seek to average over all the coordinates to obtain an overall expectation.
As stated, this does not work unfortunately.
The difficulty is that the terms $\FE \xkspis$, as we define them,
need not be identical on the substitutable paths, which is needed in order to replace the terms
$\Lkskt^2$ by $\frac 1n \Lres^2$. However, it turns out that  on the RHS it suffices to use
a restricted form of these terms, which is identical across substitutable paths.
Then, having done the averaging over substitutable paths,
one can upper bound the restricted form by the general term $\FE \xkspis$.
We then average over all paths, obtaining a bound of the form \eqref{eqn:rec-FE-defn}.

\section{Parallel Asynchronous Coordinate Descent (PACD)}\label{sect:parallel-ACD}

\renewcommand{\gkttau}{g_{k_t}^{\tau}}
\newcommand{\pktaunu}{x_{k_{\tau}}^{\nu}}
\newcommand{\tgktautau}{\tilde{g}_{k_{\tau}}^{\tau}}
\newcommand{\tgktaunu}{\tilde{g}_{k_{\tau}}^{\nu}}
\newcommand{\ptau}{x^{\tau}}
\newcommand{\kstop}{k_{stop}}

Again, we want to apply Theorem \ref{thm:meta-new}. Let $H(t):= F(\pt) + A(t)$, where
$$
A(t) = \frac {1}{16} \sum_{\tau=\max\{t-q,1\}}^t \frac{(\tau+q)-t}{q}\cdot \G \left( \Delta \xktautau \right)^2.
$$

For each time $i$, let $\alpha(i)$ denote the time of the latest update to coordinate $k_i$ which is strictly before time $i$;
let it be $0$ if no such update exists.
In Appendix \ref{sect:appendix-parallel-ACD}, we will show that if $\G ~\ge~ 4q \Lmax$, then $H$ is decreasing;
also, we will use Lemma \ref{lem:W-shift} (in a spirit similar to the CCD case), and also Inequality \eqref{eqn:F-prog}, to show that
\begin{align}
&~~~~~\sum_{i=t-2r+1}^t \left[ H(i-1) - H(i) \right]\nonumber\\
&~\ge~ \sum_{j=t-2r+1}^{t-r} ~\left[\frac {1}{3r} \sum_{k=1}^n \hWk(g_k^j,x_k^{j-1},\G,\Psi_k) ~+~ \frac{1}{2q} \cdot A(j-1)\right]\nonumber\\
&~~~~~~~~~~~~~+~\left[\frac{3}{64} \sum_{i=t-2r-q}^{t-2r} \frac{i-t-2r+q}{q} ~\G\left(\Delta x_{k_i}^i\right)^2
~+~ \frac{1}{16} \sum_{i=t-2r+1}^t \G \left(\Delta x_{k_i}^i\right)^2\right]\nonumber\\
&~~~~~~~~~~~~~~~~~~~~~~-~ \frac{2}{3\G r} \sum_{i=t-2r+1}^t ~ \sum_{j=\max\{\alpha(i)+1,t-2r+1\}}^i ~ \left(g_{k_i}^j - g_{k_i}^i\right)^2
~-~ \sum_{i=t-2r+1}^t \frac 1 \G (g_{k_i}^i - \tg_{k_i}^i)^2.\label{eqn:target-PACD}
\end{align}
In the above inequality, there are four terms.
The first term matches with what is required for applying Theorem \ref{thm:meta-new},
on setting $\alpha = n/(3r)$ and $\beta = 1/2$.
Our remaining task is to show that for some sufficiently large $\G$,
the three remaining terms, in sum, are non-negative.
To achieve this, we will show that for some sufficiently large $\G$,
\begin{equation}\label{eqn:error1}
\frac{2}{3\G r} \sum_{i=t-2r+1}^t ~ \sum_{j=\max\{\alpha(i)+1,t-2r+1\}}^i ~ \left(g_{k_i}^j - g_{k_i}^i\right)^2
~\le~ \frac{1}{64} \sum_{i=t-2r+1}^t \G \left(\Delta x_{k_i}^i\right)^2
\end{equation}
and
\begin{equation}\label{eqn:error2}
\sum_{i=t-2r+1}^t \frac 1 \G (g_{k_i}^i - \tg_{k_i}^i)^2
~\le~ \frac{3}{64} \sum_{i=t-2r-q}^{t-2r} \frac{i-t+2r+q}{q} ~\G\left(\Delta x_{k_i}^i\right)^2
~+~ \frac{3}{64} \sum_{i=t-2r+1}^t ~\G\left(\Delta x_{k_i}^i\right)^2.
\end{equation}
We call the above left hand sides the first and second type of gradient error, respectively.

\paragraph{Bounding The First Type of Gradient Error.}

The analysis is almost identical to that in the CCD case, apart from the following two changes.
The first change is to construct doubling intervals of length $2^i$ for $i=0,1,2,\cdots,(\ceil{\lg r}-1)$.
The second change is to apply the following inequality in place of \eqref{eqn:cyclic-inter1}. For each doubling interval $[\tau_1,\tau_2]$,
by a use of the Cauchy-Schwarz inequality, we show that
\begin{align*}
\sum_{k=1}^n \left( g_k^{\tau_1} - g_k^{\tau_2} \right)^2
&~\le~ L^2 \sum_{k=1}^n \left(\sum_{\tau\in[\tau_1,\tau_2]\text{ and }\ktau = k} \left| \Dx_{\ktau} \right| \right)^2\\
&~\le~ L^2 \sum_{k=1}^n \kmax \sum_{\tau\in[\tau_1,\tau_2]\text{ and }\ktau = k} \left( \Dx_{\ktau} \right)^2
~=~ L^2 \kmax \sum_{\tau=\tau_1}^{\tau_2} \left(\Dx_{\ktau}\right)^2.
\end{align*}
Thus, for \eqref{eqn:error1} to hold, we need $\G \ge \frac{16}{\sqrt{3}} L \sqrt{\kmax} \lceil \lg r \rceil$.

\paragraph{Bounding The Second Type of Gradient Error.}

By Assumption \ref{assume:PACD-q} and Definition \ref{def:Lipschitz-parameters},\\
$\left|\gkt - \tgkt\right| ~\le~ \Lmax \sum_{\tau=\max\{1,t-q\}}^{t-1} \left| \Dx_{\ktau} \right|$.
Then by the Cauchy-Schwarz inequality,
\begin{equation}\label{eqn:bound-g-minus-tg}
\left(\gkt - \tgkt\right)^2 ~\le~ (\Lmax)^2 q \sum_{\tau=t-q}^{t-1} \left( \Dx_{\ktau} \right)^2,
\end{equation}
and hence
\begin{align*}
\sum_{i=t-2r+1}^t \frac 1 \G (g_{k_i}^i - \tg_{k_i}^i)^2
&~\le~ \sum_{i=t-2r+1}^t \frac{(\Lmax)^2 q}{\G} \sum_{\tau=i-q}^{i-1} \left( \Dx_{\ktau} \right)^2\\
&~\le~ \frac{(\Lmax)^2 q}{\G}\left(\sum_{i=t-2r-q}^{t-2r} (i-t+2r+q) \cdot \left( \Dx_{\ktau} \right)^2 + \sum_{i=t-2r+1}^t ~q\cdot \left( \Dx_{\ktau} \right)^2\right).
\end{align*}
Thus, for \eqref{eqn:error2} to hold, we need $\frac{(\Lmax)^2 q}{\G} ~\le~ \frac{3\G}{64q}$,
or equivalently, $\G ~\ge~ \frac{8}{\sqrt{3}} q \Lmax$.

To summarize, we need 
$\G ~\ge~ \max\left\{\frac{16}{\sqrt{3}} L \sqrt{\kmax} \lceil \lg r \rceil ~,~ \frac{8}{\sqrt{3}} q \Lmax\right\}$.

\section{Asynchronous Tatonnement in CES Fisher Markets}\label{sect:CES}

The concept of a market equilibrium was first proposed by Walras \cite{Walras1874}.
He also proposed an algorithmic approach for finding equilibrium prices, namely to
adjust prices by tatonnement: upward if there is too much demand and downward if too little.
Since then, studies of market equilibria and tatonnement have received much attention in
economics, operations research, and most recently in computer science
\cite{ABH1959,Uzawa1960,Dohtani1993,CMV2005,CF2008,CFR2010,CCR2012,CCD2013,PY2010}.
Underlying many of these works is the issue of what are plausible price adjustment mechanisms
and in what types of markets they attain a market equilibrium.

The tatonnements studied in prior work have mostly been continuous, or discrete and synchronous.
Observing that real-world market dynamics are highly distributed and hence presumably asynchronous,
Cole and Fleischer~\cite{CF2008} initiated the study of asynchronous tatonnement with their \emph{Ongoing market model},
a market model incorporating update dynamics.

Cheung, Cole and Devanur~\cite{CCD2013} showed that tatonnement is equivalent to coordinate descent
on a convex function for several classes of Fisher markets, and consequently
that a suitable synchronous tatonnement converges toward the market equilibrium in two classes of markets:
complementary-CES Fisher markets and Leontief Fisher markets.
This equivalence also enables us to apply our amortized analysis to show that
the corresponding asynchronous version of tatonnement converges toward the market equilibrium in these two classes of markets.
We note that the tatonnement for Leontief Fisher markets analyzed in~\cite{CCD2013}
had an unnatural constraint on the step sizes; our analysis removes that constraint.

As Cole and Fleischer~\cite{CF2008} argued, any realistic price dynamics must involve out-of-equilibrium trade
in order to induce the imbalances leading to price updates.
Further, they argued that simple rules with relatively low information requirements were more plausible.
The lowest imaginable level of information would be for each seller to only know the demand for the good it was selling,
and for any price updating to occur in a non-coordinated manner, i.e., asynchronously.
Accordingly, Cole and Fleischer analyzed the performance of an asynchronous tatonnement in a market that repeated,
which they named the Ongoing market. The market also incorporated warehouses (buffers) to cope with imbalances between supply and demand.

We review a few standard notions. A \emph{Fisher market} comprises a set of $n$ goods and two sets of agents, sellers and buyers.
The sellers bring the goods to market and the buyers bring money with which to buy the goods.
The trade is driven by a collection of non-negative prices $\{p_j\}_{j=1\cdots n}$ of the goods.
Without loss of generality, we can assume that each seller brings one distinct good to the market, and she is the price-setter for this good.
By normalization, we may assume that there is one unit of each good.

Each buyer $i$ starts with $e_i$ money, and has a utility function
$u_i(x_{i1},x_{i2}, \cdots ,x_{in})$ expressing her preferences:
if she prefers bundle $\{x^a_{ij}\}_{j=1\cdots n}$ to bundle $\{x^b_{ij}\}_{j=1\cdots n}$, then
$u_i ( \{x^a_{ij}\}_{j=1\cdots n} ) > u_i ( \{x^b_{ij}\}_{j=1\cdots n} )$.
At any given prices $\{p_j\}_{j=1\cdots n}$, each buyer $i$ seeks to purchase a utility-maximizing bundle of goods costing at most $e_i$.
The \emph{demand} for good $j$, denoted by $x_j$, is the total quantity of the good sought by all buyers.
The \emph{supply} of good $j$ is the quantity of good $j$ its seller brings to the market, which we have assumed to be $1$.
The \emph{excess demand} for good $j$, denoted by $z_j$, is the demand for the good minus its supply, i.e., $z_j = x_j - 1$.
Prices $\{\ps_j\}_{j=1\cdots n}$ form a \emph{market equilibrium} if, for any good $j$ with $\ps_j>0$, $z_j = 0$,
and for any good $j$ with $\ps_j=0$, $z_j \leq 0$.

The Ongoing market simply repeats the market over a sequence of time periods
called days. Sellers are allowed to update their prices as frequently as they wish,
but at least once a day (in order to ensure progress toward convergence ---
a slower rate of updating can be captured by redefining what is a day).

A CES utility function has the form
$$u\left(x_{1},x_{2},\cdots,x_{n}\right) = \left(a_{1} (x_{1})^{\rho} + a_{2} (x_{2})^{\rho}
+ \cdots + a_{n} (x_{n})^{\rho}\right)^{1/\rho},$$
where $\rho\leq 1$ and for all $ \ell$, $a_{\ell}\geq 0$.
When $\rho$ is negative, the goods form complements,
and hence the utility function is called a complementary-CES utility function.
A Leontief utility function has the form
$$u\left(x_{1},x_{2},\cdots,x_{n}\right) = \min_{\ell \in S} \left\{b_{\ell} x_{\ell}\right\},$$
where $S$ is a non-empty subset of the goods in the market, and $\forall \ell\in S$, $b_{\ell} > 0$.

Cheung, Cole and Devanur~\cite{CCD2013} showed that tatonnement is equivalent to gradient descent on a convex function $\phi$
for Fisher markets with buyers having complementary-CES or Leontief utility functions.
To be specific, \cite{CCD2013} showed that for the convex function $\phi(p) = \sum_k p_k + \sum_i \hat{u}_i(p)$,
where $\hat{u}_i(p)$ is the optimal utility that buyer $i$ attains at prices $p$, we have that $\nabla_k \phi(p) = -z_k(p)$.
The corresponding update rule is
\begin{equation}\label{eq:sync-tat-CES-rule}
p_j' = p_j\cdot \left(1 + \lambda \cdot \min\{{z}_j,1\}\right),
\end{equation}
where $\lambda > 0$ is a suitable constant. 
As the update rule was multiplicative, they assumed that the initial prices were positive.

As argued in~\cite{CF2008}, when the economic activity is occurring over time,
it is natural to base each price update for a good on the excess demand observed by its seller since the time of the last price update to her good 
(possibly weighted toward more recent sales).
This perceived excess demand can be written as the product of the length of the time interval with
an instantaneous excess demand at some specific time in this interval, 
which yields the following modification of update rule \eqref{eq:sync-tat-CES-rule}.
\begin{equation}\label{eq:async-tat-CES-rule}
p_j' = p_j\cdot \left(1 + \lambda \cdot \min\{\tilde{z}_j,1\}\cdot (t-\alpha_j(t))\right),
\end{equation}
where $\alpha_j(t)$ denotes the time of the latest update to price $j$ \emph{strictly} before time $\tau$,
$\tilde{z}_j$ is a value between the minimum and maximum instantaneous excess demands during the time interval $(\alpha_j(t), t)$,
and $\lambda > 0$ is a suitable constant.

As we will see in Appendix \ref{app:tat-anal}, having $\lambda \le 1/37$ suffices.
In comparison, in the synchronous version, $\lambda \le 1/6$ suffices.
This implies that the step sizes of the asynchronous tatonnement
can be kept at a constant fraction of those used in its synchronous counterpart.

\begin{theorem}
\label{thm:CES-Fisher-cnvge}
For  $\lambda \le 1/37$, asynchronous tatonnement price updates using rule \eqref{eq:async-tat-CES-rule} converge toward
the market equilibrium in any complementary-CES or Leontief Fisher market.
\end{theorem}

To prove Theorem \ref{thm:CES-Fisher-cnvge}, we perform an analysis similar in spirit to that of Theorem \ref{thm:main-PACD}.
However, since there is no \emph{global} Lipschitz bounds for the function $\phi(p)$ employed in~\cite{CCD2013},
we have to settle for \emph{local} Lipschitz bounds, and this introduces extra technical difficulties. 
The proof of Theorem~\ref{thm:CES-Fisher-cnvge} is deferred to Appendices~\ref{app:tat-anal} and~\ref{app:leontief}.

In an earlier version of this paper~\cite{CC2014-arxiv},
we proved the same result using a potential function,
and we extended that analysis to account for the warehouses that are present in the Ongoing Market model;
the warehouses allowed excess over and under demands to be handled in a natural way, by drawing down or adding to stocks in the warehouses.
While we are confident the present analysis can be extended in a similar way, we have yet to do so.

\section{Discussion}

Computer Science has long been concerned with the organization and manipulation of information
in the form of well-defined problems with a clear intended outcome.
But in the last 15 years, Computer Science has gained a new dimension,
in which outcomes are predicted or described, rather than designed.
This work assesses the performance of iterative problems of both types via a rather general amortized analysis.

Iterative procedures are pervasive in optimization (for example, see~\cite{FS2000} and references therein).
As these are often applied to very large problem instances,
it is desirable to perform multiple iterations or portions of iterations in parallel.
Then the challenge is to avoid inconsistent updating by the various parts of the computation that proceed in parallel.
This can be avoided by synchronizing appropriately, but this introduces its own costs and delays.
Thus, asynchronous algorithms, which can tolerate inconsistent updates, are particularly desirable;
see~\cite{LiuW2015} and~\cite[Section 4]{Wright2015} for further discussion.

Practitioners have used parallel and asynchronous optimization algorithms for a long time.
The folklore is convergence to minimum can usually be achieved, but linear speedup is rarely achievable.
In this paper, we use an amortized approach to analyze a very general form of asynchronous coordinate descent (ACD).
The approach allows us to prove the first convergence results for general parallel ACD,
to improve existing convergence results on cyclic coordinate descent and stochastic ACD,
and further, to provide sufficient conditions for achieving linear speedup in various settings,
which are given in terms of structural parameters of the underlying convex function.
These provide an interesting theoretical counterpoint to the folklore.

Iterative updating also arises in many natural systems. Examples include
bird flocking~\cite{Chazelle2009}, influence systems~\cite{Chazelle2012},
spread of information memes across the Internet~\cite{LBK2009} and market economies~\cite{CF2008}.
Many of these problems fall into the broad category of analyzing dynamical systems.
Dynamical systems are also a staple of the physical sciences;
often the dynamics are captured via elegant deterministic sets of rules
(e.g., Newton's law of motion, Maxwell's equations for electrodynamics).
The modeling of dynamical systems with intelligent agents presents new challenges
because agent behavior may not be wholly consistent or systematic.
One issue that has received relatively little attention is the timing of agents' actions.
In most prior analyses, amenable timing schemes (e.g., synchronous or round robin updates) and perfect information retrieval were assumed,
perhaps because they were more readily analyzed.
In contrast, in this paper we give an analysis that handles the more generalized timing schemes
and imperfect information retrieval that can occur due to asynchrony,
in the context of an asynchronous tatonnement price update rule for market economies;
this is achieved by extending the coordinate descent analysis to apply to this setting.
We believe it is worthwhile to investigate whether the insight our approach provides will prove helpful in
achieving more \emph{realistic} analyses of other agent-based dynamical systems in which asynchrony is natural.

The amortized analysis of tatonnement given here has its roots in some of the earlier amortized analysis of Ongoing markets~\cite{CFR2010, CCR2012}.
Since updates in these markets can happen at arbitrary times, it is natural to model time as being continuous.
This leads to potential functions with a mix of integral and discrete terms,
where progress (i.e., the decrease in the potential function) is continuous, 
and at the events, price updates, the potential function is guaranteed not to increase.
In this paper, instead, we adapt the PACD analysis to handle arbitrary update times, but without
the need for integral terms (although an equivalent formulation with such terms is possible).
We note that the amortization in this analysis is quite distinct from the one in~\cite{CCR2012}.

\section{Acknowledgment}

We thank several anonymous reviewers for their helpful and thoughtful suggestions regarding earlier versions of this paper.
We also thank Yixin Tao for his perceptive comments.

\newpage
\appendix

\section{Some Basic Lemmas and Facts}\label{sect:appendix-general}

\subsection{Notation and Basic Facts}\label{subsect:discrete-basic}

Recall we assume that at each time, there is exactly one coordinate value being updated.
For any time $\tau$, let $\ktau$ denote the coordinate which is updated at time $\tau$, and
let $\Dp_{\ktau,\tau} := x^\tau_{\ktau} - x^{\tau-1}_{\ktau}$.

Since, for each coordinate $j$, the parameter $\G_j$ and the function $\Psi_j$ remain unchanged throughout the ACD process, to avoid clutter,
we use the shorthand
$$\Wj(d,g,x) := W(d,g,x,\G_j,\Psi_j)~~~~~~~~\hWj(g,x) := \hW(g,x,\G_j,\Psi_j)~~~~~~~~\hdj(g,x) := \hd(g,x,\G_j,\Psi_j).$$
Note that $\Wj(0,g,x)=0$; thus $\hWj(g,x) \ge 0$.

Let $[n]$ denote the set of coordinates $\{1,2,\cdots,n\}$.
In this proof, $\Psi$ will always denote a function $\rr\ra\rr$ which is univariate, proper, convex and lower semi-continuous.

It is well-known that for any $k\in [n]$, $x\in\rr^n$ and $r\in\rr$,
\begin{equation}\label{eq:upper-sandwich}
f(x+r \ve_k) \leq f(x) + \nabla_k f(x) \cdot r + \frac{L_k}{2} r^2.
\end{equation}
Also, for any $x,r\in\rr^n$,
\begin{equation}\label{eq:global-upper-sandwich}
f(x+r) \leq f(x) + \sum_{k=1}^n \nabla_k f(x) \cdot r_k + \frac{L}{2} \sum_{k=1}^n (r_k)^2.
\end{equation}

Finally, for any $\tau \ge 1$, we define the progress as follows:
\begin{equation}\label{eq:def-PRG}
\PRG(\tau-1) ~:=~ \sum_{k=1}^n ~\hWk(\nabla_k f(x^{\tau-1}),x_{k}^{\tau-1},\G,\Psik).
\end{equation}

\subsection{Some Lemmas about the Functions $\hW$ and $\hd$, and Proof of Lemma \ref{lem:W-shift}}\label{subsect:Wd}

We state five technical lemmas concerning the functions $\hW$ and $\hd$.
The following fact, which follows directly from the definition of $\hW$, will be used multiple times:
\begin{equation} \label{eqn:What-as-fn-Gamma}
\text{If }0 < \G < \G',~\text{then}~\forall g,x\in\rr,~\hW(g,x,\G,\Psi) \geq \hW(g,x,\G',\Psi).
\end{equation}

\begin{lemma}[{Three-Point Property, \cite[Lemma 3.2]{CT1993}}]\label{lem:3pp}
For any proper, convex and lower semi-continuous function $Y:\rr\ra\rr$ and for any $\xminus\in\rr$, let
$\xplus := \argmax_{x\in \rr} \left\{-Y(x) - \G (x-\xminus)^2 / 2\right\}$.
Then for any $x'\in\rr$,
$$Y(x') + \frac{\G}{2} (x'-\xminus)^2 \geq Y(\xplus) + \frac{\G}{2} (x'-\xplus)^2 + \frac{\G}{2} (\xplus - \xminus)^2.$$
\end{lemma}

\begin{lemma}[{\cite[Lemma 4]{TsengY2009}}]\label{lem:change-of-Dp-vs-change-of-g}
For any $g_1,g_2,x\in\rr$ and $\G\in\rrplus$,
$$\left| \hd(g_1,x,\G,\Psi) - \hd(g_2,x,\G,\Psi) \right| \leq \frac{1}{\G}\cdot \left|g_1 - g_2\right|.$$
\end{lemma}

\begin{lemma}\label{lem:prog-better-than-quadratic}
For any $g,x\in\rr$ and $\G\in\rrplus$, $\hW(g,x,\G,\Psi) \geq \frac{\G}{2} \left(\hd(g,x,\G,\Psi)\right)^2$.
\end{lemma}

\begin{pf}
We apply Lemma \ref{lem:3pp} with $x^- = x' = 0$ and $Y(d) = gd - \Psi(x) + \Psi(x+d)$.
Then $W(d,g,x,\G,\Psi) = -Y(d) - \G d^2/2$, and hence $x^+$, as defined in Lemma \ref{lem:3pp}, equals $\hd(g,x,\G,\Psi)$.
These yield
$$Y(0) \geq Y(\hd(g,x,\G,\Psi)) + \G\cdot \left(\hd(g,x,\G,\Psi)\right)^2.$$
Since $Y(0) = 0$ and $-Y(\hd(g,x,\G,\Psi)) = \hW(g,x,\G,\Psi) + \frac{\G}{2} \left(\hd(g,x,\G,\Psi)\right)^2$, we are done.
\end{pf}

\begin{lemma}\label{lem:prog-with-Dt}
For any $g,x\in\rr$, $\G\in\rrplus$ and $0\leq q \leq 1$,
$$W(q\cdot \hd(g,x,\G,\Psi),g,x,\G,\Psi) \geq q\cdot \hW(g,x,\G,\Psi).$$
\end{lemma}

\begin{pf}
The lemma can be proved easily using the fact that $W(d,g,x,\G,\Psi)$ is a concave function of $d$, as follows:
\begin{align*}
W(q\cdot \hd(g,x,\G,\Psi),g,x,\G,\Psi) &\geq (1-q) \cdot W(0,g,x,\G,\Psi) + q \cdot W(\hd(g,x,\G,\Psi),g,x,\G,\Psi)\\
&= (1-q) \cdot 0 + q \cdot \hW(g,x,\G,\Psi).
\end{align*}
\end{pf}

We finish this subsection by using Lemmas \ref{lem:change-of-Dp-vs-change-of-g} and \ref{lem:prog-better-than-quadratic}
to prove Lemma \ref{lem:W-shift}.

\begin{rlemma}{lem:W-shift}
For any $g_1,g_2,x\in\rr$ and $\G\in\rrplus$, $\hW(g_2,x,\G,\Psi) \leq \frac 32 \hW(g_1,x,\G,\Psi) + \frac{2}{\G} (g_1 - g_2)^2$.
\end{rlemma}

\begin{pf}
To avoid clutter, we use the shorthand $\hd(g_i) := \hd(g_i,x,\G,\Psi)$ for $i=1,2$.
\begin{align*}
& \hW(g_1,x,\G,\Psi)\\
&~=~ \max_{d\in \rr} W(d,g_1,x,\G,\Psi)\\
&~\geq~ W(\hd(g_2),g_1,x,\G,\Psi)\\
&~=~ -g_1 \cdot \hd(g_2) - \G \cdot \hd(g_2)^2 / 2 + \Psi(x) - \Psi(x + \hd(g_2))\\
&~=~ -g_2 \cdot \hd(g_2) - \G \cdot \hd(g_2)^2 / 2 + \Psi(x) - \Psi(x + \hd(g_2)) + (g_2 - g_1) \cdot \left[\hd(g_1) + (\hd(g_2) - \hd(g_1))\right]\\
&~\geq~ \hW(g_2,x,\G,\Psi) - |g_1 - g_2| \cdot \left|\hd(g_1)\right| - |g_1 - g_2| \cdot \left|\hd(g_2) - \hd(g_1)\right|\\
&~\geq~ \hW(g_2,x,\G,\Psi) - |g_1 - g_2| \cdot \left|\hd(g_1)\right| - \frac{1}{\G} (g_1 - g_2)^2\comm{By Lemma \ref{lem:change-of-Dp-vs-change-of-g}}\\
&~\geq~ \hW(g_2,x,\G,\Psi) - \frac{1}{\G} (g_1 - g_2)^2  - \frac{\G}{4} (\hd(g_1))^2  - \frac{1}{\G} (g_1 - g_2)^2\comm{AM-GM ineq.}\\
&~\geq~ \hW(g_2,x,\G,\Psi) - \frac{2}{\G} (g_1 - g_2)^2 - \frac 12 \hW(g_1,x,\G,\Psi).\comm{By Lemma \ref{lem:prog-better-than-quadratic}}
\end{align*}
\end{pf}

\subsection{Proof of Theorem \ref{thm:meta-new}}

Recall the definition of $\PRG(t-1)$ in \eqref{eq:def-PRG}.
We will use the following lemma from~\cite[Lemmas 4,6]{RiT2014}.
We provide a proof here for completeness.

\begin{lemma}[{(\cite[Lemmas 4,6]{RiT2014})}]\label{lem:good-progress}
~\\
(a) Suppose that $f,F$ are strongly convex with parameters $\mu_f,\mu_F > 0$ respectively, 
and also suppose that $\G \ge \muf$. Then
$$\PRGe(t-1) ~\geq~ \frac{\mu_F}{\mu_F + \G - \mu_f}\cdot F(\ptone).$$
(b) Let $\calR := \min_{\ps\in P^*} \|\ptone - \ps\|$. Then
$$\PRGe(t-1) ~\geq~ \min\left\{\frac 12~,~\frac{F(\ptone)}{2\G~\calR ^2}\right\}\cdot F(\ptone).$$
\end{lemma}

\begin{pf}
First of all, we show a lower bound for $\PRG(t-1)$, which will be used to prove both (a) and (b).
\begin{align*}
\PRG(t-1) &~=~ \sum_{k=1}^n ~\max_{d_k\in\rr} \left\{-\nabla_k f(\ptone) \cdot d_k - \G\cdot (d_k)^2 / 2 + \Psi_k(\ptone_k) - \Psi_k(\ptone_k + d_k)\right\}\\
&~\ge~ \max_{d\in\rr^n}~\left\{\sum_{k=1}^n\left[-\nabla_k f(\ptone) \cdot d_k - \G\cdot (d_k)^2 / 2 + \Psi_k(\ptone_k) - \Psi_k(\ptone_k + d_k)\right]\right\}.
\end{align*}

When $f$ is strongly convex with parameter $\mu_f$, for any $d\in\rr^n$,
$$f(\ptone+d) \geq f(\ptone) + \sum_{k=1}^n \nabla_k f(\ptone) \cdot d_k + \frac{\mu_f}{2} \sum_{k=1}^n (d_k)^2.$$
Thus
\begin{align}
\PRG(t-1) &~\geq~ \max_{d\in\rr^n} \left\{f(\ptone) - f(\ptone + d) - \frac{\G-\mu_f}{2} \sum_{k=1}^n (d_k)^2 + \sum_{k=1}^n\left[\Psi_k(\ptone_k) - \Psi_k(\ptone_k + d_k)\right]\right\}\nonumber\\
&~=~ \max_{d\in\rr^n} \left\{ F(\ptone) - F(\ptone + d) - \frac{\G-\mu_f}{2} \sum_{k=1}^n (d_k)^2 \right\}\nonumber\\
&~\geq~ \max_{0\leq \beta\leq 1} \left\{ F(\ptone) - F\left(\beta \ps + (1-\beta) \ptone\right)  - \frac{(\G-\mu_f)\beta^2}{2} \sum_{k=1}^n (\ptone_k - \ps_k)^2 \right\}.\label{eq:good-progress-key}
\end{align}

To prove (a), we apply the following characterization of strong convexity of $F$: for any $0\leq \beta\leq 1$,
$$F\left(\beta \ps + (1-\beta) \ptone\right) ~\leq~ \beta \cdot F(\ps) ~+~ (1-\beta) \cdot F(\ptone) ~-~ \frac{\mu_F\beta(1-\beta)}{2}\sum_{k=1}^n (\ptone_k - \ps_k)^2.$$
Note that $F(\ps) = F^* = 0$.
By \eqref{eq:good-progress-key},
\begin{align*}
\PRG(t-1) &~\geq~ \max_{0\leq \beta\leq 1} \left\{ \beta \cdot F(\ptone) ~+~ \frac{\mu_F\beta(1-\beta)-(\G-\mu_f)\beta^2}{2} \sum_{k=1}^n (\ptone_k - \ps_k)^2 \right\}\\
&~\geq~ \left.\left(\beta \cdot F(\ptone) ~+~ \frac{\mu_F\beta(1-\beta)-(\G-\mu_f)\beta^2}{2}
\sum_{k=1}^n (\ptone_k - \ps_k)^2\right)\right|_{\beta = \mu_F/(\mu_F+\G-\mu_f)}\\
&~=~ \frac{\mu_F}{\mu_F + \G - \mu_f}\cdot F(\ptone).
\end{align*}
Note that the constraint  $\beta \le 1$ forces $\G \ge \muf$.

To prove (b), let $\ps$ denote a point in $P^*$ such that $\|\ptone - \ps\| \leq \calR$,
where $P^*$ is the set of minimum points for $F$.
By the convexity of $F$, when $0\leq \beta\leq 1$, $F(\beta \ps + (1-\beta)\ptone) \leq (1-\beta)\cdot F(\ptone)$.
Since $f$ is convex, $\mu_f\geq 0$ always. From \eqref{eq:good-progress-key},
$$\PRG(t-1) ~\geq~ \max_{0\leq \beta\leq 1} \left\{ \beta \cdot F(\ptone) - \frac{\G \beta^2}{2} \calR ^2 \right\}.$$
The R.H.S.~of the above inequality is a maximization of a quadratic function of $\beta$, which can be easily solved to yield the lower bound in (b).
\end{pf}

The following proof is a generalization of the above result to account for rounds of length $r>1$, and our amortization function $A(t)$.

\begin{pfof}{Theorem \ref{thm:meta-new}}
~\\
\noindent\textbf{Proof of (i).}~By the second assumption and Lemma \ref{lem:good-progress},
\begin{align*}
\sum_{i=t-2r+1}^t \left[ H(i-1) - H(i) \right]
&~\ge~
\sum_{i=t-2r+1}^{t-r} \left[~\frac \alpha n ~ \PRG(i-1) ~+~ \frac{\beta}{q} \cdot A(i-1)~\right]\\
&~\ge~ 
\sum_{i=t-2r+1}^{t-r} \left[~\frac \alpha n \cdot \frac{\mu_F}{\mu_F + \G - \mu_f} ~ F(x^{i-1}) ~+~ \frac{\beta}{q} \cdot A(i-1)~\right]\\
&~\ge~ \sum_{i=t-2r+1}^{t-r} \delta \cdot H(i-1),
\end{align*}
where $\delta ~:=~ \min\left\{ \frac \alpha n \cdot \frac{\mu_F}{\mu_F + \G - \mu_f} ~,~ \frac{\beta}{q}\right\}$.

By the assumption that $H$ is decreasing,
$$\sum_{i=t-2r+1}^{t-r} \delta \cdot H(i-1) ~\ge~ \frac \delta 2\cdot \sum_{i=t-2r+1}^{t} H(i-1).$$

Combining all the above yields
$$
\sum_{i=t-2r+1}^{t} H(i) ~~\le~~ \left( 1-\frac \delta 2 \right) \cdot \sum_{i=t-2r+1}^{t} H(i-1)
~~=~~ \left( 1-\frac \delta 2 \right) \cdot \sum_{i=t-2r}^{t-1} H(i).
$$
For any $t\ge 2r$, iterating the above inequality $(t-2r+1)$ times yields
$$
\sum_{i=t-2r+1}^{t} H(i) ~~\le~~ \left( 1-\frac \delta 2 \right)^{t-2r+1} ~ \sum_{i=0}^{2r-1} H(i).
$$
Since $H$ is decreasing, the summation in LHS is at least $2r \cdot H(t)$,
while the summation in RHS is at most $2r \cdot H(0)$. Thus,
$$
H(t) ~~\le~~ \left( 1-\frac \delta 2 \right)^{t-2r+1} \cdot H(0).
$$
To finish the proof, note that since $A(t)$ is non-negative by assumption, $F(\pt) \le H(t)$,
and note that since $A(0) = 0$, $H(0) = F(\pc)$.

\bigskip

\noindent\textbf{Proof of (ii).}~By the second assumption and Lemma \ref{lem:good-progress},
\begin{align*}
\sum_{i=t-2r+1}^t \left[ H(i-1) - H(i) \right]
&~\ge~
\sum_{i=t-2r+1}^{t-r} \left[~\frac \alpha n ~ \PRG(i-1) ~+~ \frac{\beta}{q} \cdot A(i-1)~\right]\\
&~\ge~ 
\sum_{i=t-2r+1}^{t-r} \left[~\frac \alpha n \cdot \min\left\{\frac 12~,~\frac{F(x^{i-1})}{2\G~\calR ^2}\right\}\cdot F(x^{i-1}) ~+~ \frac{\beta}{q} \cdot A(i-1)~\right].
\end{align*}

For each $i$, there are two possible cases:
\begin{itemize}
\item If $F(x^{i-1}) \le A(i-1)$, then $A(i-1) ~\ge~ \frac{H(i-1)}{2}$, thus
\begin{align*}
\frac \alpha n \cdot \min\left\{\frac 12~,~\frac{F(x^{i-1})}{2\G~\calR ^2}\right\}\cdot F(x^{i-1}) ~+~ \frac{\beta}{q} \cdot A(i-1)
&~\ge~ \frac{\beta}{2q}\cdot H(i-1).
\end{align*}
\item If $F(x^{i-1}) > A(i-1)$, then $F(x^{i-1}) ~>~ \frac{H(i-1)}{2}$, thus
\begin{align*}
\frac \alpha n \cdot \min\left\{\frac 12~,~\frac{F(x^{i-1})}{2\G~\calR ^2}\right\}\cdot F(x^{i-1}) ~+~ \frac{\beta}{q} \cdot A(i-1)
&~>~ \frac \alpha {2n} \cdot \min\left\{\frac 12~,~\frac{H(i-1)}{4\G~\calR ^2}\right\}\cdot H(i-1).
\end{align*}
\end{itemize}
Since $H$ is a decreasing function, $H(i-1) ~\le~ H(0) = F(\pc)$.
Thus, unconditionally, we have
\begin{align*}
\frac \alpha n \cdot \min\left\{\frac 12~,~\frac{F(x^{i-1})}{2\G~\calR ^2}\right\}\cdot F(x^{i-1}) ~+~ \frac{\beta}{q} \cdot A(i-1)
&~\ge~ \min \left\{ \frac{\beta}{2q} ~,~ \frac{\alpha}{4n} ~,~ \frac{H(i-1)}{4\G ~\calR ^2} \right\} \cdot H(i-1)\\
&~\ge~ \min \left\{ \frac{\beta}{2q~F(\pc)} ~,~ \frac{\alpha}{4n~F(\pc)} ~,~ \frac{1}{4\G ~\calR ^2} \right\} \cdot H(i-1)^2.
\end{align*}
Note that the term $\min \left\{ \frac{\beta}{2q~F(\pc)} ~,~ \frac{\alpha}{4n F(\pc)} ~,~ \frac{1}{4\G ~\calR ^2} \right\}$
is independent of $i$. We let $\varepsilon$ denote it.

Now, we have
\begin{align*}
\sum_{i=t-2r+1}^t \left[ H(i-1) - H(i) \right] &~\ge~ \sum_{i=t-2r+1}^{t-r} \varepsilon~ H(i-1)^2\\
&~\ge~ \frac{\varepsilon}{r} \left(\sum_{i=t-2r+1}^{t-r} H(i-1)\right)^2\comm{by the Power Mean Inequality}\\
&~\ge~ \frac{\varepsilon}{r} \left(\frac 12 \sum_{i=t-2r+1}^{t} H(i-1)\right)^2\comm{$H$ is a decreasing function}\\
&~=~ \frac{\varepsilon}{4r}\left(\sum_{i=t-2r+1}^{t} H(i-1)\right)^2
\end{align*}
For brevity, let $S_\tau ~:=~ \sum_{i=\tau-2r+1}^\tau H(i)$. Then the above inequality translates to
$$
S_{t-1} - S_t ~\ge~ \frac{\varepsilon}{4r} \left(S_{t-1}\right)^2.
$$
Note that $S_{t-1} ~\ge~ S_t ~\ge~ 0$.
Dividing both sides by $S_{t-1} \cdot S_t$ yields
$$
\frac{1}{S_t} - \frac{1}{S_{t-1}} ~\ge~ \frac{\varepsilon}{4r} \frac{S_{t-1}}{S_t} ~\ge~ \frac{\varepsilon}{4r}.
$$
Iterating the above inequality $t-2r+1$ times yields
$$
\frac{1}{S_t} - \frac{1}{S_{2r-1}} ~\ge~ \frac{\varepsilon}{4r} (t-2r+1)
$$
and hence
$$
S_t ~\le~ \frac{1}{\frac{1}{S_{2r-1}} + \frac{\varepsilon}{4r} (t-2r+1)}
~=~ \frac{S_{2r-1}}{1 + \frac{\varepsilon}{4} \cdot \frac{S_{2r-1}}{r} \cdot (t-2r+1)}.
$$
Since $H$ is a decreasing function, $S_t ~\ge~ 2r\cdot H(t)$, while $2r \cdot H(2r-1) ~\le~ S_{2r-1} ~\le~ 2r\cdot H(0) ~=~ 2r \cdot F(\pc)$. Thus,
$$
H(t) ~\le~ \frac{F(\pc)}{1 + \frac{\varepsilon}{4} \cdot \frac{2r\cdot H(2r-1)}{r} \cdot (t-2r+1)}
~=~ \frac{F(\pc)}{1 + \frac{\varepsilon}{2} \cdot H(2r-1) \cdot (t-2r+1)}.
$$
We finish the proof by noting that $F(\pt) ~\le~ H(t)$.
\end{pfof}

\subsection{Proofs of Lemmas \ref{lem:F-prog-one} and \ref{lem:F-prog-two}}

Lemmas \ref{lem:F-prog-one} and \ref{lem:F-prog-two} follows directly from the lemma below.
In the following lemma we allow the use of distinct $\G_j$  for each coordinate, as this generalization will
be used in the tatonnement analysis. 

\begin{lemma}\label{lem:discrete-improvement-of-F}
Suppose there is an update to coordinate $j$ at time $t$ according to rule \eqref{eq:update-rule},
and suppose that $\G_j \ge L_j$.
Let $g_j = \nabla_j f(\ptone)$ and $\tg_j = \nabla_j f(\tx)$.
Then
\begin{align*}
F(\ptone) - F(\pt) &~\geq~ \frac{\G_j}{2} (\Dp_{j,t})^2 - |g_j - \tg_j|\cdot |\Dp_{j,t}|\\
\text{and}~~~~~~~~~~~~~~~~~~~~~~~
F(\ptone) - F(\pt) &~\geq~ \hW(g_j,\ptone_j,\G_j,\Psi_j) - \frac{1}{\G_j}(g_j - \tg_j)^2.~~~~~~~~~~~~~~~~~~~~~~~~~
\end{align*}
\end{lemma}

\begin{pf}
To avoid clutter, we use the shorthand $d_j := \hdj(g_j,\ptone_j)$ and $\td_j := \hdj(\tg_j,\ptone_j)$.
By update rule \eqref{eq:update-rule}, $\td_j = \Dp_{j,t}$.
\begin{align*}
F(\pt) &~=~ f(\pt) + \Psi_j (x_j^t) + \sum_{k\neq j} \Psi_k(x_k^t)\\
&~\leq~ f(\ptone) + g_j \td_j + \frac{\G_j}{2} (\td_j)^2 + \Psi_j (\ptone_j + \td_j) + \sum_{k\neq j} \Psi_k(\ptone_k) \comm{By \eqref{eq:upper-sandwich}, \eqref{eq:update-rule}, and the assumption $\G_j\ge L_j$}\\
&~=~ F(\ptone) + \tg_j \td_j  + \frac{\G_j}{2} (\td_j)^2 - \Psi_j(\ptone_j) + \Psi_j (\ptone_j + \td_j) + (g_j - \tg_j)\td_j\\
&~=~ F(\ptone) - \hWj(\tg_j,\ptone_j) + (g_j - \tg_j)\td_j.
\end{align*}

Hence,
$$F(\ptone) - F(\pt) ~\geq~ \hWj(\tg_j,\ptone_j) - (g_j - \tg_j) \td_j.$$
	
Then we can apply Lemma \ref{lem:prog-better-than-quadratic} to prove the first inequality in Lemma \ref{lem:discrete-improvement-of-F}:
$$F(\ptone) - F(\pt) ~\geq~ \hW_j(\tg_j,\ptone_j) - (g_j - \tg_j) \td_j ~\geq~ \frac{\G_j}{2} (\td_j)^2 - |g_j - \tg_j|\cdot |\td_j|.$$

We prove the second inequality in Lemma \ref{lem:discrete-improvement-of-F} as follows:
\begin{align*}
F(\ptone) - F(\pt) &~\geq~ \hW_j(\tg_j,\ptone_j) - (g_j - \tg_j) \td_j\\
&~\geq~ W_j(d_j,\tg_j,\ptone_j) - (g_j - \tg_j) \td_j\\
&~=~ W_j(d_j,g_j,\ptone_j) + (g_j - \tg_j) d_j - (g_j - \tg_j) \td_j\\
&~=~ \hWj(g_j,\ptone_j) + (g_j - \tg_j) (d_j - \td_j)\\
&~\geq~ \hWj(g_j,\ptone_j) - |g_j - \tg_j| \cdot |d_j - \td_j|\\
&~\geq~ \hWj(g_j,\ptone_j) - \frac{1}{\G_j}(g_j - \tg_j)^2 \comm{By Lemma \ref{lem:change-of-Dp-vs-change-of-g}}.
\end{align*}
\end{pf}

\newpage

\section{Stochastic Asynchronous Coordinate Descent (SACD)}\label{sect:appendix-stochastic-ACD}
 
\renewcommand{\Ap}{A^+}
\renewcommand{\Am}{A^-}

\newcommand{\Dmaxt}{\Delta_{\max}^t}
\newcommand{\Dmaxu}{\Delta_{\max}^u}

\newcommand{\Dmint}{\Delta_{\min}^t}
\newcommand{\Dminu}{\Delta_{\min}^u}

\newcommand{\xv}{x_v}
\newcommand{\gmaxkupiu}{g_{\max,k_u}^{\pi,u}}
\newcommand{\gmaxkupipiu}{g_{\max,k_u}^{\pi,\pi',u}}
\newcommand{\gmaxkupitRu}{g_{\max,k_u}^{\pi,t,R,u}}
\newcommand{\gmaxkupitphiu}{g_{\max,k_u}^{\pi,t,\phi,u}}
\newcommand{\gminkupiu}{g_{\min,k_u}^{\pi,u}}
\newcommand{\gminkupipiu}{g_{\min,k_u}^{\pi,\pi',u}}
\newcommand{\gminkupitRu}{g_{\min,k_u}^{\pi,t,R,u}}
\newcommand{\gminkupitphiu}{g_{\min,k_u}^{\pi,t,\phi,u}}
\newcommand{\gkspis}{ {g}_{k_s}^{\pi,s}}
\newcommand{\gkspipis}{ {g}_{k_s}^{\pi,\pi',s}}
\newcommand{\tgkspis}{ \tilde{g}_{k_s}^{\pi,s}}
\newcommand{\tgkspipis}{ \tilde{g}_{k_s}^{\pi,\pi',s}}
\newcommand{\gkupiu}{ {g}_{k_u}^{\pi,u}}
\newcommand{\gkupipiu}{ {g}_{k_u}^{\pi,\pi',u}}
\newcommand{\tgkupiu}{ \tilde{g}_{k_u}^{\pi,u}}
\newcommand{\tgkupipiu}{ \tilde{g}_{k_u}^{\pi,\pi',u}}
\newcommand{\pip}{\pi'}
\newcommand{\pku}{p_{k_u}}
\newcommand{\pkv}{p_{k_v}}
\newcommand{\ubar}{\overline{\{u\}}}
\newcommand{\xkss}{x_{k_s}^s}
\newcommand{\xkuu}{x_{k_u}^u}
\newcommand{\xkspipis}{x_{k_s}^{\pi,\pi',s}}
\newcommand{\xkspiRs}{x_{k_s}^{\pi,R,s}}
\newcommand{\xkspius}{x_{k_s}^{\pi,\{u\},s}}
\newcommand{\xkspiphis}{x_{k_s}^{\pi,\phi,s}}
\newcommand{\Ru}{R\cup \{u\}}
\newcommand{\xkspiRus}{x_{k_s}^{\pi,\Ru,s}}
\newcommand{\xkupiRu}{x_{k_u}^{\pi,R,u}}
\newcommand{\xkupipiu}{x_{k_u}^{\pi,\pi',u}}
\newcommand{\xkpupiRpu}{x_{k'_u}^{\pi,R',u}}
\newcommand{\xkupipphiu}{x_{k_u}^{\pi',\phi,u}}
\newcommand{\xkpupiRu}{x_{k'_u}^{\pi,R,u}}
\newcommand{\xkupiRpu}{x_{k_u}^{\pi,R',u}}
\newcommand{\xkpupiphiu}{x_{k'_u}^{\pi,\phi,u}}
\newcommand{\xkpupipRu}{x_{k'_u}^{\pi',R,u}}
\newcommand{\xkupiphiu}{x_{k_u}^{\pi,\phi,u}}

\newcommand{\FEst}{\Delta^{\text{FE},t}_s~}
\newcommand{\FEu}{\Delta^{\text{FE}}_u} 
\newcommand{\FEut}{\Delta^{\text{FE},t}_u~} 
\newcommand{\FEuRt}{\Delta^{\text{FE},R,t}_u~} 
\newcommand{\FEsRut}{\Delta^{\text{FE},\Ru,s}_u~} 
\newcommand{\FEtpi}{\Delta^{\text{FE}}_{t+i}}
\newcommand{\FEtpip}{\Delta^{\text{FE}}_{t+i+1}}
\newcommand{\tone}{t-1}

We begin by fixing a path $\pi$ from the root to a leaf to analyze the update to $ \xktpit$.
Because the commit ordering need not match the start ordering of updates,
when considering the updates of other variables that can influence the updates of  $ \xktpit$
we will want to exclude certain updates, which is done via the following definition of $\Rbar$.
For any set $R \subset [t-q,t+q]$, let
$$\Rbar := \left\{ s \le t+q \,\Big|\, \exists~r\in R \text{ such that }s \ge r \text{ and }k_s = k_r \right\}.$$
In fact, the only sets $R$ we consider are $R = \phi$ or $R$ is a singleton.

Let ${\Delta_{\max}^t} \xkupiRu $ denote the maximum value that $\Delta \xkupiu$ can assume
when the first $t-q -1$ updates on path $\pi$ have been fixed,
assuming the update does not read any of the updates at times in $\Rbar$,
nor any of the variables updated at times $v > t+q$.
Let ${\Delta_{\min}^t} \xkupiRu $ denote the analogous minimum value.
Let $\gmaxkupitRu$ denote the value of $\gkupiu$ used to evaluate ${\Dmaxt} \xkupiRu $
and $\gminkupitRu$ denote the value of $\gkupiu$ used to evaluate ${\Dmint} \xkupiRu $.
Note that $\Delta_{\max}^t \xkupiRu = {\Delta_{\min}^t} \xkupiRu = \Delta \xkupiu$
if $u < t-q$.
Also, $\Delta_{\max}^t \xkupiRu \ge \Delta_{\max}^t \xkupiphiu$ and
${\Delta_{\min}^t} \xkupiRu \le {\Delta_{\min}^t} \xkupiphiu$.

Then, by Lemma~\ref{lem:change-of-Dp-vs-change-of-g},
\begin{align*}
& \left( \Dmaxt \xkupiphiu - \Dmint \xkupiphiu \right)^2
 \le ~ \frac {1} {\G^2} \left( \gmaxkupitphiu - \gminkupitphiu \right)^2
 \\
& \le ~ \frac {1} {\G^2} 
\left[ \sum_{ s\in [\min\{t-q,u-q\},\min\{ t+q, u+q\}] \setminus \ubar}  
   \Lksku \cdot \max \left\{ 
            \left|  \Dmaxt \xkspius - \Dmint \xkspius \right|, 
            \left| \Dmaxt \xkspius \right|,  
 \right. \right.
\\
&~~~~~~~~~~~ ~~~~~~~~~~~ ~~~~~~~~~~~ ~~~~~~~~~~~ ~~~~~~~~~~~ ~~~~~~~~~~~ ~~~~~~~~~~~ ~~~~~~~~~~~ 
~~~~~~~~~~~~~~~~~~~~~~
\left.\left.
            \left| \Dmint \xkspius \right|
           \right\}
\right]^2
\\
& \le ~ \frac {2q} {\G^2} 
\cdot \sum_{ s\in [\min\{t-q,u-q\}, \min\{ t+q, u+q\}] \setminus \ubar}  
    \Lksku^2 \cdot \max \left\{   
             \left(  \Dmaxt \xkspius - \Dmint \xkspius \right)^2, 
             \left( \Dmaxt \xkspius \right)^2, \right.\\
& ~~~~~~~~~~~ ~~~~~~~~~~~ ~~~~~~~~~~~ ~~~~~~~~~~~ ~~~~~~~~~~~ ~~~~~~~~~~~ ~~~~~~~~~~~ ~~~~~~~~~~~ 
~~~~~~~~~~~~~~~~~~~~~~\left.
              \left( \Dmint \xkspius \right)^2 
       \right\}.
\end{align*}
We are interested in those $k_u$ whose update is not determined when the update of $k_t$ starts
and whose update can affect the update of $\pkt$. Consequently,
for a fixed $u$, we only consider $t$ in the range $[u-q, u+q]$, for a larger $t$ means that
$\Delta \xkupiu$ is already fixed, and a smaller $t$ means that the update of $x_{k_u}$ does not affect the update of $\pkt$.

For a fixed $u$, we maximize over  $t \in[u-q,u+q]$ and then
average over a collection $k'_u$ of substitutable $k_u$
(letting $\pip$ denote the substitutable path in which $k'_u$ replaces $k_u$), which gives
\begin{align*}
& \expect{ \max_{t}\left\{ \left( \Dmaxt \xkpupiphiu - \Dmint \xkpupiphiu \right)^2 \right\}~ \Big| ~ k'_u ~\text{and}~ k_u ~\text{are substitutable}~ } \\
	&~~~~~~~~~~~ \le ~ \frac {2q} {(n-q)\G^2} 
\sum_{ s\in [u-2q,  u+q] \setminus \ubar}  
    \Lres^2  \cdot \max_{t}  \left\{ 
             \left(  \Dmaxt \xkspius - \Dmint \xkspius \right)^2, \right. \\
&~~~~~~~~~~~ ~~~~~~~~~~~ ~~~~~~~~~~~ ~~~~~~~~~~~ ~~~~~~~~~~~  ~~~~~~~~~~~ 
~~~~~~~~~~~~~~~~~~~~~~\left.
             \left( \Dmaxt \xkspius \right)^2, 
              \left( \Dmint \xkspius \right)^2 
       \right\}.
\end{align*}
This averaging is legitimate because on the RHS the paths $\pi$ being considered in the averaging
all have the same values for $ \Dmaxt \xkspius $ 
and for $ \Dmint \xkspius  $ as their computation does not involve updates in $\ubar$.

Since $\max\left\{\left| \Dmint \xkspius \right|, \left| \Dmaxt \xkspius \right| \right\} ~
\le ~ \max\left\{\left| \Dmint \xkspiphis \right|, \left| \Dmaxt \xkspiphis \right| \right\}$,
and as \\
$ \Delta  \xkspis \in \left[\Dmint \xkspiphis ~,~ \Dmaxt \xkspiphis \right]$,
we have $\left| \Dmint \xkspius \right|, \left| \Dmaxt \xkspius \right| 
~\le ~ \left|\Delta  \xkspis \right| +  \left(  \Dmaxt \xkspiphis - \Dmint \xkspiphis \right)$;
it follows that
$$\left(\Dmint \xkspius \right)^2, \left(\Dmaxt \xkspius \right)^2 \le 2\left( \Delta  \xkspis \right)^2 + 2 \left(  \Dmaxt \xkspiphis - \Dmint \xkspiphis \right)^2.$$
Thus
\begin{align*}
& \expect{ \max_{t}\left\{ \left( \Dmaxt \xkpupiphiu - \Dmint \xkpupiphiu \right)^2 \right\}~ \Big| ~ k'_u ~\text{and}~ k_u ~\text{are substitutable}~ } \\
& ~~~~~~~~\le ~ \frac {2q\Lres^2} {(n-q)\G^2} 
\sum_{ s\in  [u-2q,  u+q] \setminus{ \{u\} } }  
  \max_t \left[
             2\left(  \Dmaxt \xkspiphis - \Dmint \xkspiphis \right)^2 
+ 2\left( \Delta  \xkspis \right)^2
     \right].
\end{align*}

Now, we extend the expectation to all paths $\pi$.
Let  $\left(\FEu \right)^2$ denote the resulting expectation for level $u$:
\begin{align*}
\left(\FEu \right)^2 ~=~ \expectpi{\max_{t} \left( \Dmaxt \xkupiphiu - \Dmint \xkupiphiu \right)^2}.
\end{align*}
Also, let $\left( \Es \right)^2 = \expectpi{\left( \Delta  \xkspis \right)^2}$.
Then, 
\begin{align}
\label{eqn:redib-FE}
(\FEu )^2 ~\le 
~\frac {4q\Lres^2} {(n-q)\G^2} 
	\sum_{ s\in  [u-2q, u+q] \setminus{\{u\}  }}
     \left[
             \left(  \FEs   \right)^2 + \left(\Es \right)^2
     \right]
~\le 
~\frac {4q\Lres^2} {(n-q)\G^2} 
\sum_{ s\in  [u-2q, u+2q] \setminus{\{u\} } }
     \left[
             \left(  \FEs   \right)^2 + \left(\Es \right)^2
     \right].
\end{align}

Let $\nu = \frac  {16q^2\Lres^2} {(n-q)\G^2} $ and choose $\G$ so that $\nu <1 $.
Then we can repeatedly redistribute the terms $(\FEs)^2$ to level $s$ and apply the upper bound in~\eqref{eqn:redib-FE} recursively,
ad infinitum.
We see that level $s$ receives a total charge of at most 
\begin{equation}
\label{eqn:rec-FE-defn}
\nu \cdot \left(\Es \right)^2 \frac {1}{1 - \nu}.
\end{equation}

\subsection{Rate of Convergence}
\eqref {eqn:F-prog} gives in expectation 
\begin{equation*}
F(\ptone) - F(\pt) ~ \ge ~\frac 12~\sum_{k_t=1}^n~\left[\hW_{k_t}(\gkt,x_{k_t}^{t-1},\G,\Psik) + \frac 18 \G \left( \Delta x_{k_t}^t \right)^2 - \frac {1}{\G} (\gkt - \tgkt)^2\right].
\end{equation*}
We want to bound $\frac {1}{\G} (\gkt - \tgkt)^2$ by $\frac 18 \G \left( \Delta x_{k_t}^t \right)^2$ in expectation in an
amortized sense.
By the discussion above, $\expectpi{\frac{1}{\G} (\gkt - \tgkt)^2}$ is bounded by 
\begin{equation}
\label{eqn:g-bound-as-FE-E}
\sum_{\stackrel{ s\in [t-2q, t+2q]} {s\ne t} }\frac {\nu}{4q} \left[\G\cdot \left(\FEs\right)^2
~+~ \G\cdot \left(\Es \right)^2 \right].
\end{equation}

Let $F(t)= \expect{F(\pt)}$.
Considering the recursively unwound form of~\eqref{eqn:redib-FE},
i.e., with the RHS having only terms of the form $\left(\Es \right)^2$,
shows that some of the cost in \eqref{eqn:g-bound-as-FE-E} at time $t$
has already been paid for (namely for the terms with $s<t$),
while some of the progress will be paid for in the future (for terms with $s>t$).
Accordingly we introduce two functions $\Ap(t)$ and $\Am(t)$;
$\Ap(t)$ will be the exact value of the progress already achieved that is needed in the future,
and $\Am(t)$ is the exact value of the already desired progress that will be paid for in the future.
$\Ap(t), \Am(t) \ge 0$ for all $t$, $\Ap(0) = \Am(0) = 0$ and $\Ap(T) = \Am(T)=0$,
for a run that lasts for exactly $T$ iterations.
Let $H(t) = F(t) + (1 + \gamma) \Ap(t) - \Am(t)$, for a suitable parameter $\gamma > 0$.
We will show that 
\begin{equation}
\label{eqn:stoc-prog}
H(\tone) - H(t) \ge \frac{1}{n} \left[ \frac 12 \sum_{k_t=1}^n \hW_{k_t}(\gkt,x_{k_t}^{t-1},\G,\Psik)
+ \frac 12 (1 + \gamma) \Ap(\tone) \right].
\end{equation}
The standard coordinate descent convergence bound then applies to $H$, namely:

\begin{theorem}
\label{thm:stoch-conv}
If \eqref{eqn:stoc-prog} holds, then:

\noindent
\emph{(i)}~If $F$ is strongly convex with parameter $\muF$,
	and $f$ has strongly convex parameter $\muf$,
	then
\begin{equation}
\label{eqn:stoc-conv-str-convex}
	H(t) ~\le~ \left[ 1 - \frac {1}{2n} \cdot \frac{\muF} {\muF + \G - \muf} \right]^{t} \cdot H(0).
\end{equation}
	
\noindent
\emph{(ii)}~Now suppose that $F$ is convex.
	Let $R$ be the radius of the level set for $\xc$. 
	Then, for $t \ge 0$,
\begin{equation}
	\label{eqn:stoc-conv-wkly-convex}
	H(t) ~\le~ \frac{H(0)} {1 + \frac{1}{2n} \cdot \min \left\{ 1 ~,~ \frac{H(0)} {\G ~\calR ^2} \right\} \cdot ~ t}.
\end{equation}
\end{theorem}
As $F(T)=H(T)$ and $F(0)=H(0)$,
\eqref{eqn:stoc-conv-str-convex} and \eqref{eqn:stoc-conv-wkly-convex}
continue to hold with $F(T)$ replacing $H(T)$ and $F(0)$ replacing $H(0)$.

In order to determine the right value for $\gamma$ we will need to derive some bounds
on $\Ap(t)$ and $\Am(t)$.
By \eqref{eqn:rec-FE-defn}, the total demand for $\left(\Es\right)^2$
is at most $\frac{\nu}{1-\nu}\cdot \G \cdot \left(\Es\right)^2$.
We write the $\left(\Es\right)^2$ demand at time $t$ as $m(s,t)\cdot \G \cdot \left(\Es\right)^2$,
where $m(s,t)$ is a suitable function.
We will derive tight bounds on $m(s,t)$ based on \eqref{eqn:g-bound-as-FE-E}, assuming that $s$ and $t$
range unboundedly over $(-\infty,+\infty)$, which are upper bounds for their actual values in a $T$ iteration computation.
As \eqref{eqn:g-bound-as-FE-E} is symmetric about $t$, it follows that $m(t,t+i) = m(t,t-i)$, for all $i\ge 1$.
Further, the bound on the total demand for $\left(\Es\right)^2$, 
namely $\frac{\nu}{1 - \nu}\cdot \G\cdot \left(\Es\right)^2$, continues to apply
as the earlier argument did not assume any particular bounds on $s$ and $t$.
We will show the following bounds:

\begin{lemma}
\label{lem:m-bounds}
Suppose that $\nu \le 1/3$. Then,
\begin{enumerate}
	\item[(A)] $m(t,t) \le \frac {\nu}{4q} \cdot \frac{\nu}{1 - \nu}$;
	\item[(B)] for any $k\ge 0$ and $1\le i \le 2q$, $m(t, t +k \cdot 2q+i) \le \frac {1}{4q} \cdot \frac{\nu^{k+1}}{1 - \nu}$;
	\item[(C)] for any $i\ge 1$, $m(t,t+i) \ge m(t,t+i+1)$;
	\item[(D)] for any $i\ge 1$, $m(t,t+2q+i) \le \nu \cdot m(t,t+i)$; and
	\item[(E)] $m(t,t) \ge m(t,t+2q+1)$.
\end{enumerate}
\end{lemma}

\begin{pf}
To derive the bound, consider a fixed $t$ and the demand for it generated at each level $s$.
By \eqref{eqn:g-bound-as-FE-E}, this demand comes in the form of two terms,
bounded by $\frac {\nu}{4q} \cdot \G \cdot \left( \Et \right)^2$, for $t< s \le t+2q$,
and $\sum_{ u\in[s-2q,s+2q]\setminus{\{s\} } } \frac {\nu}{4q} \cdot \G \cdot \left( \FEu \right)^2$, 
for all $s \ge t$, with the terms in the sum being redistributed recursively.

Let us view \eqref{eqn:g-bound-as-FE-E} as an iterative process, defined as follows.
First, define $m^1$ as below:
$$
m^1(t,s) ~:=~
\begin{dcases}
\frac {\nu}{4q},& \text{ if }t < s \le t+2q\\
&\\
0,&\text{ otherwise,}
\end{dcases}
$$
and for $j\ge 1$, define $m^{j+1}$ recursively as
$$
m^{j+1}(t,s) ~:=~
\begin{dcases}
\frac {\nu}{4q} ~+~ \sum_{u\in[s-2q,s+2q] \setminus{\{s\} } } \frac {\nu}{4q}\cdot m^j(t,u),& \text{ if }t < s \le t+2q\\
&\\
\sum_{u\in[s-2q,s+2q] \setminus{\{s\} } } \frac {\nu}{4q}\cdot m^j(t,u),&\text{ otherwise.}
\end{dcases}
$$
It is easy to verify by induction that $m^j(t,s)$ satisfies the bounds (A) and (B) for all $j$.
Clearly, $\lim_{j\ra \infty} m^j(t,s)$ is an upper bound on $m(t,s)$, and so it follows that $m(t,s)$
satisfies these bounds too.

To prove bounds (C), (D) and (E),
we again look at $m^j(t,s)$ and verify by induction that it satisfies the bounds,
and thus the bounds hold for the function $m$ too.

We will use following bound (F), which can be proved easily by induction:
for $j\ge 1$, $m^{j+1}(t,s) \ge m^j(t,s)$.

To prove (C), we separate into three cases.
If $i > 2q$,
\begin{align*}
m^{j+1}(t, t+i) - m^{j+1}(t, t+i+1) 
& = ~ \frac {\nu}{4q} \left[ m^{j}(t,t+i-2q) -m^j(t,t+i) + m^j(t,t+i+1) - m^j(t,t+i+2q+1) \right]
 \\
& \ge ~\frac {\nu}{4q} \left[ m^{j}(t,t+i-1 ) -m^j(t,t+i) + m^j(t,t+i+1) - m^j(t,t+i+2) \right] ~\ge ~ 0,
\end{align*}
in which the inequality holds by induction hypothesis.
If $i = 2q$,
\begin{align*}
m^{j+1}(t, t+2q) - m^{j+1}(t, t+2q+1) 
& \ge ~ \frac {\nu}{4q} \left[1+ m^{j}(t,t) -m^j(t,t+2q) + m^j(t,t+2q+1) - m^j(t,t+4q+1) \right]\\
& \ge 0 ~~~~\text{if~} 1 \ge \frac {1}{4q}m(t,t+2q), ~\text{i.e., if}~1 \ge \frac {1}{4q}\cdot \frac{\nu}{1 - \nu};
\nu \le \frac 45~\text{suffices}.
\end{align*}
If $1\le i \le 2q$, we use the definition of $m^j$, bounds (A) and (B) and the assumption that $\nu \le 1/3$:
\begin{align*}
m^{j+1}(t, t+i) - m^{j+1}(t, t+i+1) 
& = ~ \frac {\nu}{4q} \left[ m^{j}(t,t+i-2q) -m^j(t,t+i) + m^j(t,t+i+1) - m^j(t,t+i+2q+1) \right]
\\
& \ge ~\frac {\nu}{4q} \left[ \frac{\nu}{4q} ~-~ \frac{\nu}{4q}\cdot \frac{1}{1-\nu} ~+~ \frac{\nu}{4q}
~-~ \frac{\nu}{4q}\cdot \frac{\nu}{1-\nu} \right]\\
& = ~ \left(\frac {\nu}{4q}\right)^2 \cdot \left( 2 - \frac{1+\nu}{1-\nu} \right) ~\ge~ 0.
\end{align*}

Bound (D) is straightforward given bounds (C) and (F):
\begin{align*}
m^{j+1}(t,t+2q+i) ~=~ \frac {\nu}{4q} \sum_{u\in[t+i,t+4q+i] \setminus{\{t+2q+i\} } } 
m^j(t,u)  ~\le ~ \nu \cdot m^j(t,t+i)  ~\le ~ \nu \cdot m^{j+1}(t,t+i).
\end{align*}

Bound (E) is also straightforward given bound (C):
\begin{align*}
m^{j+1}(t,t) &~=~ 
2~\sum_{u\in [t+1,t+2q]} \frac {\nu}{4q}\cdot m^j(t,u)\\
&~\ge~ \sum_{u\in [t+1,t+2q]} \frac {\nu}{4q}\cdot m^j(t,u) ~+~
\sum_{u\in [t+2q+2,t+4q+1]} \frac {\nu}{4q}\cdot m^j(t,u) ~=~ m^{j+1}(t,t+2q+1).
\end{align*}
\end{pf}

\begin{cor}
\label{cor:Ap-bound}
$\Ap(t) \le \sum_{1 \le s \le t} \frac{2q}{1 - \nu}\cdot m(s,t+1) \cdot \G  \cdot \left(\Es\right)^2$.
\end{cor}
\begin{pf}
By bounds (C) and (D) of Lemma~\ref{lem:m-bounds},
\begin{align*}
\Ap(t) ~=~ \sum_{s\le t}\sum_{v>t} m(s,v) \cdot \G  \cdot \left(\Es\right)^2
&~\le~ \sum_{s \le t} 2q\cdot m(s,t+1) \left[1 + \nu +\nu^2 + \ldots \right] \cdot \G \cdot \left(\Es\right)^2\\
&~\le~  \sum_{s\le t} \frac{2q}{1 - \nu}\cdot m(s,t+1)\cdot \G \cdot \left(\Es\right)^2.
\end{align*}
\end{pf}

Next we show~\eqref{eqn:stoc-prog}.
\begin{lemma}
\label{lem:stoc-prog}
Suppose that $\nu \le \frac {1}{10}$,
$ \gamma = \frac {q}{n (1 - \nu) - q}$,
and $q \le \frac  {1 - \nu}{5} n$.
Then
$$H(\tone) - H(t) \ge \frac{1}{2n} \left[  \sum_{k_t=1}^n \hW_{k_t}(\gkt,x_{k_t}^{t-1},\G,\Psik) + \Ap(\tone) \right].$$
\end{lemma}

\begin{pf}
\begin{align*} 
H(\tone) - H(t) &~\ge~ \left[F(\tone) + (1 + \gamma)\Ap(\tone) - \Am(\tone)\right]
~-~ \left[F(t) + (1 +\gamma)\Ap(t) - \Am(t)\right] \\
&~\ge~  \frac{1}{2n} \expect{\sum_{k=1}^n\hWk(\gkt,\pktm,\G,\Psik) } 
~+~ \frac {\G}{8} \left(\Et\right)^2
~-~ \frac {1}{\G} \expect{ \left( \gkt - \tgkt \right)^2}\\
&~~~~~~~+~ (1+\gamma)\cdot \left[ \Ap(\tone) - \Ap(t)\right] ~-~ \left[ \Am(\tone) - \Am(t)  \right]\\
&~\ge~ \frac{1}{2n} \expect{\sum_{k=1}^n\hWk(\gkt,\pktm,\G,\Psik) } 
~+~ \frac {\G}{8} \left(\Et\right)^2
~-~ \frac {1}{\G} \expect{ \left( \gkt - \tgkt \right)^2}\\
&~~~~~~~+~ (1+\gamma) \left[ \sum_{s=1}^{t-1} m(s,t)\cdot \G\cdot \left(\Es\right)^2
~-~ \sum_{v>t} m(t,v)\cdot \G\cdot \left( \Et \right)^2\right]\\
&~~~~~~~+~ \sum_{v=t+1}^{T} m(v,t)\cdot \G\cdot (E_v^x)^2
 ~-~ \sum_{s=1}^{t-1}m(t,s)\cdot \G \cdot \left( \Et \right)^2.
\end{align*}
Now, note that
$$
\frac {1}{\G} \expect{ \left( \gkt - \tgkt \right)^2} ~\le~ \sum_{s=1}^T m(s,t)\cdot \G \cdot \left(\Es\right)^2,
$$
and hence
\begin{align*}
H(\tone) - H(t)
&~\ge~ \frac{1}{2n} \expect{\sum_{k=1}^n\hWk(\gkt,\pktm,\G,\Psik) }
~+~ \gamma \cdot \sum_{s=1}^{t-1} m(s,t)\cdot \G\cdot \left(\Es\right)^2 \\
&~~~~~~~~+~ \left(\frac 18 - (1+\gamma)\cdot \sum_{v>t} m(t,v) ~-~ m(t,t)  ~-~\sum_{s=1}^{t-1}m(t,s) \right)~\G\cdot \left( \Et \right)^2 \\
&~\ge ~ \frac{1}{2n} \expect{\sum_{k=1}^n\hWk(\gkt,\pktm,\G,\Psik) } + \gamma \cdot \frac {1 - \nu} {2q} \cdot \Ap(t-1) ~~~~\comm{by Corollary~\ref{cor:Ap-bound}} \\
&~~~~~~~~~~~~+~ \left(\frac 18 ~-~ \sum_{v \ge 1} m(t,v) ~-~ \gamma \sum_{v > t} m(t,v)\right)~\G\cdot \left( \Et \right)^2.
\end{align*}
	
We want $\frac{\gamma (1 - \nu)} {2q} \ge \frac {1}{2n}(1 + \gamma)$; i.e., $\gamma\cdot n (1 - \nu) \ge q + q\gamma$;
$ \gamma \ge \frac {q}{n (1 - \nu) - q}$ suffices.

By symmetry, $\sum_{v>t} m(t,v) \le \frac 12 \cdot \frac {\nu}{1 - \nu}$.
Thus, it suffices that $\frac {1}{8} \ge \frac {\nu}{1 - \nu} \left(1 + \frac {\gamma}{2}\right)$;
with $\nu \le \frac {1}{10}$,
it suffices that $ \frac 18 \ge \frac 1{9} (1 + \frac 12 \gamma)$, 
which holds if $\gamma \le \frac {1}{4} $, i.e., if $n(1 - \nu) - q \ge 4q$, or $q \le \frac{1 - \nu}{5} n$.
\end{pf}

Recall that $\nu = \frac  {16 q^2\Lres^2} {(n-q)\G^2}$. 
We have shown that \eqref{eqn:stoc-prog}
holds if $\nu \le \frac 1{10}$, 
and hence so do the conclusions of Theorem~\ref{thm:stoch-conv}.
This holds if $q \le \frac{\G \sqrt{n-q}} {4 \sqrt {10} \cdot \Lres}$ and $q \le \frac {9}{50} n$.

\paragraph{Fully general non-consistent coordinates}
Here we address the possibility that the start and commit times of a single coordinate need
not be consistent.
We address this by providing each coordinate with a \emph{virtual} start time, where the virtual start times
are consistent with the commit times.
For each coordinate $x_k$, the virtual start times are simply the actual start times of $x_k$ in sorted order,
with the $i$-th such time becoming the virtual start time for the update instance of
$x_k$ with the $i$-th commit time in sorted order.
The effect, for an update whose virtual start time come before its actual start time, is to
increase the number of earlier updates that would interfere with it, but only up to at most $2q$.
To see this, suppose the update at start time $t$ receives the virtual start time $s<t$.
Then the updates at times $s$ and $t$ could interfere with each other, and so $s\ge t-q$.
The only new updates that could
interfere with the update with real start time $t$ are the up to $q$ updates with start
times in $[s-q,s-1]$, the updates before time $s$ that might interfere with the update with real
start time $s$.

So it suffices to replace $q$ with $2q$ in the previous bounds.

\newpage

\section{Parallel Asynchronous Coordinate Descent (PACD)}\label{sect:appendix-parallel-ACD}

\subsection{Proof of Inequality \eqref{eqn:target-PACD}}

First of all, by \eqref{eqn:F-prog},
\begin{align*}
&\sum_{i=t-2r+1}^t \left[ H(i-1) - H(i) \right]\\
&~\ge~ \sum_{i=t-2r+1}^t \left[ ~\frac 12 \hW_{k_i}(g_{k_i}^i , x_{k_i}^{i-1} , \G , \Psi_{k_i})
~+~ \frac 18 ~\G \left(\Delta x_{k_i}^i\right)^2
~-~ \frac 1 \G (g_{k_i}^i - \tg_{k_i}^i)^2~ \right] + A(t-2r) - A(t)
\end{align*}

Recall that for each time $i$,  $\alpha(i)$ denotes the time of the latest update to coordinate $k_i$ strictly before time $i$,
and is $0$ if no such update exists.
Note that $i - \alpha(i) \le r$ by Assumption \ref{assume:PACD-frequent}.
Then by Lemma \ref{lem:W-shift},
$$
\hW_{k_i}(g_{k_i}^i , x_{k_i}^{i-1} , \G , \Psi_{k_i}) ~\ge~ 
\frac 1r \sum_{j=\max\{\alpha(i)+1,t-2r+1\}}^i \left[ ~\frac 23 \cdot \hW_{k_i}(g_{k_i}^j , x_{k_i}^{j-1} , \G , \Psi_{k_i})
~-~ \frac{4}{3\G} \cdot \left(g_{k_i}^j - g_{k_i}^i\right)^2 ~\right],
$$
and hence
\begin{align*}
&\sum_{i=t-2r+1}^t \left[ H(i-1) - H(i) \right]\\
&~\ge~ \sum_{i=t-2r+1}^t ~\frac 1{3r} 
\sum_{j=\max\{\alpha(i)+1,t-2r+1\}}^i ~\hW_{k_i}(g_{k_i}^j , x_{k_i}^{j-1} , \G , \Psi_{k_i})
~-~ \frac{2}{3\G r} \sum_{i=t-2r+1}^t ~ \sum_{j=\max\{\alpha(i)+1,t-2r+1\}}^i ~ \left(g_{k_i}^j - g_{k_i}^i\right)^2\\
&~~~~~~~~~~~~~~+~\frac 18 \sum_{i=t-2r+1}^t ~\G \left(\Delta x_{k_i}^i\right)^2
~-~ \sum_{i=t-2r+1}^t \frac 1 \G (g_{k_i}^i - \tg_{k_i}^i)^2 ~~+~~ A(t-2r) - A(t)\\
&~\ge~ \sum_{j=t-2r+1}^{t-r} ~\frac {1}{3r} \sum_{k=1}^n \hWk(g_k^j,x_k^{j-1},\G,\Psi_k)
~+~ \left[\frac 18 \sum_{i=t-2r+1}^t ~\G \left(\Delta x_{k_i}^i\right)^2 + A(t-2r) - A(t)\right]\\
&~~~~~~~~~~~~~~-~ \frac{2}{3\G r} \sum_{i=t-2r+1}^t ~ \sum_{j=\max\{\alpha(i)+1,t-2r+1\}}^i ~ \left(g_{k_i}^j - g_{k_i}^i\right)^2
~-~ \sum_{i=t-2r+1}^t \frac 1 \G (g_{k_i}^i - \tg_{k_i}^i)^2\\
&~\ge~ \sum_{j=t-2r+1}^{t-r} ~\left[\frac {1}{3r} \sum_{k=1}^n \hWk(g_k^j,x_k^{j-1},\G,\Psi_k) ~+~ \frac{1}{2q} \cdot A(j-1)\right]\\
&~~~~~~~~~~~~~+~\left[\frac 18 \sum_{i=t-2r+1}^t ~\G \left(\Delta x_{k_i}^i\right)^2
~-~ \frac{1}{2q} \sum_{j=t-2r+1}^{t-r} A(j-1)  ~+~ A(t-2r) - A(t)\right]\\
&~~~~~~~~~~~~~~~~~~~~~~-~ \frac{2}{3\G r} \sum_{i=t-2r+1}^t ~ \sum_{j=\max\{\alpha(i)+1,t-2r+1\}}^i ~ \left(g_{k_i}^j - g_{k_i}^i\right)^2
~-~ \sum_{i=t-2r+1}^t \frac 1 \G (g_{k_i}^i - \tg_{k_i}^i)^2.
\end{align*}

By a direct expansion using the definition of $A$,
\begin{align*}
&\frac 18 \sum_{i=t-2r+1}^t ~\G \left(\Delta x_{k_i}^i\right)^2
~-~ \frac{1}{2q} \sum_{j=t-2r+1}^{t-r} A(j-1)  ~+~ A(t-2r) - A(t)\\
& ~\ge~ \frac{3}{64} \sum_{i=t-2r-q}^{t-2r} \frac{i-t+2r+q}{q} ~\G\left(\Delta x_{k_i}^i\right)^2
~+~ \frac{1}{16} \sum_{i=t-2r+1}^t \G \left(\Delta x_{k_i}^i\right)^2\\
& ~\ge~ \frac{3\G}{64q}\left[ \sum_{i=t-2r-q}^{t-2r} \frac{i-t+2r+q}{q} ~\G\left(\Delta x_{k_i}^i\right)^2
~+~ \sum_{i=t-2r+1}^t \G \left(\Delta x_{k_i}^i\right)^2 \right],
\end{align*}
and we are done.

\subsection{Proving that $H$ is Decreasing}

For any $t\ge 1$, by \eqref{eqn:F-prog},
\begin{align*}
&~~~~~H(t-1) - H(t)\\
&~\ge~ \frac 12 \hW_{k_t} (g_{k_t}^t,x_{k_t}^{t-1},\G,\Psi_{k_t}) + \frac 18 (\Delta x_{k_t}^t)^2
- \frac 1\G (g_{k_t}^t - \tg_{k_t}^t)^2 + \frac{1}{16} \sum_{i=\max\{1,t-q\}}^{t-1} \frac{1}{q} \cdot \G (\Delta x_{k_i}^i)^2
- \frac{1}{16} (\Delta x_{k_t}^t)^2\\
&~\ge~ -\frac{(\Lmax)^2 q}{\G} \sum_{i=\max\{1,t-q\}}^{t-1} (\Delta x_{k_i}^i)^2
+ \frac{\G}{16q} \sum_{i=\max\{1,t-q\}}^{t-1} (\Delta x_{k_i}^i)^2.\comm{by \eqref{eqn:bound-g-minus-tg}}
\end{align*}
When $\frac{\G}{16q}~\ge~ \frac{(\Lmax)^2 q}{\G}$, or equivalently $\G ~\ge~ 4q \Lmax$, we are done.

\subsection{Enforcing Assumptions \ref{assume:PACD-q} and \ref{assume:PACD-frequent}}

To enforce Assumption~\ref{assume:PACD-q}, one can use two fetch-and-add global counters,
one for counting the number of updates that have been started, and the other counting the number of updates that have finished.
A core is allowed to initiate a new iteration only when the two counters differ by at most $q/2$.
Then, if the updates require a similar amount of computation and communication,
and assuming the asynchronous effects are not too variable,
it is plausible that there will not be much busy waiting for $q$ a small multiple of the number of cores.

Next, we describe one modest overhead method to enforce Assumption~\ref{assume:PACD-frequent} in the following setting:
the coordinates are distributed among the processors so that they have similar work loads and
each processor iterates over its coordinates cyclically.

We choose $\kstop = \floor{r/2n}$.
We maintain a global counter of the number of updates that have committed using the fetch and add operation.
We will be partitioning the updates into pseudo-rounds which have length at most $r$.

Each processor keeps a count of the number of updates it has committed since the current pseudo-round began
(we will explain how the start of a pseudo-round is detected shortly).
If this count reaches $\kstop$ times its number of variables it waits for the start of the next pseudo-round before
performing more updates.
Thus two updates of any variable must occur within $2(n-1)\floor{r/2n} +1 \le r$ updates of each other.
Further, in any round of $r$ updates, there can be at most $\kmax \doteq 2\kstop$ updates to any one variable.

To identify the start of a new pseudo-round we use 3 shared counters, whose roles rotate over cycles of 3 pseudo-rounds.
The current counter keeps track of the number of variables that have been 
updated at least once during the current pseudo-round, by means of fetch and add operations.
A new pseudo-round begins when the current counter attains the value $n$.
At that point, we switch to updating the second counter, and as each commit in this round is performed the
third counter is reset to $0$.

Note that if keeping counters for each processor's collection of variables is impractical, 
one could instead follow the above procedure for each variable separately.

\newpage

\section{Tatonnement Analysis for Complementary-CES Fisher Markets}
\label{app:tat-anal}

\newcommand{\alphak}{\alpha_k}
\newcommand{\alphaks}{\alpha_{k_s}}
\newcommand{\alphakt}{\alpha_{k_t}}
\newcommand{\alphaktau}{\alpha_{k_{\tau}}}
\newcommand{\dt}{\delta t}
\newcommand{\dtau}{\delta t_{\tau}}
\newcommand{\Dtau}{\Delta  t_{\tau}}
\newcommand{\Dsig}{\Delta t_{\sigma}}
\newcommand{\dsig}{\delta t_{\sigma}}
\renewcommand{\gktt}{g_{k_t}^t}
\newcommand{\gksigsig}{g_{k_{\sigma}}^{\sigma}}
\newcommand{\gksigtau}{g_{k_{\sigma}}^{\tau}}
\renewcommand{\gktautau}{g_{k_{\tau}}^{\tau}}
\newcommand{\tgktt}{\tilde{g}_{k_t}^t}
\newcommand{\tgksigtau}{\tilde{g}_{k_{\sigma}}^{\tau}}
\renewcommand{\tgktautau}{\tilde{g}_{k_{\tau}}^{\tau}}
\newcommand{\ks}{k_s}
\newcommand{\Gktt}{\G_{k_t}^t}
\newcommand{\Gktautau}{\G_{k_{\tau}}^{\tau}}
\newcommand{\Gksigtau}{\G_{k_{\sigma}}^{\tau}}
\newcommand{\Gktausig}{\G_{k_{\tau}}^{\sigma}}
\newcommand{\hWksig}{\widehat{W}_{k_{\sigma}}}
\newcommand{\Itau}{I_{\tau}}
\newcommand{\lktaut}{\ell_k(\tau,t)}
\newcommand{\lksigtau}{\ell_k(\sigma,\tau)}
\newcommand{\lksignutau}{\ell_k(\sigma,\nutau)}
\newcommand{\nutau}{\nu(\tau)}
\renewcommand{\ptone}{p^{t-1}}
\renewcommand{\pt}{p^{t}}
\renewcommand{\ptau}{p^{\tau}}
\newcommand{\pnutau}{p^{\nu(\tau)}}
\newcommand{\pnut}{p^{\nu(t)}}
\renewcommand{\ptp}{p^{t+1}}
\newcommand{\pkssig}{p_{k_s}^{\sigma}}
\newcommand{\pksigsig}{p_{k_\sigma}^{\sigma}}
\newcommand{\pkstau}{p_k^{\tau}}
\newcommand{\pkst}{p_k^t}
\renewcommand{\pkt}{p_{k_t}}
\renewcommand{\pktt}{p_{k_t}^t}
\renewcommand{\pkttm}{p_{k_t}^{t-}}
\renewcommand{\pkttau}{p_{k_t}^{\tau}}
\renewcommand{\pktautau}{p_{k_{\tau}}^{\tau}}
\newcommand{\pktausig}{p_{k_{\tau}}^{\sigma}}
\newcommand{\pksigsigm}{p_{k_{\sigma}}^{\sigma-}}
\renewcommand{\pktautaum}{p_{k_{\tau}}^{\tau-}}
\newcommand{\pti}{p^{t_i}}
\newcommand{\ptip}{p^{t_{i+1}}}
\newcommand{\Lksigkssig}{L_{k_{\sigma}k_s}^{[\sigma,\sigma+1]}}
\newcommand{\Lktaukstau}{L_{k_{\tau}k_s}^{[\tau,\tau+1]}}
\newcommand{\Lktkt}{L_{k_t,k}^{[t,t+1]}}
\newcommand{\Lktauktau}{L_{k_{\tau},k}^{[\tau,\tau+1]}}
\newcommand{\tprev}{t_{prev}}
\newcommand{\tmin}{t_{\min}}
\newcommand{\tsig}{t_{\sigma}}
\newcommand{\ttau}{t_{\tau}}
\renewcommand{\tt}{t_t}
\newcommand{\ti}{t_i}
\newcommand{\tip}{t_{i+1}}

In the tatonnement analysis, the control variables are prices, so it is more natural to use $p$
(instead of $x$ in the ACD setting) to denote their values.

\subsection{Comparison with the PACD Analysis}

At a high level, the analysis is similar to that for PACD.
Again we define a function $A(t)\ge 0$ and $H(t) = \phi(\pt) + A(t)$
(note that $\phi$ plays the role of $F$ here), with $A(0) = 0$.
The basic idea is to show that if there are a series of updates 
between times $t$ and $t'$, where $t-3\le t'\le t-2$,
then $H(t') - H(t) \ge H(t')\cdot \exp(-\varepsilon)$ for a suitable constant $\varepsilon > 0$.
This then implies that for any $t\ge 2$, $F(x^t) ~\le~ H(t) ~\le~ e^{-\varepsilon \cdot \Theta(t)}\cdot H(\pc) ~=~ F(x^\circ)$.

There are a few modifications to the PACD analysis.
First, the function $\phi(p)$ does not admit \emph{global} Lipschitz parameters.
We have to settle for using \emph{local} Lipschitz parameters,
and also using \emph{adaptive} step sizes in the update rule.
Second, we will use individual $L_{jk}$ parameters to bound gradient errors,
instead of using the lump parameter $\Lmax$ --- this is crucial for showing that $\lambda$,
the parameter in update rule \eqref{eq:async-tat-CES-rule}, can be as large as $\Theta(1)$.

Third, for tatonnement, the two types of gradient errors in the PACD setting can be merged into one.
To see why, note that in the market setting,
each seller of a good does \emph{not} observe the prices of the other goods,
but rather the excess demand of her own good which is a function of all prices,
while the excess demand is the same of the gradient of $\phi$ at current prices.
Thus, each update to a good uses an excess demand value of the good
which lies between the minimum and maximum values since the last update to the same good
--- this difference between the used value and the accurate value is the source of the second type of error in the PACD setting.
When translated to the coordinate descent setting,
the update uses a gradient value which is between the minimum and maximum values since the last update to the same coordinate,
and as we will see, it can be bounded in exactly the same way as the first type of error.
In contrast, in the PACD setting, due to \emph{inconsistent reads}, the used gradient value need not lie between
the minimum and maximum gradient values since the last update to the same coordinate.

We note that $\phi$ is a smooth function so $\Psik \equiv 0$,
and consequently we omit it from the arguments of the function $\hWk$.

\subsection{Analysis}

\newcommand{\deltau}{\delta_\tau}
\newcommand{\delnu}{\delta_\nu}
\renewcommand{\pktaunum}{p_{\ktau}^{\nu-}}
\newcommand{\gkmaxtau}{g_k^{\max,\tau}}
\newcommand{\gkmintau}{g_k^{\min,\tau}}
\newcommand{\gktaumaxtau}{g_{\ktau}^{\max,\tau}}
\newcommand{\gktaumintau}{g_{\ktau}^{\min,\tau}}
\newcommand{\numofupdatek}{U_k}

\paragraph{Notations, Definitions and Two Lemmas.}
Recall update rule \eqref{eq:async-tat-CES-rule}. Note that, for the purposes of our analysis,
for each update we now need to know the elapsed time since the previous update to the same coordinate,
or since time $0$ if it is the first update to that coordinate; for the update at time $\tau$, we denote this by $\Dtau$.
We let $\alphak(\tau)$ denote the time of the most recent update to $\pk$
\emph{strictly} before time $\tau$, or time $0$ if there is no previous update to this price.
We let $\alpha(\tau)$ denote the time of the most recent update to \emph{any} coordinate
\emph{strictly} before time $\tau$, or time $0$ if there is no previous update to this price.
And we let $\deltau ~:=~ \tau - \alpha(\tau)$.

Recall that the coordinates being updated are prices $\pj$.
We let $\pkttm$ denote the value of coordinate $\pkt$ right before it is updated at time $t$.
In our analysis, when we write $\sum_{\tau\in I}$, where $I$ is some time interval,
the summation is summing over all updates occurred in the time interval $I$.

For each update $\tau$, let $\gkmaxtau$ and $\gkmintau$ denote
the maximum and minimum of accurate gradient values along coordinate $k$ in the time interval $(\alpha_{\ktau}(\tau),\tau)$.

We need the following two lemmas. The first lemma is a variant of Lemma~\ref{lem:discrete-improvement-of-F},
which takes account of the time intervals.
The second lemma can be derived easily from the Power-Mean inequality.

\begin{lemma}\label{lem:discrete-improvement-of-F-with-time}
Suppose there is an update to coordinate $j$ at time $t$ according to rule \eqref{eq:async-tat-CES-rule},
and suppose that $\G_j \ge L_j$. Let $\tau = \alpha(t)$.
Let $g_j = \nabla_j f(\ptau)$ and $\tg_j = \nabla_j f(\tilde{p})$.
Then
\begin{align*}
F(\ptau) - F(\pt) &~\geq~ \frac{\G_j}{2}\frac{ (\Delta p_{j}^t)^2} {\Dt} - |g_j - \tg_j|\cdot |\Delta p_j^t|\\
\text{and}~~~~~~~~~~~~~~~~~~~~~~~
F(\ptau) - F(\pt) &~\geq~ \hW(g_j,\ptau_j,\G_j,\Psi_j)\cdot \Dt - \frac{1}{\G_j}(g_j - \tg_j)^2 \cdot \Dt.~~~~~~~~~~~~~~~~~~~~~~~~~
\end{align*}
\end{lemma}
This is proved in exactly the same way as  Lemma~\ref{lem:discrete-improvement-of-F}.
We can then deduce the following variant of \eqref{eqn:F-prog}.
\begin{equation}
\label{lem:F-progress-for-tat}
F(\ptau) - F(\pt)  ~\ge~ \frac 12 \hW(g_j,\ptau_j,\G_j,\Psi_j)\cdot \Dt ~+~ \frac{\G_j}{8}\frac{ (\Delta p_{j}^t)^2} {\Dt}
~-~ \frac {1}{\G_j}(g_j - \tg_j)^2 \cdot \Dt.
\end{equation}

\begin{lemma}\label{lem:power-mean-ine}
Suppose that $w_1,w_2,\cdots,w_\ell$ and $y_1,y_2,\cdots,y_\ell$ are non-negative numbers. Then
$$
\left(\sum_{j=1}^\ell w_j x_j\right)^2 ~~\le~~ \left(\sum_{j=1}^\ell w_j\right)\left( \sum_{j=1}^\ell w_j \cdot (x_j)^2 \right).
$$
\end{lemma}

\paragraph{Analysis Details.}
We define the function $A(t)$ as follows.
\[
A(t) = \frac{1}{3} \sum_{\tau\in (t-1,t]} ~~
\sum_{k\neq \ktau}~
(2-\min\{1,\numofupdatek(\tau,t)\}) \cdot 
\Lktauktau \cdot \frac{p_k^\tau}{p_{\ktau}^{\tau-}} \cdot \frac {\left( \Delta \pktautau \right)^2} {\Delta \ttau},
\]
where $\numofupdatek(\tau,t)$ denotes the number of updates to coordinate $k$ in the time interval $(\tau,t]$,
and each $L$ parameter of the format $L_{jk}^{[\tau_a,\tau_b]}$ is an upper bound
on the Lipschitz gradient parameter $L_{jk}$ of the function $\phi$,
within a rectangular hull of those prices which might appear in the time interval $[\tau_a,\tau_b]$ only.

For some $t\ge 2$ at which there is an update,
we let $t_a$ denote the time of the latest update \emph{strictly} before time $(t-2)$;
we let $t_a=0$ if no such update exists.
We let $t_b$ denote the time of the earliest update in the time interval $[t-1,t]$.
\begin{align}
&\sum_{\tau\in (t_a,t]} \left[H(\alpha(\tau)) - H(\tau)\right]\nonumber\\
&\ge~ \sum_{\tau\in (t_a,t]} \left[
~\frac 12 \hWktau(\gktautau,\pktautaum,\Gktautau) \cdot \Dtau ~+~ \frac {\Gktautau}8 \cdot \frac{(\Delta p_{\ktau}^\tau)^2}{\Dtau}
~-~\frac{1}{\Gktautau} (\gktautau - \tgktautau)^2\cdot \Dtau
~\right] ~+~ A(t_a) - A(t)\nonumber\\
&~\stackrel{(*)}{\ge}~ \sum_{\tau\in (t_a,t]} \left[
~\frac 12 \hWktau(\gktautau,\pktautaum,\Gktautau)\cdot \Dtau ~+~ \frac {\Gktautau}8 \cdot \frac{(\Delta p_{\ktau}^\tau)^2}{\Dtau}
~-~\frac{1}{\Gktautau} (\gktaumaxtau - \gktaumintau)^2\cdot \Dtau
~\right] ~+~ A(t_a) - A(t).\label{eq:tat-target-first}
\end{align}
Inequality $(*)$ holds because in the tatonnement setting,
both the accurate gradient (excess demand) $\gktautau$ and
the inaccurate gradient (excess demand) $\tgktautau$
must lie in between $\gktaumaxtau$ and $\gktaumintau$.
Note that in the PACD setting which allows inconsistent read, this does \emph{not} hold.\footnote{In PACD,
cores retrieve outdated prices and use them to compute gradients,
while in the tatonnement setting, sellers observe outdated gradients (excess demands) directly.}

By Lemma \ref{lem:W-shift}, for each $\tau\in (t_a,t]$,
\begin{align}
&\hWktau(\gktautau,\pktautaum,\Gktautau)\cdot \Dtau\nonumber\\
&~\ge~ \sum_{\nu\in (\max\{t_a,\alpha_{\ktau}(\tau)\} , \tau ]} ~ \left[~
\frac 23 \hWktau(\gktaunu,\pktaunum,\Gktautau) \cdot \delnu
~-~ \frac 43 \cdot \frac{1}{\Gktautau} \cdot \left( \gktaunu - \gktautau \right)^2\cdot \delnu ~\right]\nonumber\\
&~\ge~ \frac 23 \left(\sum_{\nu\in (\max\{t_a,\alpha_{\ktau}(\tau)\} , \tau ]} ~ \hWktau(\gktaunu,\pktaunum,\Gktautau) \cdot \delnu\right)
~-~ \frac 43 \cdot \frac{1}{\Gktautau}\cdot ( \gktaumaxtau - \gktaumintau )^2\cdot \Dtau.\label{eq:tat-target-second}
\end{align}

For any $k$ and any time $\nu$, we let $\G_k^\nu$ denote the step size used by the update to coordinate $k$ on or after time $\nu$.
Combining \eqref{eq:tat-target-first} and \eqref{eq:tat-target-second} yields
\begin{align*}
&\sum_{\tau\in (t_a,t]} \left[H(\alpha(\tau)) - H(\tau)\right]\nonumber\\
&\ge~ \frac 13 \sum_{\tau\in (t_a,t]} ~ \sum_{\nu\in (\max\{t_a,\alpha_{\ktau}(\tau)\} , \tau ]}
~\hWktau(\gktaunu,\pktaunum,\Gktautau) \cdot \delnu\\
&~~~~~~~~~~~
~+~ \left(\frac 18 \sum_{\tau\in (t_a,t]} \Gktautau \frac{(\Delta p_{\ktau}^\tau)^2}{\Dtau} + A(t_a) - A(t)\right)
~-~ \frac {5}{3} \sum_{\tau\in (t_a,t]} \frac{1}{\Gktautau} \cdot ( \gktaumaxtau - \gktaumintau )^2\cdot \Dtau\\
&\ge~ \frac 13 \sum_{\nu\in (t_a,t_b]} \delnu\cdot \sum_{k=1}^n \hWk(g_k^{\nu},p_k^{\nu-},\G_k^\nu)\\
&~~~~~~~~~~~
~+~ \left(\frac 18 \sum_{\tau\in (t_a,t]} \Gktautau \frac{(\Delta p_{\ktau}^\tau)^2}{\Dtau} + A(t_a) - A(t)\right)
~-~ \frac {5}{3} \sum_{\tau\in (t_a,t]} \frac{1}{\Gktautau} \cdot ( \gktaumaxtau - \gktaumintau )^2\cdot \Dtau\\
&=~ \sum_{\nu\in (t_a,t_b]} \delnu\cdot \left(\frac 13~\sum_{k=1}^n \hWk(g_k^{\nu},p_k^{\nu-},\G_k^\nu) ~+~ \frac 1{12} \cdot  A(\alpha(\nu))\right)\\
&~~~~~
~+~ \left(\frac 18 \sum_{\tau\in (t_a,t]} \Gktautau \frac{(\Delta p_{\ktau}^\tau)^2}{\Dtau} - \frac 1{12} \sum_{\nu\in (t_a,t_b]} \delnu\cdot A(\alpha(\nu)) + A(t_a) - A(t)\right)
~-~ \frac {5}{3} \sum_{\tau\in (t_a,t]} \frac{1}{\Gktautau} \cdot ( \gktaumaxtau - \gktaumintau )^2\cdot \Dtau.
\end{align*}

In \cite{CCD2013}, it was proved that the function $\phi$ is strongly convex, and that the maximum $\G$ value throughout the tatonnement
is upper bounded by a finite constant which depends on the starting price $p^\circ$.\footnote{Their
argument concerned the synchronous setting, but it can be reused without change for the asynchronous setting.}
We denote the finite upper bound on all $\G$'s by $\bG$, and the strong convexity parameter of $\phi$ by $\mu_\phi$.
We let $\varepsilon ~:=~ \mu_\phi / (\mu_\phi + \bG)$.
Then by Lemma \ref{lem:good-progress},
$$\sum_{k=1}^n \hWk(g_k^{\nu},p_k^{\nu-},\G_k^\nu) ~\ge~ \varepsilon \cdot F(\alpha(\nu)).$$
By replacing $\varepsilon$ with $\min\{\varepsilon/3~,~1/12\}$, we have
\begin{align}
&\sum_{\tau\in (t_a,t]} \left[H(\alpha(\tau)) - H(\tau)\right]\nonumber\\
&~\ge~ \sum_{\nu\in (t_a,t_b]} \varepsilon \cdot \delnu\cdot H(\alpha(\nu))
~+~ \left(\frac 18 \sum_{\tau\in (t_a,t]} \Gktautau \frac{(\Delta p_{\ktau}^\tau)^2}{\Dtau} - \frac 1{12} \sum_{\nu\in (t_a,t_b]} \delnu\cdot A(\alpha(\nu)) ~+~ A(t_a) - A(t)\right)\nonumber\\
&~~~~~~~~~~~~~~~~~~~~~~~~~~~~~~~~~~~~~~~~~~~~
~-~ \frac{5}{3} \sum_{\tau\in (t_a,t]} \frac{1}{\Gktautau} \cdot ( \gktaumaxtau - \gktaumintau )^2\cdot \Dtau.\label{eq:tat-target-final}
\end{align}

We are going to prove that the final two terms, in sum, are non-negative.
Then the above inequality is equivalent to saying that the progress from time $t_a$ to time $t$
is at least the progress achieved by the following process: over the time interval $(t_a,t_b]$,
which is divided into subintervals indexed by $\nu$, of respective length $\delnu$,
such that in each interval the progress is reducing $H$ by a factor of at least $\exp(-\varepsilon \cdot \delnu)$.
Since the length of the time interval $(t_a,t_b]$ is at least $1$,
$$
H(t) ~\le~ H(t_a) \cdot \exp\left( - \varepsilon\cdot \sum_{\nu\in (t_a,t_b]} \delnu\right)
~\le~ H(t_a) \cdot \exp(-\varepsilon).
$$
Iterating the above inequality from time $t$ down to at least time $2$
yields
$$F(p^t) ~\le~ H(t) ~\le~ \exp(-\varepsilon\cdot(t/3-1)) \cdot H(0) ~=~ \exp(-\varepsilon\cdot(t/3-1)) \cdot F(p^\circ),$$
as desired.

We handle the remaining tasks in the next two subsections, namely 
showing that when each $\G_k^t$ is sufficiently large, $H$ is decreasing,
and the last two terms in \eqref{eq:tat-target-final}, in sum, are non-negative.

\newcommand{\Lktaukttautauo}{L_{\ktau,k_t}^{[\tau,\tau+1]}}

\subsubsection{$H$ is a Decreasing Function}\label{sect:CES-H-decreasing}

For any time $\tau$ at which there is an update, 
from \eqref{eq:tat-target-first}
and the definition of $A$, we have
\begin{align*}
H(\alpha(\tau)) - H(\tau)
&~\ge~ \frac 12 \hWktau(\gktautau,\pktautaum,\Gktautau) \cdot \Dtau ~+~ \frac {\Gktautau}8 \cdot \frac{(\Delta p_{\ktau}^\tau)^2}{\Dtau}
~-~\frac{1}{\Gktautau}  \left(\gktaumaxtau - \gktaumintau\right)^2 \cdot \Dtau\\
&~~~~~~~~~~~~
~+~ \frac 13 \sum_{\nu\in (\alpha_{\ktau}(\tau),\tau)} 
L_{k_\nu,\ktau}^{[\nu,\nu+1]} \cdot \frac{p_{\ktau}^{\tau-}}{p_{k_\nu}^{\nu-}}
\cdot \frac{\left(\Delta p_{k_\nu}^\nu\right)^2}{\Delta t_\nu}
~-~ \frac 23 \sum_{k\neq \ktau} L_{\ktau,k}^{[\tau,\tau+1]} \cdot \frac{p_k^\tau}{p_{\ktau}^{\tau-}}
\cdot \frac{\left(\Delta p_{\ktau}^\tau\right)^2}{\Dtau}.
\end{align*}

Next,
\begin{align}
\left( \gktaumaxtau - \gktaumintau \right)^2
&~\le~ \left( \sum_{\nu\in (\alpha_{\ktau}(\tau),\tau)} L_{k_\nu,\ktau}^{[\nu,\tau]} \left| \Delta p_{k_\nu}^\nu \right| \right)^2\nonumber\\
&~=~ \left( \sum_{\nu\in (\alpha_{\ktau}(\tau),\tau)} \left(L_{k_\nu,\ktau}^{[\nu,\tau]} \cdot \Delta t_\nu \cdot \frac{p_{k_\nu}^{\nu-}}{p_{\ktau}^{\nu-}}\right)
\cdot \left(\frac{\left| \Delta p_{k_\nu}^\nu \right|}{\Delta t_\nu} \cdot \frac{p_{\ktau}^{\nu-}}{p_{k_\nu}^{\nu-}}\right) \right)^2\nonumber\\
&~\le~ \left( \sum_{\nu\in (\alpha_{\ktau}(\tau),\tau)} L_{k_\nu,\ktau}^{[\nu,\tau]} \cdot \Delta t_\nu
\cdot \frac{p_{k_\nu}^{\nu-}}{p_{\ktau}^{\nu-}} \right)
\left(
\sum_{\nu\in (\alpha_{\ktau}(\tau),\tau)} L_{k_\nu,\ktau}^{[\nu,\tau]} \cdot \frac{p_{\ktau}^{\nu-}}{p_{k_\nu}^{\nu-}}
\cdot \frac{\left( \Delta p_{k_\nu}^\nu \right)^2}{\Delta t_\nu}
\right)
\comm{by Lemma \ref{lem:power-mean-ine}}\nonumber\\
&~\le~ \left( \sum_{\nu\in (\alpha_{\ktau}(\tau),\tau)} L_{k_\nu,\ktau}^{[\alpha_{\ktau}(\tau),\tau]} \cdot \Delta t_\nu
\cdot \frac{p_{k_\nu}^{\nu-}}{p_{\ktau}^{\tau-}} \right)
\left(
\sum_{\nu\in (\alpha_{\ktau}(\tau),\tau)} L_{k_\nu,\ktau}^{[\nu,\nu+1]} \cdot \frac{p_{\ktau}^{\tau-}}{p_{k_\nu}^{\nu-}}
\cdot \frac{\left( \Delta p_{k_\nu}^\nu \right)^2}{\Delta t_\nu}
\right).\label{eq:tat-gradient-error}
\end{align}

Combining all the above yields
\begin{align*}
&~~~~~H(\alpha(\tau)) - H(\tau)\\
&~\ge~ \frac 12 \hWktau(\gktautau,\pktautaum,\Gktautau) \cdot \Dtau
~+~ \left( \frac {\Gktautau}8 - \frac 23 \sum_{k\neq \ktau} L_{\ktau,k}^{[\tau,\tau+1]} \cdot \frac{p_k^\tau}{p_{\ktau}^{\tau-}} \right)
\frac{\left(\Delta p_{\ktau}^\tau\right)^2}{\Dtau}\\
&~~~~~~~~~~~~+~ \left( \frac 13 - \frac{1}{\Gktautau} \sum_{\nu\in (\alpha_{\ktau}(\tau),\tau)} L_{k_\nu,\ktau}^{[\alpha_{\ktau}(\tau),\tau]} \cdot \Delta t_\nu
\cdot \frac{p_{k_\nu}^{\nu-}}{p_{\ktau}^{\tau-}} \right)\left(
\sum_{\nu\in (\alpha_{\ktau}(\tau),\tau)} L_{k_\nu,\ktau}^{[\nu,\nu+1]} \cdot \frac{p_{\ktau}^{\tau-}}{p_{k_\nu}^{\nu-}}
\cdot \frac{\left( \Delta p_{k_\nu}^\nu \right)^2}{\Delta t_\nu}
\right).
\end{align*}

Thus, for $H$ to be decreasing, it suffices that
$$
\Gktautau ~\ge~ \frac{16}{3} \sum_{k\neq \ktau} L_{\ktau,k}^{[\tau,\tau+1]} \cdot \frac{p_k^\tau}{p_{\ktau}^{\tau-}}~~~~~~~~\text{and}~~~~~~~~
\Gktautau ~\ge~ 3 \sum_{\nu\in (\alpha_{\ktau}(\tau),\tau)} L_{k_\nu,\ktau}^{[\alpha_{\ktau}(\tau),\tau]} \cdot \Delta t_\nu
\cdot \frac{p_{k_\nu}^{\nu-}}{p_{\ktau}^{\tau-}}.
$$
But we will impose the stronger requirement that
$$
\Gktautau ~\ge~ 6 \sum_{\nu\in (\alpha_{\ktau}(\tau),\tau)} L_{k_\nu,\ktau}^{[\alpha_{\ktau}(\tau),\tau]} \cdot \Delta t_\nu
\cdot \frac{p_{k_\nu}^{\nu-}}{p_{\ktau}^{\tau-}}.
$$

\subsubsection{The Sum of the Last Two Terms in \eqref{eq:tat-target-final} is Non-negative}

It remains to show that the sum of the last two terms in \eqref{eq:tat-target-final} is non-negative, i.e.,
\begin{align*}
\frac 18 \sum_{\tau\in (t_a,t]} \Gktautau \frac{(\Delta p_{\ktau}^\tau)^2}{\Dtau} - \frac 1{12} \sum_{\nu\in (t_a,t_b]} \delnu\cdot A(\alpha(\nu)) ~+~ A(t_a) - A(t) ~~\ge~~
\frac{5}{3} \sum_{\tau\in (t_a,t]} \frac{1}{\Gktautau} \cdot ( \gktaumaxtau - \gktaumintau )^2\cdot \Dtau.
\end{align*}
We first simplify the LHS using the definition of $A$:
\begin{align*}
&~~~~~~\frac 18 \sum_{\tau\in (t_a,t]} \Gktautau \frac{(\Delta p_{\ktau}^\tau)^2}{\Dtau}
~-~ \frac 1{12} \sum_{\nu\in (t_a,t_b]} \delnu\cdot A(\alpha(\nu))
~+~ A(t_a) ~-~ A(t)\\
&~\ge~ \frac 18 \sum_{\tau\in (t_a,t]} \Gktautau \frac{(\Delta p_{\ktau}^\tau)^2}{\Dtau}
~-~ \frac 1{12} ~\sum_{\nu\in (t_a-1,t_b]}~ \frac 23~ \sum_{k\neq k_\nu}~ L_{k_\nu,k}^{[\nu,\nu+1]} \cdot \frac{p_k^{\nu}}{p_{k_\nu}^{\nu-}}
\cdot \frac{\left(\Delta p_{k_\nu}^\nu\right)^2}{\Delta t_\nu}\\
&~~~~~~~~~~
~+~ \frac{1}{3} ~\sum_{\nu\in (t_a-1,t_a]}~ \sum_{k\neq k_\nu}~ L_{k_\nu,k}^{[\nu,\nu+1]} \cdot \frac{p_k^{\nu}}{p_{k_\nu}^{\nu-}}
\cdot \frac{\left(\Delta p_{k_\nu}^\nu\right)^2}{\Delta t_\nu}
~-~ \frac{2}{3} ~\sum_{\nu\in (t-1,t]}~ \sum_{k\neq k_\nu}~ L_{k_\nu,k}^{[\nu,\nu+1]} \cdot \frac{p_k^{\nu}}{p_{k_\nu}^{\nu-}}
\cdot \frac{\left(\Delta p_{k_\nu}^\nu\right)^2}{\Delta t_\nu}\\
&~\ge~ \left( \frac 13 - \frac 1{18} \right)
\sum_{\nu\in (t_a-1,t_a]}~ \sum_{k\neq k_\nu}~ L_{k_\nu,k}^{[\nu,\nu+1]} \cdot \frac{p_k^{\nu}}{p_{k_\nu}^{\nu-}}
\cdot \frac{\left(\Delta p_{k_\nu}^\nu\right)^2}{\Delta t_\nu}
~~+~~ \sum_{\tau\in (t_a,t]}\left[\frac{\Gktautau}{8} - \left( \frac{1}{18} + \frac 23 \right) \sum_{k\neq k_\tau}~L_{\ktau,k}^{[\tau,\tau+1]} \cdot \frac{p_k^{\tau}}{p_{\ktau}^{\tau-}}\right]~\frac{(\Delta p_{\ktau}^\tau)^2}{\Dtau}.
\end{align*}

By imposing the requirement that
$\Gktautau ~\ge~ 8 ~ \sum_{k\neq k_\tau}~L_{\ktau,k}^{[\tau,\tau+1]} \cdot \frac{p_k^{\tau}}{p_{\ktau}^{\tau-}}$, we have
$$
\frac 18 \sum_{\tau\in (t_a,t]} \Gktautau \frac{(\Delta p_{\ktau}^\tau)^2}{\Dtau}
~-~ \frac 1{12} \sum_{\nu\in (t_a,t_b]} \delnu\cdot A(\alpha(\nu))
~+~ A(t_a) ~-~ A(t)
~\ge~ \frac 5{18} \sum_{\tau\in(t_a-1,t]} ~\sum_{k\neq \ktau}~ L_{k_\tau,k}^{[\tau,\tau+1]} \cdot \frac{p_k^{\tau}}{p_{\ktau}^{\tau-}}
~\frac{(\Delta p_{\ktau}^\tau)^2}{\Dtau}.
$$

On the other hand, by \eqref{eq:tat-gradient-error} and by the condition imposed on the $\G$s in the last subsection,
\begin{align*}
\frac{5}{3} \sum_{\tau\in (t_a,t]} \frac{1}{\Gktautau} \cdot ( \gktaumaxtau - \gktaumintau )^2\cdot \Dtau
&~\le~ \frac{5}{3} \sum_{\tau\in (t_a,t]}~ \frac{1}{6}
~\sum_{\nu\in ( \alpha_{\ktau}(\tau) , \tau )}~ L_{k_\nu,\ktau}^{[\nu,\nu+1]} \cdot \frac {\left( \Delta p_{k_\nu}^\nu \right)^2} {\Delta t_\nu}\\
&~\le~ \frac {5}{18} \sum_{\tau\in(t_a-1,t]} ~\sum_{k\neq \ktau}~ L_{k_\tau,k}^{[\tau,\tau+1]} \cdot \frac{p_k^{\tau}}{p_{\ktau}^{\tau-}}
~\frac{(\Delta p_{\ktau}^\tau)^2}{\Dtau}.
\end{align*}

\subsubsection{Upper Bounds on the Local Lipschitz Parameters, and Determining the $\G$'s}

In the last two subsections, we have imposed the requirements
\begin{equation}\label{eq:Gamma-requirement}
\Gktautau ~\ge~ 8 ~ \sum_{k\neq k_\tau}~L_{\ktau,k}^{[\tau,\tau+1]} \cdot \frac{p_k^{\tau}}{p_{\ktau}^{\tau-}} ~~~~~~\text{and}~~~~~~
\Gktautau ~\ge~ 6 \sum_{\nu\in (\alpha_{\ktau}(\tau),\tau)} L_{k_\nu,\ktau}^{[\alpha_{\ktau}(\tau),\tau]} \cdot \Delta t_\nu
\cdot \frac{p_{k_\nu}^{\nu-}}{p_{\ktau}^{\tau-}},
\end{equation}
where the $L$ parameters are the \emph{local} Lipschitz parameters of the function $\phi$.
Our remaining tasks are to derive lower bounds on the two summations above.

Suppose in a Fisher market with buyers having complementary-CES utility functions,
each buyer $i$ has a budget of $e_i$, and her CES utility function has parameter $\rho_i$.
For each $i$, let $\theta_i := \rho_i / (\rho_i - 1)$.
As we have discussed in Section \ref{sect:CES}, at any given prices $p\in (\rrplus)^n$,
buyer $i$ computes the demand-maximizing bundle of goods costing at most $e_i$;
we let $x_{i\ell}(p)$ denote the buyer $i$'s demand for good $\ell$ at prices $p$.

In a Fisher market with buyers having complementary-CES utility functions,
Properties 1 and 2 below are well-known. Property 3 was proved in~\cite{CCD2013}.
\begin{enumerate}
	\item For any $k\neq j$,
	$$\left|\frac{\partial^2 \phi}{\partial p_j ~ \partial p_k}(p) \right| ~=~ \sum_i \frac{\theta_i~x_{ij}(p)~x_{ik}(p)}{e_i} ~\leq~ \sum_i \frac{x_{ij}(p)~x_{ik}(p)}{e_i}.$$
	\item Given positive prices $p$, for any $0 < r_1 < r_2$,
	let $p'$ be prices such that for all $\ell$, $r_1 p_\ell \leq p'_\ell \leq r_2 p_\ell$.
	Then for all $\ell$, $\frac{1}{r_2} x_\ell(p) \leq x_\ell(p') \leq \frac{1}{r_1} x_\ell(p)$.
	\item If for each $\ell$, $\left|\frac{\Dp_\ell}{p_\ell}\right| \leq \frac 16$, then
	$$\phi(p+\Dp) - \phi(p) - \sum_\ell \nabla_\ell\phi(p) \cdot \Dp_\ell ~\le~ \sum_\ell \frac{1.5 x_\ell}{p_\ell} (\Dp_\ell)^2.$$
\end{enumerate}

\begin{lemma}\label{lem:ces-bounds}
	If the parameter $\lambda$ in update rule \eqref{eq:async-tat-CES-rule} is at most $1/20$, then
$$
\sum_{k\neq k_\tau} L_{k_\tau,k}^{[\tau,\tau+1]} \cdot \frac{p_k^{\tau}}{p_{\ktau}^{\tau-}} ~\le~ \frac{5}{4}\cdot \frac{x_{\ktau}(p^{\tau-})}{p_{\ktau}^{\tau-}}
~~~~~~\text{and}~~~~~~
\sum_{\nu\in (\alpha_{\ktau}(\tau),\tau)} L_{k_\nu,\ktau}^{[\alpha_{\ktau}(\tau),\tau]} \cdot \Delta t_\nu
\cdot \frac{p_{k_\nu}^{\nu-}}{p_{\ktau}^{\tau-}}
~\le~ \frac{11}{4}\cdot \frac{x_{\ktau}(p^{\tau-})}{p_{k_\tau}^{\tau-}}.
$$
\end{lemma}

\begin{pf}
Since $\lambda \leq 1/20$, it is easy to observe that for any $\nu\in [\tau,\tau+1]$
and for any $k$ (including coordinate $\ktau$),
\begin{equation}\label{eq:price-range}
e^{-2/20}\cdot p_k^{\tau-} ~~\le~~ p_k^{\nu} ~~\le~~ e^{2/19}\cdot p_k^{\tau-}.
\end{equation}
Accordingly, we let
$$\tP := \left\{(\tp_1,\tp_2,\cdots,\tp_n)\,\left|\,\forall k\in [n],~
e^{-2/20} \cdot p_k^{\tau-} ~\le~ \tp_k ~\le~ e^{2/19} \cdot p_k^{\tau-}\right.\right\}.$$

\begin{align}
\sum_{k\neq k_\tau} L_{k_\tau,k}^{[\tau,\tau+1]} \cdot \frac{p_k^{\tau}}{p_{\ktau}^{\tau-}}
&~\le~ \frac{1}{p_{\ktau}^{\tau-}}~\sum_{k\neq \ktau} \left(\max_{\tp\in\tP}~\left|\frac{\partial^2 \phi}{\partial p_{\ktau} ~ \partial p_k} (\tp)\right|\right) \cdot p_k^{\tau}\nonumber\\
&~\le~ \frac{1}{p_{\ktau}^{\tau-}}~\sum_{k\neq \ktau} \sum_i \frac{(e^{2/20} x_{i\ktau}(p^{\tau-}))\cdot(e^{2/20} x_{ik}(p^{\tau-}))}{e_i}
\cdot p_k^{\tau-}  \comm{By Properties 1 and 2}\nonumber\\
&~\le~ \frac{e^{1/5}}{p_{\ktau}^{\tau-}}~\sum_i~x_{i\ktau}(p^{\tau-})~\sum_{k\neq k_\tau} \frac{x_{ik}(p^{\tau-})\cdot p_k^{\tau-}}{e_i}\nonumber\\
&~\le~ \frac{e^{1/5}}{p_{\ktau}^{\tau-}}~\sum_i~x_{i\ktau}(p^{\tau-})\comm{The second summation is at most $1$, due to the budget constraint}\nonumber\\
&~=~ e^{1/5}\cdot \frac{x_{\ktau}(p^{\tau-})}{p_{\ktau}^{\tau-}} ~~\le~~
\frac 54 \cdot \frac{x_{\ktau}(p^{\tau-})}{p_{\ktau}^{\tau-}}.\label{eq:first-part}
\end{align}

For the time range $\nu\in [\alpha_{\ktau}(\tau),\tau]$, the inequality \eqref{eq:price-range} holds also.

\begin{align*}
\sum_{\nu\in (\alpha_{\ktau}(\tau),\tau)}~L_{k_\nu,\ktau}^{[\alpha_{\ktau}(\tau),\tau]} \cdot \Delta t_\nu
\cdot \frac{p_{k_\nu}^{\nu-}}{p_{\ktau}^{\tau-}}
&~\le~ \frac{1}{p_{\ktau}^{\tau-}}~\sum_{\nu\in (\alpha_{\ktau}(\tau),\tau)}~L_{k_\nu,\ktau}^{[\alpha_{\ktau}(\tau),\tau]} \cdot \Delta t_\nu
\cdot e^{2/19}\cdot p_{k_\nu}^{\tau-}\\
&~\le~ \frac{e^{2/19}}{p_{\ktau}^{\tau-}}~\sum_{k\neq \ktau}~L_{k,\ktau}^{[\alpha_{\ktau}(\tau),\tau]} \cdot p_k^{\tau-} \cdot
\sum_{\stackrel{\nu\in(\alpha_{\ktau}(\tau),\tau)}{k_\nu=k}} \Delta t_\nu\\
&~\le~ \frac{2e^{2/19}}{p_{\ktau}^{\tau-}}~\sum_{k\neq \ktau}~L_{k,\ktau}^{[\alpha_{\ktau}(\tau),\tau]} \cdot p_k^{\tau-}.~~~~~~\text{(Observe that the $\sum_\nu \Delta t_\nu$ term}\\
&~~~~~~~~~~~~~~~~~~~~~~~~~~~~~~~~~~~~~~~~~~~~~~~~~~~~
\text{above is at most $2$)}
\end{align*}
The summation $\sum_{k\neq \ktau}~L_{k,\ktau}^{[\alpha_{\ktau}(\tau),\tau]} \cdot p_k^{\tau-}$
above can be bounded as in \eqref{eq:first-part}, yielding an upper bound of $e^{1/5}\cdot x_{\ktau}(p^{\tau-})$.
Noting that $2\cdot e^{2/19}\cdot e^{1/5} ~\le~ \frac{11}{4}$, we are done.
\end{pf}

To conclude, by \eqref{eq:Gamma-requirement} and Lemma \ref{lem:ces-bounds}, we need that
$$
\Gktautau ~\ge~ \frac{33}{2} \cdot \frac{x_{\ktau}(p^{\tau-})}{p_{\ktau}^{\tau-}}.
$$
Note that in update rule \eqref{eq:async-tat-CES-rule},
it is equivalent to that $\Gktautau ~=~ \frac{1}{\lambda\cdot p_{\ktau}^{\tau-}} \cdot \max\{\tilde{z}_{\ktau},1\}$.
Thus, we need that
$$
\lambda ~\le~ \frac{2}{33}\cdot \frac{\max\{\tilde{z}_{\ktau},1\}}{x_{\ktau}(p^{\tau-})} ~\le~
\frac{2}{33}\cdot \frac{\max\{\tilde{z}_{\ktau},1\}}{e^{2/19}\cdot \tilde{x}_{\ktau}} ~\le~
\frac{2}{33\cdot e^{2/19}} \cdot \frac 12,
$$
or slightly stronger, $\lambda ~\le~ 1/37$.

\newpage

\section{Tatonnement Analysis for Leontief Fisher Markets}
\label{app:leontief}

It is well-known that Leontief utility functions can be considered as the ``limit'' of CES utility functions as $\rho\ra -\infty$.
We note that the three properties listed in Appendix \ref{sect:CES} also hold for Leontief utility functions.
However, we cannot apply the analysis in Appendix \ref{sect:CES} to the Leontief Fisher markets for two reasons.
Firstly, while $\phi$ remains convex, it is no longer strongly convex, so $\mu_\phi = 0$.
Secondly, recall that $\Gktautau ~\ge~ \frac{33}{2} \cdot \frac{x_{\ktau}(p^{\tau-})}{p_{\ktau}^{\tau-}}$.
For Leontief Fisher markets, it is possible that some good $j$ has zero equilibrium price.
Under this scenario, for convergence to the equilibrium, $\G_j^t$ has to grow towards infinity.

Here, we provide additional arguments which build on top of the result that $H(t)$ decreases with $t$, to show that tatonnement
with update rule \eqref{eq:async-tat-CES-rule} still converges toward the market equilibrium.
However, this result does not provide a bound on the rate of convergence.

Tatonnement in Leontief Fisher markets was first analysed by Cheung, Cole and Devanur~\cite{CCD2013}.
They gave a bound on the convergence rate, but with a less natural update rule --- in their update rule,
$\Gjt$ increases with the number of buyers in the market, and is also a function of the demands for all the goods,
both of which seem unnatural, while the $\Gjt$ used here is independent of
the number of buyers and depends only on the demand for good $j$.

\subsection{Analysis}

\begin{lemma}\label{lem:one-big-move-implies-big-progress}
Let $\alpha_j(t),t$ be the times at which two consecutive updates to $p_j$ occur. Let $\Delta t = t - \alpha_j(t)$
Then $H(\alpha_j(t)) - H(t) ~\ge~ \frac{\G_j^t}{8} \cdot \frac{\left(\Delta p_j^t\right)^2}{\Delta t}$.
\end{lemma}

\begin{pf}
Same as in Section \ref{sect:CES-H-decreasing},
except that instead of using \eqref{lem:F-progress-for-tat}, we use the following inequality instead:
$$
F(p^t) - F(p^{\alpha(t)})  ~\ge~ \frac{\G_j^t}{4}\frac{\left(\Delta p_j^t\right)^2}{\Delta t}
~-~ \frac {1}{\G_j^t}(g_j - \tg_j)^2 \cdot \Dt.
$$
\end{pf}

In~\cite{CCD2013}, they showed that there exists a finite positive number $U$ which is an upper bound on all the prices
throughout the tatonnement process.

\begin{lemma}\label{lem:big-moves-in-one-time-unit-implies-big-overall-progress}
	Suppose that there are consecutive updates to $p_j$ at times $\tau_0 < \tau_1 < \cdots < \tau_m$, where $\tau_m - \tau_0 \leq 2$.
	If $\left|p_j^{\tau_0} - p_j^{\tau_m}\right| \geq \epsilon$, where $\epsilon \leq 1$, then
	$H(\tau_0) - H(\tau_m) ~\ge~ \epsilon^2\cdot \min\left\{\frac{1}{16},\frac{1}{64\lambda U}\right\}$.
\end{lemma}
\begin{pf}
	\newcommand{\tz}{\tilde{z}}
	For $q=1,2,\cdots,m$, let $\Delta\pjq$ be the change made to $p_j$ by the update at time $\tau_q$, 
	and let $\tzjq$ be the $\tz$-value used for the update,
	i.e., $\G_j^{\tau_q} = \frac{\max\{1,\tzjq\}}{\lambda p_j^{\tau_q^-}}$ and
	$\Delta\pjq = \lambda p_j^{\tau_q^-}\cdot \min\{1,\tzjq\}\cdot \Dt_q$.
	
	If $\tzjq < 1$, then
	$$\frac{\G_j^{\tau_q}(\D\pjq)^2}{\Dt_q} = \frac{1}{\lambda p_j^{\tau_q^-}}\frac{(\D\pjq)^2}{\Dt_q}\geq \frac{1}{\lambda U}\frac{(\D\pjq)^2}{\Dt_q}.$$
	
	If $\tzjq \geq 1$, then
	$$\frac{\G_j^{\tau_q}(\D\pjq)^2}{\Dt_q} ~=~ \frac{\tzjq}{\lambda p_j^{\tau_q^-}}\cdot \lambda^2 \left(p_j^{\tau_q^-}\right)^2\cdot \Dt_q
	~=~ \lambda p_j^{\tau_q^-} \tzjq \cdot \Dt_q \geq |\D\pjq|.$$
	By Lemma \ref{lem:one-big-move-implies-big-progress},
	\begin{align*}
	H(\tau_0) - H(\tau_m) ~=~ \sum_{q=1}^m \left(H(\tau_{q-1}) - H(\tau_q)\right) &~\ge~ \frac 18~\sum_{q=1}^m \frac{\G_j^{\tau_q} (\D\pjq)^2}{\Dt_q}\\
	&~\ge~ \frac{1}{8\lambda U}\sum_{q:\tzjq < 1} \frac{(\D\pjq)^2}{\Dt_q} ~+~ \frac 18 \sum_{q:\tzjq\geq 1} |\D\pjq|.
	\end{align*}
	
	By the assumption $|p_j^{\tau_0} - p_j^{\tau_m}| \geq \epsilon$, $\sum_{q=1}^m |\D\pjq| \geq \epsilon$.
	Let $\sigma := \epsilon^{-1} \sum_{q:\tzjq\geq 1} |\D\pjq|$. Then $\sum_{q:\tzjq<1} |\D\pjq| \geq \max\{0,(1-\sigma)\epsilon\}$.
	By the Cauchy-Schwarz inequality,
	\begin{align*}
	\left[\max\{0,(1-\sigma)\epsilon\}\right]^2 \leq \left(\sum_{q:\tzjq<1} |\D\pjq|\right)^2 &= \left(\sum_{q:\tzjq<1} \left|\frac{\D\pjq}{\sqrt{\Dt_q}}\right|\cdot \sqrt{\Dt_q}\right)^2\\
	&\leq \left(\sum_{q:\tzjq<1} \frac{(\D\pjq)^2}{\Dt_q}\right)\left(\sum_{q:\tzjq<1} \Dt_q\right)\\
	&\leq 2\sum_{q:\tzjq<1} \frac{(\D\pjq)^2}{\Dt_q},
	\end{align*}
	as $\tau_m - \tau_0 \leq 2$. Then $\sum_{q:\tzjq<1} \frac{(\D\pjq)^2}{\Dt_q} \geq \frac{1}{2}\left[\max\{0,(1-\sigma)\epsilon\}\right]^2$ and hence
	$$H(\tau_0) - H(\tau_m) \geq \frac{1}{16\lambda U} \left[\max\{0,(1-\sigma)\epsilon\}\right]^2 ~+~ \frac{\sigma \epsilon}{8}.$$
	By considering the following three cases: $\sigma \geq 1$, $1>\sigma\geq 1/2$ or $\sigma < 1/2$,
	it is not difficult to show that the minimum value of R.H.S.~of
	the above inequality is at least $\epsilon^2\cdot \min\left\{\frac{1}{16},\frac{1}{64\lambda U}\right\}$.
\end{pf}

\begin{cor}\label{cor:small-move-in-one-time-unit}
	For any $\epsilon > 0$, there exists a finite time $T_\epsilon$ such that for any good $j$,
	any $t\geq T_\epsilon$, and any $0\leq \Dt \leq 1$, $|p_j^t - p_j^{t+\Dt}| \leq \epsilon$.
\end{cor}
\begin{pf}
Suppose not, then by Lemma \ref{lem:big-moves-in-one-time-unit-implies-big-overall-progress},
$H$ drops by at least $\epsilon^2\cdot \min\left\{\frac{1}{16},\frac{1}{64\lambda U}\right\}$
infinitely often. But $H(0)$ is finite and $H$ remains positive throughout, a contradiction.
\end{pf}

\bigskip

\begin{pfof}{Theorem \ref{thm:CES-Fisher-cnvge} for the Leontief case}
The proof comprises four steps. We need the following definitions:
for any two price vectors $p^A$ and $p^B$, let $d(p^A,p^B)$ denote the $L_1$ norm distance between the two price vectors,
i.e., $d(p^A,p^B) = \sum_j |p^A_j - p^B_j|$.
For any two sets of price vectors $P^A$ and $P^B$, let $d(P^A,P^B) := \inf_{p^A\in P^A,\, p^B\in P^B} d(p^A,p^B)$.

\medskip

\noindent\textbf{Step 1.} Let $\Omega$ be the set of limit points of a tatonnement process. We show that $\Omega$ is non-empty and connected.

\smallskip

Since all prices remain bounded by $U$ throughout the tatonnement process, $\Omega$ is non-empty.

Suppose $\Omega$ is not connected. Let $\Omega_a$ denote a connected component of $\Omega$
that is well separated from $\Omega_b = \Omega\setminus \Omega_a$,
i.e., $d(\Omega_a,\Omega_b) = \ep ' > 0$ (if there is no such $\Omega_a$ then $\Omega$ is connected).
By the definition of limit points, there exists a finite time such that
thereafter the prices in the tatonnement process are always within an $\epsilon'/4$-neighborhood of either $\Omega_a$ or $\Omega_b$.
This forces an infinite number of updates, each separated by at least one time unit,
such that each update makes a change to a price by at least at least $\epsilon'/(2n)$.
This contradicts Corollary \ref{cor:small-move-in-one-time-unit}.

\medskip
	
\noindent\textbf{Step 2.} Recall that a market equilibrium is a price vector $\ps$ at which for each $j$,
$\ps_j>0$ implies $z_j(\ps)=0$ and $\ps_j=0$ implies $z_j(\ps)\leq 0$.
We define a pseudo-equilibrium: a price vector $\tp$ is a pseudo-equilibrium if for each $j$,
$\tp_j>0$ implies $z_j(\tp)=0$.
Note that every market equilibrium is a pseudo-equilibrium. We show that all limit points in $\Omega$ are pseudo-equilibria.

\smallskip

\newcommand{\ringp}{\grave{p}}
Suppose not. Let $p'\in \Omega$ be a price vector which is not a pseudo-equilibrium,
i.e., there exists $j$ such that $p'_j>0$ but $z_j(p')\neq 0$.
Let $\ep'' ~\le~ p'_j |z_j(p')|/32$ be a positive number such that for any price vector $\ringp$ in the $\ep''$-neighborhood of $p'$,
we must have $\ringp_j\geq p'_j/2$ and $z_j(\ringp)$ lie between $z_j(p')/2$ and $z_j(p')$.

By the definition of limit points, the tatonnement process enters the $(\ep''/2)$-neighborhood of $p'$ infinitely often.
By Corollary \ref{cor:small-move-in-one-time-unit}, there exists a finite time such that subsequently,
every time the tatonnement process enters the $\ep''/2$-neighborhood of $p'$,
it stays in the $\ep''$-neighborhood of $p'$ for at least three time units.
Within the first two time units, $p_j$ is updated for at least once,
and by update rule \eqref{eq:async-tat-CES-rule},
such updates will make a total change to $p_j$ of at least $\lambda (p'_j/2) (\left|z_j(p')\right|/2) ~\ge~ 8\ep''$,
which forces quitting the $\ep''$-neighborhood of $p'$ strictly before the three time unit interval, a contradiction.
	
\medskip

\noindent\textbf{Step 3.} We show that the excess demands at all limit points in $\Omega$ are identical.

\smallskip

For every subset of goods $S$, let $\Omega_S = \{p'\in \Omega\,|\,p'_k > 0 \Leftrightarrow k\in S\}$.
For each buyer, there are two cases:
\begin{itemize}
	\item \textbf{the buyer wants at least one good in $S$, say good $\ell$:}\\
	Observe that by the definition of pseudo-equilibrium and Step 2, every price vector in $\Omega_S$,
	excluding the zero prices in the price vector,
	is a market equilibrium for the sub-Leontief-market comprising the goods in $S$.
	Codenotti and Varadarajan \cite{CV2004} pointed out that the demands for the goods in $S$ of each buyer are identical at every market equilibrium of the sub-Leontief market,
	and hence also in the original Leontief market.
	So the buyer demands the same positive but finite amount of good $\ell$ at every price vector in $\Omega_S$
	in the original market.
	Also note that the buyer always demands the goods in the original market in a fixed proportion.
	This forces the demands for the goods not in $S$ of the buyer to also be identical at every price vector in $\Omega_S$.
	\item \textbf{the buyer wants no good in $S$:}\\
	Then the buyer demands an infinite amount of each good that she wants, and demands zero amount of each good that she does not want.
\end{itemize}
In either case, the buyer's demands for each good at every price vector in $\Omega_S$ are identical,
and hence also the total demand for each good.

Then consider a graph $G$ with each vertex corresponding to a subset of goods $S$ such that $\Omega_S$ is non-empty,
and two vertices $S_1,S_2$ being adjacent if and only if $d\left(\Omega_{S_1},\Omega_{S_2}\right) = 0$.
Since excess demands are a continuous function\footnote{The range of the excess demand functions is the extended real line $\rr \cup \{+\infty\}$; continuity of the excess demand function is w.r.t.~the usual topology on the extended real line. To be specific, if $z_k(p)=+\infty$ for some $p$ and $k$, then for any $M\in \rr$, there exists an $\ep_M>0$ such that $z_k(p)\geq M$ in the $\ep_M$-neighborhood of $p$.} of prices,
if $S_1$ and $S_2$ are adjacent, then the excess demands for all goods at every price vector in $S_1\cup S_2$ are identical.
By Step 1, the graph $G$ is connected, thus the excess demands at all limit points in $\Omega$ are identical.

\bigskip

\noindent\textbf{Step 4.} We show that every limit point in $\Omega$ is indeed a market equilibrium.

\smallskip

Suppose not, i.e., there exists a limit point $p'$ in $\Omega$ which is a pseudo-equilibrium but not a market equilibrium,
i.e., there exists $k$ such that $p'_k=0$ but $z_k(p')>0$.
By Step 3, $z_k$ is positive at every limit point in $\Omega$, and hence every $p_k$ at every limit point must be zero.
By the definition of limit points, for any $\epsilon>0$, beyond a finite time, the tatonnement process must stay within the $\epsilon$-neighborhood of $\Omega$ thereafter.
By choosing a sufficiently small $\epsilon$, $z_k$ is bounded away from zero in the $\epsilon$-neighborhood of $\Omega$, and hence $p_k$ increases indefinitely and eventually $p_k$ becomes so large that the tatonnement process must leave the $\epsilon$-neighborhood of $\Omega$, a contradiction.
\end{pfof}

\newpage

\section{Toward a CCD Lower Bound?}
\label{app:eigenvalues}

\newcommand{\ao}{a_1}
\newcommand{\ai}{a_i}
\newcommand{\at}{a_2}
\newcommand{\an}{a_n}
\newcommand{\bo}{b_1}
\newcommand{\bi}{b_i}
\newcommand{\bt}{b_2}
\newcommand{\btnmi }{b_{2n-i}}
\newcommand{\ben}{b_n}
\newcommand{\bnp}{b_{n+1}}
\newcommand{\btn}{b_{2n}}
\newcommand{\bn}{b_n}
\newcommand{\lam}{\lambda}
\newcommand{\An}{A_n}
\newcommand{\Bn}{B_n}
\newcommand{\Dn}{D_n}
\newcommand{\Mn}{M_n}
\newcommand{\lmax}{\lambda_{\max}}

We present a family of strongly convex functions of $4n$ coordinates for which the 
relationship between the squares of the gradient differences and the squares of the changes in the coordinates, $(\Dx_{k_t})^2$, 
given in \eqref{eqn:g-to-x-log-bound} are tight, for suitable choices of the $\Dx_{k_t}$.
This suggests that the $\log n$ factor in our bound is not artificial.

We will be considering quadratic functions whose Hessians are  $4n\times 4n$ matrices of the form $\Mn = \Dn + \An$, 
 where $\Dn$ is a diagonal matrix with entry $2\lmax$ and $\An$ is a symmetric matrix  with eigenvalues all bounded in magnitude by $\lmax$.
It then follows that $\Mn$ has eigenvalues in the range $[\lmax,3\lmax]$.
We will show that $\lmax \le 4.5$, though our computation on matrices of size up to $1200$ strongly suggest it is about $3.68$.

$\An$ is the following matrix.
Its first row has the form 
\[
(0,\ao,0,\at,\ldots,0,\an,0,-\an,0,\ldots,0,-\at,0-\ao).
\]
Succesive rows are obtained by rotating one position to the right and flipping the signs. We choose $\ai =\frac 1i$. For instance,
\begin{equation*}
A_2 = 
\left(
\begin{array}{rrrrrrrr}
0 & 1 & 0 & \frac 12 & 0 & -\frac 12 & 0 & -1 \\
1 & 0 & -1 & 0 & -\frac 12 & 0 & \frac 12 & 0 \\
0 & -1 & 0 & 1 & 0 & \frac 12 & 0 & -\frac 12 \\
\frac 12 & 0 & 1 & 0 & -1 & 0 & -\frac 12 & 0 \\
0 & -\frac 12 & 0 & -1 & 0 & 1 & 0 & \frac 12 \\
-\frac 12 & 0 & \frac 12 & 0 & 1 & 0 & -1 & 0 \\
0 & \frac 12 & 0 & -\frac 12 & 0 & -1 & 0 & 1 \\
-1 & 0 & -\frac 12 & 0 & \frac 12 & 0 & 1 & 0 
\end{array}
\right)
\end{equation*}

Next, we prove the claimed bound on $\An$'s eigenvalues.
Define $\Bn = \An^2$.
Let $\Bn$'s first row be written as
\[
(\bo,0,\bt,0,\ldots,\ben,0,\bnp,0,\ldots,\btn,0).
\]
Note that if $\lam$ is an eigenvalue of $\An$ then $\lam^2$ is an eignevalue of $\Bn$. 
We will show that $\Bn$'s eigenvalues are all bounded by $2\pi^2$ in magnitude
and hence $\An$'s eigenvalues are at most $\sqrt 2 \pi$ in magnitude  (note that as $\An$ is symmetric, all its eigenvalues are real).

Note that 
\[
\bo = 2 \sum_{i=1}^n \frac{1}{i^2} \le 2\cdot \frac{\pi^2}{6} = \frac {\pi^2}{3}.
\]
For $2 \le i \le n$,
\begin{align*}
\bi & =~ - \sum_{h=1}^{i-1} \frac ih \cdot \frac{1}{i-h} 
+ \sum_{h=i}^n \frac 1h \cdot \frac{1}{h-i+1}
- \sum_{h=1}^{i-1} \frac{1}{n-h+1} \cdot \frac{1}{n-i+1+h}
+ \sum_{h=i}^n \frac{1}{n-h+1} \cdot {1}{n+i-h} \\
&= ~- \frac 1i \sum_{h=1}^{i-1} \left( \frac 1h + \frac {1}{i-h} \right)
+ \frac {1}{i-1} \sum_{h=i}^n \left( \frac {1}{h-i+1} - \frac 1h \right)
- \frac {1}{2n-i+2} \sum_{h=1}^{i-1} \left( \frac{1}{n-h+1} + \frac {1}{n-i+1+h} \right)\\
&~~~~~~~~~~~~~~~~+ \frac {1}{i-1} \sum_{h=i}^n \left( \frac {1}{n-h+1} - \frac {1}{n+i-h} \right) \\
& = ~- \frac 2i \sum_{h=1}^{i-1} \frac 1h 
+\frac {1}{i-1} \left[ \sum_{j=1}^{i-1} \frac 1j - \sum_{j=n-i+2}^n \frac 1j \right]
- \frac{2}{2n-i+2} \sum_{h=1}^{i-1} \frac {1}{n-h+1}
+ \frac {1}{i-1} \left[ \sum_{j=1}^{i-1} \frac 1j - \sum_{j=n-i+2}^n \frac 1j \right] \\
& =~ \frac {2}{i(i-1)} \sum_{h=1}^{i-1} \frac 1h 
- \frac{2}{i-1} \sum_{j=n-i+2}^n \frac 1j
- \frac {2}{2n-i+2} \sum_{h=1}^{i-1} \frac {1}{n-h+1}.
\end{align*}
By the mirror symmetry of $\Bn$'s first row, $2 \le i \le n$, 
$\btnmi = \bi$.
Finally, 
\[
\btn = -2 \sum_{i=1}^n \frac 1i \cdot \frac {1}{n+1-i} < 0.
\]

Observe that $\sum_{i=1}^{2n} \bi = \left[ \sum_{i=1}^n ( \ai + (-\ai)\right]^2 = 0$.
Thus, $\sum_{i=1}^{2n} |\bi|$ is bounded by twice the sum of the
positive parts in the expressions for the $\bi$, namely:
\begin{align*}
\sum_{i=1}^{2n} |\bi|
& \le~ 2\left[ \bo + 2\sum_{i=2}^n \frac {2}{i(i-1)} \sum_{h=1}^{i-1} \frac 1h \right] \\
& \le \frac{2\pi^2}{3} 
+ 8 \sum_{i=2}^n \left[ \frac {1}{i-1} \sum_{h=1}^{i-1} \frac 1h 
- \frac 1i \sum_{h=1}^i \frac 1h + \frac {1}{i^2} \right] \\
& \le~ \frac{2\pi^2}{3} 
+8 \left[ 1 - \frac 1n \sum_{h=1}^n \frac 1n + \sum_{i=1}^n \frac{1}{i^2} \right] \\
& \le ~ \frac{2\pi^2}{3}  + \frac{8\pi^2}{6} ~=~ 2\pi^2.
\end{align*}
Note that $\sqrt 2 \pi \approx 4.44$.

\medskip

We finish by demonstrating the need for a logarithmic term.
Suppose $\Delta \pkt = 1$ for $1\le t \le 8n$.
In the expression $Q$ in \eqref{eqn:target-cyclic-final} (but now we have $4n$ coordinates instead of $n$),
for each $i\ge 4n$,
each sum $\sum_{j=i-4n+1}^i \left(g_{k_i}^j - g_{k_i}^i\right)^2$ 
includes the terms
$$\sum_{j=1}^{2n} \left(\sum_{h=1}^j \frac 1h\right)^2 ~\ge~ \sum_{j=1}^{2n} (\ln j)^2 ~\ge~ n \ln^2 n,$$
thus the sum on the RHS of $Q$ is at least $\frac{2}{3\G n}\cdot 4n\cdot n \ln^2 n = \frac{8}{3\G} n \ln^2 n$,
while the sum on the LHS of $Q$ is $8 \G n$,
demonstrating that for this analysis to work we need $\G^2 = \Omega(\ln^2 n)$.

Of course, this does not exclude the possibility that a different analysis could produce a tighter bound.

\newpage
\bibliographystyle{plain}
\bibliography{acd_references}

\begin{thebibliography}{10}

\bibitem{AN2010}
Carlos Al{\'{o}}s{-}Ferrer and Nick Netzer.
\newblock The logit-response dynamics.
\newblock {\em Games and Economic Behavior}, 68(2):413--427, 2010.

\bibitem{ABH1959}
Kenneth~J. Arrow, H.~D. Block, and Leonid Hurwicz.
\newblock On the stability of competitive equilibrium, ii.
\newblock {\em Econometrica}, 27(1):82--109, 1959.

\bibitem{ADG2014}
Haim Avron, Alex Druinsky, and A.~Gupta.
\newblock Revisiting asynchronous linear solvers: Provable convergence rate
  through randomization.
\newblock In {\em {IEEE} IPDPS}, pages 198--207, 2014.

\bibitem{BT2013}
Amir Beck and Luba Tetruashvili.
\newblock On the convergence of block coordinate descent type methods.
\newblock {\em {SIAM} Journal on Optimization}, 23(4):2037--2060, 2013.

\bibitem{BT1989}
Dimitri~P. Bertsekas and John~N. Tsitsiklis.
\newblock {\em Parallel and Distributed Computation: Numerical Methods}.
\newblock Prentice-Hall, 1989.

\bibitem{BT2000}
Dimitri~P. Bertsekas and John~N. Tsitsiklis.
\newblock Gradient convergence in gradient methods with errors.
\newblock {\em SIAM J. Optimization}, 10(3):627--642, 2000.

\bibitem{Blume1993}
Lawrence~E. Blume.
\newblock The statistical mechanics of strategic interaction.
\newblock {\em Games and Economic Behavior}, 5(3):387--424, 1993.

\bibitem{Borkar1998}
Vivek~S. Borkar.
\newblock Asynchronous stochastic approximations.
\newblock {\em SIAM J. Control and Optimization}, 36(3):662--663, 1998.

\bibitem{Chazelle2009}
Bernard Chazelle.
\newblock Natural algorithms.
\newblock In {\em SODA}, pages 422--431, 2009.

\bibitem{Chazelle2012}
Bernard Chazelle.
\newblock The dynamics of influence systems.
\newblock In {\em FOCS}, pages 311--320, 2012.

\bibitem{CT1993}
Gong Chen and Marc Teboulle.
\newblock Convergence analysis of a proximal-like minimization algorithm using
  {B}regman's function.
\newblock {\em SIAM J. Optimization}, 3(3):538--543, 1993.

\bibitem{CC2014-arxiv}
Yun~Kuen Cheung and Richard Cole.
\newblock Amortized analysis on asynchronous gradient descent.
\newblock {\em CoRR}, abs/1412.0159, 2014.

\bibitem{CCD2013}
Yun~Kuen Cheung, Richard Cole, and Nikhil~R. Devanur.
\newblock Tatonnement beyond gross substitutes? gradient descent to the rescue.
\newblock In {\em STOC}, pages 191--200, 2013.

\bibitem{CCR2012}
Yun~Kuen Cheung, Richard Cole, and Ashish Rastogi.
\newblock Tatonnement in ongoing markets of complementary goods.
\newblock In {\em EC}, pages 337--354, 2012.

\bibitem{CMV2005}
Bruno Codenotti, Benton McCune, and Kasturi~R. Varadarajan.
\newblock Market equilibrium via the excess demand function.
\newblock In {\em STOC}, pages 74--83, 2005.

\bibitem{CV2004}
Bruno Codenotti and Kasturi~R. Varadarajan.
\newblock Efficient computation of equilibrium prices for markets with leontief
  utilities.
\newblock In {\em ICALP}, pages 371--382, 2004.

\bibitem{CF2008}
Richard Cole and Lisa Fleischer.
\newblock Fast-converging tatonnement algorithms for one-time and ongoing
  market problems.
\newblock In {\em STOC}, pages 315--324, 2008.

\bibitem{CFR2010}
Richard Cole, Lisa Fleischer, and Ashish Rastogi.
\newblock Discrete price updates yield fast convergence in ongoing markets with
  finite warehouses.
\newblock {\em CoRR}, 2010.

\bibitem{CV1995}
Corinna Cortes and Vladimir Vapnik.
\newblock Support-vector networks.
\newblock {\em Mach. Learn.}, 20(3):273--297, September 1995.

\bibitem{Dohtani1993}
Akitaka Dohtani.
\newblock Global stability of the competitive economy involving complementary
  relations among commodities.
\newblock {\em Journal of Mathematical Economics}, 22(1):73 -- 83, 1993.

\bibitem{DAW2012}
John~C. Duchi, Alekh Agarwal, and Martin~J. Wainwright.
\newblock Dual averaging for distributed optimization: Convergence analysis and
  network scaling.
\newblock {\em {IEEE} Trans. Automat. Contr.}, 57(3):592--606, 2012.

\bibitem{FS2000}
Andreas Frommer and Daniel~B. Szyld.
\newblock On asynchronous iterations.
\newblock {\em Journal of Computational and Applied Mathematics},
  123(1-2):201--216, 2000.
\newblock Numerical Analysis 2000. Vol. III: Linear Algebra.

\bibitem{HuaY2015}
Xiaoqin Hua and Nobuo Yamashita.
\newblock Iteration complexity of a block coordinate gradient descent method
  for convex optimization.
\newblock {\em {SIAM} Journal on Optimization}, 25(3):1298--1313, 2015.

\bibitem{LSZ2009}
John Langford, Alex~J. Smola, and Martin Zinkevich.
\newblock Slow learners are fast.
\newblock In {\em NIPS}, pages 2331--2339, 2009.

\bibitem{LBK2009}
Jure Leskovec, Lars Backstrom, and Jon~M. Kleinberg.
\newblock Meme-tracking and the dynamics of the news cycle.
\newblock In {\em KDD}, pages 497--506, 2009.

\bibitem{LiuW2015}
Ji~Liu and Stephen~J. Wright.
\newblock Asynchronous stochastic coordinate descent: Parallelism and
  convergence properties.
\newblock {\em {SIAM} Journal on Optimization}, 25(1):351--376, 2015.

\bibitem{LWRBS2015}
Ji~Liu, Stephen~J. Wright, Christopher R{\'{e}}, Victor Bittorf, and Srikrishna
  Sridhar.
\newblock An asynchronous parallel stochastic coordinate descent algorithm.
\newblock {\em Journal of Machine Learning Research}, 16:285--322, 2015.

\bibitem{LX2013}
Zhaosong Lu and Lin Xiao.
\newblock On the complexity analysis of randomized block-coordinate descent
  methods.
\newblock {\em CoRR}, abs/1305.4723, 2013.

\bibitem{MS2012}
Jason~R. Marden and Jeff~S. Shamma.
\newblock Revisiting log-linear learning: Asynchrony, completeness and
  payoff-based implementation.
\newblock {\em Games and Economic Behavior}, 75(2):788--808, 2012.

\bibitem{MVdGB2008}
Lukas Meier, Sara Van De~Geer, and Peter Bühlmann.
\newblock The group lasso for logistic regression.
\newblock {\em Journal of the Royal Statistical Society: Series B (Statistical
  Methodology)}, 70(1):53--71, 2008.

\bibitem{Nesterov2004}
Yu. Nesterov.
\newblock {\em Introductory Lectures on Convex Optimization: A Basic Course}.
\newblock Springer US, 2004.

\bibitem{Nesterov2012}
Yu. Nesterov.
\newblock Efficiency of coordinate descent methods on huge-scale optimization
  problems.
\newblock {\em SIAM J. Optimization}, 22(2):341--362, 2012.

\bibitem{PalLa2015}
Siddharth Pal and Richard~J. La.
\newblock Simple learning in weakly acyclic games and convergence to {N}ash
  equilibria.
\newblock 2015.
\newblock \url{http://www.ece.umd.edu/~hyongla/PAPERS/ALLERTON15.pdf}.

\bibitem{PY2010}
Christos~H. Papadimitriou and Mihalis Yannakakis.
\newblock An impossibility theorem for price-adjustment mechanisms.
\newblock {\em PNAS}, 5(107):1854--1859, 2010.

\bibitem{RT2014}
Michael~G. Rabbat and Konstantinos~I. Tsianos.
\newblock Asynchronous decentralized optimization in heterogeneous systems.
\newblock In {\em {IEEE} Conference on Decision and Control, {CDC}}, pages
  1125--1130, 2014.

\bibitem{RiT2014}
Peter Richt{\'{a}}rik and Martin Tak{\'{a}}c.
\newblock Iteration complexity of randomized block-coordinate descent methods
  for minimizing a composite function.
\newblock {\em Math. Program.}, 144(1-2):1--38, 2014.

\bibitem{ST2013}
Ankan Saha and Ambuj Tewari.
\newblock On the nonasymptotic convergence of cyclic coordinate descent
  methods.
\newblock {\em {SIAM} Journal on Optimization}, 23(1):576--601, 2013.

\bibitem{SGV1998}
Craig Saunders, Alexander Gammerman, and Volodya Vovk.
\newblock Ridge regression learning algorithm in dual variables.
\newblock In {\em Proceedings of the Fifteenth International Conference on
  Machine Learning}, ICML '98, pages 515--521, San Francisco, CA, USA, 1998.
  Morgan Kaufmann Publishers Inc.

\bibitem{Tib94}
Robert Tibshirani.
\newblock Regression shrinkage and selection via the lasso.
\newblock {\em Journal of the Royal Statistical Society, Series B},
  58:267--288, 1994.

\bibitem{Tseng2001}
Paul Tseng.
\newblock Convergence of a block coordinate descent method for
  nondifferentiable minimization.
\newblock {\em Journal of Optimization Theory and Applications},
  109(3):475--494, 2001.

\bibitem{TsengY2009}
Paul Tseng and Sangwoon Yun.
\newblock A coordinate gradient descent method for nonsmooth separable
  minimization.
\newblock {\em Math. Program.}, 117(1-2):387--423, 2009.

\bibitem{TR2012}
Konstantinos~I. Tsianos and Michael~G. Rabbat.
\newblock Distributed dual averaging for convex optimization under
  communication delays.
\newblock In {\em American Control Conference, {ACC}}, pages 1067--1072, 2012.

\bibitem{TBA1986}
John~N. Tsitsiklis, Dimitri~P. Bertsekas, and Michael Athans.
\newblock Distributed asynchronous deterministic and stochastic gradient
  optimization algorithms.
\newblock {\em IEEE Transactions on Automatic Control}, 31(9):803--812, 1986.

\bibitem{Uzawa1960}
Hirofumi Uzawa.
\newblock Walras' tatonnement in the theory of exchange.
\newblock {\em Review of Economic Studies}, 27(3):182--194, 1960.

\bibitem{Walras1874}
L\'{e}on Walras.
\newblock {\em El{\'{e}}ments d'{\,{E}}conomie Politique Pure}.
\newblock Corbaz, 1874.
\newblock (Translated as: \emph{Elements of Pure Economics.} Homewood, IL:
  Irwin, 1954.).

\bibitem{Wright2015}
Stephen~J. Wright.
\newblock Coordinate descent algorithms.
\newblock {\em Math. Program.}, 151(1):3--34, 2015.

\end{thebibliography}

\end{document}